\documentclass[a4paper,11pt,leqno]{amsart}
\usepackage{amssymb,amsmath,latexsym, comment}
\usepackage{epsfig,caption}
\usepackage{mathrsfs}
\usepackage{hyperref}
\usepackage{color,graphics,subfigure}

\usepackage{epic,eepic,wrapfig,ifthen}
\usepackage{url}
\newtheorem{thm}{Theorem}[section]
\newtheorem{cor}[thm]{Corollary}
\newtheorem{lemma}[thm]{Lemma}
\newtheorem{prop}[thm]{Proposition}

\newtheorem{remark}[thm]{Remark}

\numberwithin{equation}{section}

\newcommand{\formula}[2][nolabel]
{\ifthenelse{\equal{#1}{nolabel}}
 {\begin{align*} #2 \end{align*}}
 {\ifthenelse{\equal{#1}{}}
  {\begin{align} #2 \end{align}}
  {\begin{align} \label{#1} #2 \end{align}}
 }
}

\def\qed{{\hfill $\Box$ \bigskip}}

\DeclareMathOperator*{\esssup}{ess\,sup}
\DeclareMathOperator*{\essinf}{ess\,inf}

\renewcommand{\bar}{\overline}
\newcommand{\cal}[1]{\mathcal{#1}}
 \def\sB {{\cal B}} \def\sC {{\cal C}}
\def\sD {{\cal D}} \def\sE {{\cal E}} \def\sF {{\cal F}}
  \def\sI {{\cal I}}
  
 \def\sN {{\cal N}} 
 \def\sQ {{\cal Q}}

\def\bP {{\mathbb P}}  \def\bR {{\mathbb R}}

\def\R {{\mathbb R}}
\def\N {{\mathbb N}}

\def\P{{\mathbb P}}
\def\E{{\mathbb E}}
\def\L{{\mathbb L}}
\def\e{{\mathbf e}}

\def\1{{\bf 1}}

\def\nn{\nonumber}

\def\wt{\widetilde}
\def\wh{\widehat}

\def\eps{{\varepsilon}}

\def\oR {{\overline{ \mathbb R}}}
\def\uB {{B_+}}

\def\rF {{\bar\sF}}
\def\rY {\bar Y}
\def\rP {\bar P}
\allowdisplaybreaks

\makeatletter
\@addtoreset{equation}{section}

\makeatother

\begin{document}
\bibliographystyle{plain}

\title[Heat kernel estimates for Dirichlet forms degenerate at the boundary]
{ \bf    Heat kernel estimates for Dirichlet forms degenerate at the boundary}

\author{Soobin Cho, \quad Panki Kim, \quad Renming Song \quad and \quad Zoran Vondra\v{c}ek}

\address[Cho]{Department of Mathematics, University of Illinois at Urbana-Champaign, Urbana, IL 61801, USA}
\curraddr{}
\email{soobinc@illinois.edu}

\address[Kim]{Department of Mathematical Sciences and Research Institute of Mathematics,
	Seoul National University,	Seoul 08826, Republic of Korea}\thanks{This research is  supported by the National Research Foundation of Korea(NRF) grant funded by the Korea government(MSIP) (No. RS-2023-00270314).
}
\curraddr{}
\email{pkim@snu.ac.kr}

\address[Song]{
	Department of Mathematics, University of Illinois, Urbana, IL 61801,
	USA}
\curraddr{}\thanks{Research supported in part by a grant from
	the Simons Foundation (\#960480, Renming Song)}
\email{rsong@illinois.edu}

\address[Vondra\v{c}ek]
{
	Department of Mathematics, Faculty of Science, University of Zagreb, Zagreb, Croatia, 
	and Dr.~Franjo Tu{\dj}man Defense and Security University, Zagreb, Croatia
}
\curraddr{}\thanks{ Research supported in part by the Croatian Science Foundation under the projects IP-2018-01-4197 and IP-2022-10-2277. 
	(Zoran Vondra\v{c}ek)}
\email{vondra@math.hr}

\date{}

\begin{abstract}

The goal of this paper is to establish sharp two-sided estimates on the heat kernels of two types of purely discontinuous symmetric Markov processes in  
the upper half-space of $\R^d$  
with jump kernels degenerate at the boundary.
The jump kernels are of the form $J(x,y)=\sB(x,y)|x-y|^{-\alpha-d}$, $\alpha\in (0,2)$, where the function $\sB$ depends on four parameters and may vanish at the boundary. Our results are the first sharp two-sided estimates for the heat kernels of non-local operators with jump kernels degenerate at the boundary.

The first type of processes
are conservative Markov processes
on $\overline{\R}^d_+$ 
with jump kernel $J(x,y)$. Depending on the regions where 
the parameters belong, the heat kernels estimates have three different forms, two of them are qualitatively different from all 
previously known heat kernel estimates.

The second type of processes 
are the processes above killed either by a critical potential or upon hitting the boundary of the half-space. We establish that their heat kernel estimates have the approximate factorization property with survival probabilities decaying as a power of the distance to the boundary, where the power depends on the constant in the critical potential.
\end{abstract}

\maketitle

\bigskip



\bigskip
\tableofcontents

\section{Introduction and main results}

In this paper, we study both conservative and non-conservative purely discontinuous (self-similar) Markov processes in 
the upper half-space of $\R^d$
with jump kernels of the form $J(x,y)=|x-y|^{-d-\alpha}\sB(x,y)$,  
$\alpha\in (0,2)$. The function $\sB(x,y)$ may tend to 0 when $x$ or $y$ tends to 
the boundary of the half-space, and so the jump kernel may be degenerate.
Our main  focus is to establish sharp two-sided estimates of the transition densities of the processes  (or the heat kernels of the corresponding non-local operators). 
Heat kernel estimates for non-local operators have been the subject of many papers in the last twenty years, 
see 
\cite{BKKL19, BGR10, BGR14, CKK11, CKS-jems, CKS14, CK03, CK08,  CKW20a, CKW21, GHH18, GHH} 
and the references therein.
In all of the papers mentioned above,
the function $\sB$ is assumed to be bounded between two positive constants, which can be viewed as a uniform ellipticity condition for non-local operators. To the best of our knowledge, the current paper is the first one to study 
sharp two-sided heat kernel estimates when the jump kernel is degenerate.
Boundary Harnack principle and sharp two-sided Green function estimates for 
purely discontinuous Markov processes with degenerate jump kernels, which can be viewed as the elliptic counterpart of results of the current paper, have been recently obtained in \cite{KSV-jump, KSV-Green, KSV-nokill}.

The first type of processes we look at are conservative jump processes on 
$\overline{\R}^d_+:=\{x=(\wt{x}, x_d): \wt{x}\in \R^{d-1},  x_d\ge 0\}$ with jump kernel 
$J(x,y)=|x-y|^{-d-\alpha}\sB(x,y)$, where the function $\sB(x,y)$ 
is symmetric, homogeneous, horizontally translation invariant, and is allowed to approach 0 at the boundary at arbitrary fixed polynomial rate in terms of 
some non-negative parameters $\beta_1, \beta_2, \beta_3, \beta_4$. The generator of such a process is the non-local operator
$$
L^{\sB}_\alpha f(x)=\mathrm{p.v.} \int_{{\R}^d_+}(f(y)-f(x))|x-y|^{-d-\alpha}\sB(x,y)\, dy.
$$
When $\sB(x,y)$ is bounded between two positive constants, 
the heat kernel estimates for such 
processes  are of the form $\min\{t^{-d/\alpha}, tJ(x,y)\}$. This was first established in the pioneering work \cite{CK03}, even for the case of metric measure spaces. 
This form of the  estimates reflects the fact that the main contribution to the heat kernel comes from one (big) jump from $x$ to $y$. 
This feature has been observed in all subsequent studies, see 
\cite{BKKL19, CK08, CKK11, CKW20a, CKW21, GHH18} and the references therein.
In our setting of 
jump kernel degenerate at the boundary, 
there are two novel features in the heat kernel estimates. The first one appears 
in the cases when 
$d=1$ or 
the involved parameters satisfy $\beta_2<\alpha+\beta_1$, in which case the heat kernel is comparable to $\min\{t^{-d/\alpha}, tJ(x+t^{1/\alpha}\mathbf{e}_d, y+t^{1/\alpha}\mathbf{e}_d)\}$, where $\mathbf{e}_d=(\wt{0}, 1)$. In words, the form of the heat kernel estimates
shows that the main contribution to the heat kernel at time $t$  comes from 
one jump from the point $t^{1/\alpha}$ units above $x$ to the point $t^{1/\alpha}$ units above $y$. 
Due to the fact the jump kernel vanishes at the boundary, it is very unlikely that the process will make one (big) jump from (or to) a point very close to the boundary. 
The second feature is more striking and indicates a sort of a phase-transition 
at the level $\beta_2=\alpha+\beta_1$: 
When $d\ge 2$ and 
the parameters satisfy $\beta_2\ge \alpha+\beta_1$, 
in addition to the already mentioned part $\min\{t^{-d/\alpha}, tJ(x+t^{1/\alpha}\mathbf{e}_d, y+t^{1/\alpha}\mathbf{e}_d)\}$, the sharp heat kernel estimates include a part 
which reflects a significant contribution to the heat kernel coming from two jumps connecting $x$ and $y$. 
The precise description is given in Theorem \ref{t:HKE}.

The second type of processes we look at are the ones described in the previous paragraph but killed either by a critical potential 
$\kappa x_d^{-\alpha}$ or upon hitting the boundary of 
$\R^d_+:=\{x=(\wt{x}, x_d): \wt{x}\in \R^{d-1},  x_d > 0\}$ (the latter happens only when $\alpha\in (1,2)$). 
The generator of such a process is the non-local operator
$L^{\sB} f(x)=L^{\sB}_\alpha f(x)- \kappa x_d^{-\alpha}f(x)$. 
We study the effect of such  killings on the heat kernel. 
When $\sB(x,y)$ is bounded between two positive constants,
the effect of killing on the heat kernel of the process 
is, after intensive research during the last fifteen years, fairly well understood.
In most cases, the heat kernel of the killed process has the so called approximate factorization property:
It is comparable to the product of the heat kernel of non-killed (original) process and  
the survival probabilities starting from points $x$ and $y$. 
In case of smooth open sets, the survival probability can be expressed in terms of the distance between the point and the boundary, see 
\cite{BGR10, BGR14, CK, CKS-jems, CKS14, CKSV, GKK} and the references therein.
In our setting of jump kernel degenerate at the boundary, we establish the same property: The heat kernel of the killed process enjoys the approximate factorization property with 
survival probabilities decaying as the $q$-th power 
of the distance to the boundary, where $q$ is in one-to-one correspondence with the constant $\kappa \ge 0$, see Theorem \ref{t:HKE-kappa}. 
When $\kappa=0$, this theorem generalizes \cite{CKS-ptrf} (for half spaces) where the factorization property for censored stable process is established. Due to the quite complicated form of the heat kernel estimates in Theorem \ref{t:HKE}, obtaining the factorization property in Theorem \ref{t:HKE-kappa} is a formidable task.

We now introduce the precise setup and state 
the main results of this paper.
This setup was first introduced in \cite{KSV-jump} and was  motivated by  the results of \cite{KSV1,  KSV2} on subordinate killed L\'evy processes.  
In fact, subordinate killed L\'evy processes, whose analytical counterparts are fractional powers of Dirichlet fractional Laplacians, 
are the main natural examples of Markov processes with jump kernels satisfying the assumptions below.

Let $d \ge 1$ and  $0<\alpha<2$. Recall that
$\R^d_+=\{(\wt x, x_d) \in \R^d: x_d>0\}$ and $\overline\R^d_+=\{(\wt x, x_d) \in \R^d: x_d\ge0\}$. 
Here and below, 
$\L^{b}(a):=\log^{b}(e+a)$,
 $a\wedge b:=\min \{a, b\}$, $a\vee b:=\max\{a, b\}$, and $a\asymp b$ means that $c\le b/a \le c^{-1}$ for some $c\in (0,1)$.
We will write $\L^{1}(a)$ as $\L(a)$.
For $b_1,b_2,b_3,b_4 \ge 0$, let
\begin{align}\label{def:wtB}
	&B_{b_1, b_2, b_3, b_4}(x,y)\\
	&:= \Big(\frac{x_d \wedge y_d}{|x-y|} \wedge 1  \Big)^{b_1} \Big(\frac{x_d \vee y_d}{|x-y|} \wedge 1 \Big)^{b_2} \L^{b_3} \Big(\frac{(x_d \vee y_d) \wedge |x-y|}{x_d \wedge y_d \wedge |x-y|}  \Big) \, \L^{b_4} \Big(\frac{ |x-y|}{(x_d \vee y_d) \wedge |x-y|} \Big).	\nn
\end{align}

\begin{remark}\label{r:d=1}
	{\rm When $d=1$, we have $x \vee y = x\wedge y + |x-y| \ge |x-y|$ for all $x,y \in [0,\infty)$ so that the second and the fourth terms in \eqref{def:wtB} are positive constants. Thus, the parameters $b_2$ and $b_4$ in \eqref{def:wtB} are irrelevant if $d=1$.}
\end{remark}

Let $J(x,y)=|x-y|^{-d-\alpha} \sB(x,y)$, where  
 the function $\sB:\overline{\R}^d_+\times \overline{\R}^d_+\to [0,\infty)$ 
will be
assumed to satisfy some or all of the following hypotheses:

\medskip

\setlength{\leftskip}{4mm}

\noindent {\bf (A1)} $\sB(x,y)=\sB(y,x)$ for all $x,y \in \overline{\R}^d_+$.

\smallskip

\noindent {\bf (A2)} If $\alpha \ge 1$, then there exist $\theta>\alpha-1$ and 
$C_1>0$ such that
\begin{equation*}
	|\sB(x,x)-\sB(x,y)| \le 
	C_1\left(\frac{|x-y|}{x_d \wedge y_d} \right)^\theta, \quad x,y \in \R^d_+.
\end{equation*}

\noindent {\bf (A3)}{(I)} 
There exist $C_2\ge1$ 
and  $\beta_1,\beta_2,\beta_3,\beta_4 \ge 0$, with $\beta_1>0$ if $\beta_3>0$, and $\beta_2>0$ if $\beta_4>0$, such that
\begin{equation*}
	C_2^{-1} B_{\beta_1, \beta_2, \beta_3, \beta_4}(x,y) \le \sB(x,y) \le C_2, \quad x,y \in \overline{\R}^d_+.
\end{equation*}

\hspace{2.2mm} {(II)} 
There exists $C_3>0$ such that
\begin{equation*}
	\sB(x,y) \le C_3B_{\beta_1, \beta_2, \beta_3, \beta_4}(x,y), \quad x,y \in \overline{\R}^d_+,
\end{equation*}
where $\beta_1, \beta_2, \beta_3, \beta_4$ are the same constants as in {(I)}.

\noindent {\bf (A4)}  For all $x,y \in \overline{\R}^d_+$ and $a>0$, $\sB(ax,ay)=\sB(x,y)$.  In case $d \ge 2$, for all $x,y \in \overline{\R}^d_+$ and $\wt z\in \R^{d-1}$,  $\sB(x+(\wt z, 0), y + (\wt z,0)) = \sB(x,y)$.

\setlength{\leftskip}{0mm}

\medskip

These four hypotheses were introduced in \cite{KSV-jump}, and with the same notation as above repeated in \cite{KSV-Green, KSV-nokill}. Condition \textbf{(A2)} is not needed in Theorems \ref{t:HKE} and \ref{t:green-1}, while in Theorems \ref{t:HKE-kappa} and \ref{t:green-2} it is used 
through several results from \cite{KSV-jump, KSV-Green, KSV-nokill}.

\textit{Throughout this paper, we always assume that  $\sB(x,y)$ satisfies {\bf (A1)}, {\bf (A3)}{\rm (I)} and {\bf (A4)}.} 

\smallskip

Consider the following symmetric form
\begin{equation*}
	\sE^0(u,v):= \frac{1}{2}\int_{\R^d_+} \int_{\R^d_+}  (u(x)-u(y))(v(x)-v(y))
	J(x, y)dydx.
\end{equation*} 
Since $\sB(x,y)$ is bounded,  $C_c^\infty(\oR^d_+)$ is closable in $L^2(\oR^d_+,dx) = L^2(\R^d_+,dx)$  by Fatou's lemma. Let  $\rF$ be the closure of  $C_c^\infty(\oR^d_+)$ in $L^2(\R^d_+,dx)$ under the norm  
$(\sE^0_1)^{1/2}$ where $\sE^0_1:=\sE^0+(\cdot,\cdot)_{L^2(\R^d_+,dx)}$. 
 Then  $(\sE^0, \rF)$ is a regular Dirichlet form on $L^2(\R^d_+,dx)$. Let  $\rY=(\rY_t,t \ge 0; \P_x,x \in \overline\R^d_+\setminus \sN')$ be the Hunt 
process associated with  $(\sE^0, \rF)$,  where $\sN'$ is an  exceptional set.

Here is our first main result. The heat kernel estimates are expressed in different but equivalent forms,
each providing a different viewpoint. 
 Recall that  $\e_d=(\wt 0, 1) \in \R^d$.

\begin{thm}\label{t:HKE}
	Suppose that {\bf (A1)}, {\bf (A3)} and {\bf (A4)}  hold. Then the process $\rY$ can be refined to be a
	conservative   Feller process with strong Feller property starting from every point in $\oR^d_+$ and  has a jointly continuous heat kernel $\bar p:(0,\infty) \times \oR^d_+ \times \oR^d_+ \to (0,\infty)$.
 Moreover, the heat kernel $\bar p$ has the following estimates: 
\noindent
(a) When $d=1$, we have that  for all $(t,x,y) \in (0,\infty) \times \oR_+ \times \oR_+$,
\begin{align}\label{e:allcase1}	
		&\bar p(t,x,y) \asymp 	 	 t^{-1/\alpha} \wedge tJ(x+t^{1/\alpha},y+t^{1/\alpha} ).
\end{align}

\noindent
(b) When $d \ge 2$,  we have that for all $(t,x,y) \in (0,\infty) \times \oR^d_+ \times \oR^d_+$,
	\begin{align}\label{e:allcase}	
		&\bar p(t,x,y) \asymp 	 	 t^{-d/\alpha} \wedge
		\Big[tJ(x+t^{1/\alpha}\e_d,y+t^{1/\alpha}\e_d )\\
		&\,+{\bf 1}_{\{\beta_2 \ge \alpha+\beta_1\}} t^2 
		 \int_{(x_d \vee y_d\vee t^{1/\alpha}) \wedge (|x-y|/4)}^{|x-y|/2} 
		    \! \! \! \! \!   
		J(x+t^{1/\alpha}\e_d, 
		 x+r \e_d  ) \, J(x+ r \e_d  ,y+t^{1/\alpha}\e_d)r^{d-1} dr\Big].\nn
\end{align}
	Furthermore the heat kernel 	estimates in \eqref{e:allcase} can be rewritten in terms of   
	the function $B$ 
explicitly by considering three cases separately: 

\smallskip

	\noindent (i) If $\beta_2<\alpha+\beta_1$, then for all $(t,x,y) \in (0,\infty) \times \oR^d_+ \times \oR^d_+$,
	\begin{align}\label{e:Case1}
			\bar p(t,x,y) &
			\asymp \Big( t^{-d/\alpha} \wedge \frac{t 	B_{\beta_1, \beta_2, \beta_3, \beta_4}(x+t^{1/\alpha}\e_d,y+t^{1/\alpha}\e_d)	}{|x-y|^{d+\alpha}}\Big)
				\\
		&  
		\asymp  \Big( t^{-d/\alpha} \wedge \frac{t}{|x-y|^{d+\alpha}}\Big)	B_{\beta_1, \beta_2, \beta_3, \beta_4}(x+t^{1/\alpha}\e_d,y+t^{1/\alpha}\e_d). \nn
	\end{align}

	\noindent (ii) If $\beta_2>\alpha+\beta_1$, then for all $(t,x,y) \in (0,\infty) \times \oR^d_+ \times \oR^d_+$,
	\begin{align}\label{e:Case2}
		&\bar p(t,x,y) 	\asymp	\Big( t^{-d/\alpha} \wedge \frac{t}{|x-y|^{d+\alpha}}\Big)	\Big[	B_{\beta_1, \beta_2, \beta_3, \beta_4}(x+t^{1/\alpha}\e_d,y+t^{1/\alpha}\e_d)\\
		&+ \Big(1 \wedge  \frac{t}{|x-y|^\alpha}\Big)B_{\beta_1, \beta_1, 0, \beta_3}(x+t^{1/\alpha}\e_d,y+t^{1/\alpha}\e_d)  \L^{\beta_3} \Big( \frac{|x-y|}{((x_d \wedge y_d)+ t^{1/\alpha}) \wedge |x-y|} \Big) \Big] .	\nn
		\end{align}

	\noindent (iii) If $\beta_2=\alpha+\beta_1$, then for all $(t,x,y) \in (0,\infty) \times \oR^d_+ \times \oR^d_+$,
	\begin{align}\label{e:Case3}	
		&\bar p(t,x,y) \asymp 		\Big( t^{-d/\alpha} \wedge \frac{t}{|x-y|^{d+\alpha}}\Big) \Big[ B_{\beta_1, \beta_2, \beta_3, \beta_4}(x+t^{1/\alpha}\e_d,y+t^{1/\alpha}\e_d) \\	
			& \!\!+ \!  \Big( 1 \wedge  \frac{t}{|x-y|^\alpha} \Big) B_{\beta_1, \beta_1, 0, \beta_3+\beta_4+1}(x\!+\!t^{1/\alpha}\e_d,y\!+\!t^{1/\alpha}\e_d) \L^{\beta_3} 
			\Big( \frac{|x-y|}{((x_d \wedge y_d)+ t^{1/\alpha}) \wedge |x-y|}\Big)\Big].  \nn
		\end{align}
\end{thm}

Note that if $x$ or $y$  
is close to the boundary, then $J(x+t^{1/\alpha}\e_d,y+t^{1/\alpha}\e_d )$ is not comparable to $J(x,y)$ in our setting. 
Thus,  
even in the case $d=1$ or $\beta_2<\alpha+\beta_1$, 
the form of the heat kernel estimates is different from the usual form. 
The appearance of the two terms in the brackets on the right-hand sides of \eqref{e:Case2}--\eqref{e:Case3}
reflects the fact that the dominant contribution to the heat kernel may come from either one jump or two jumps. Moreover, it
is  easy to see that neither of these two terms dominates the other one for all $(t, x, y)$. Below we illustrate this feature 
when $d\ge 2$ and $\beta_2>\alpha+\beta_1$.

\begin{figure}[b!]
	\centering
	\subfigure[One jump regime]{\includegraphics[width=
		9.2cm]{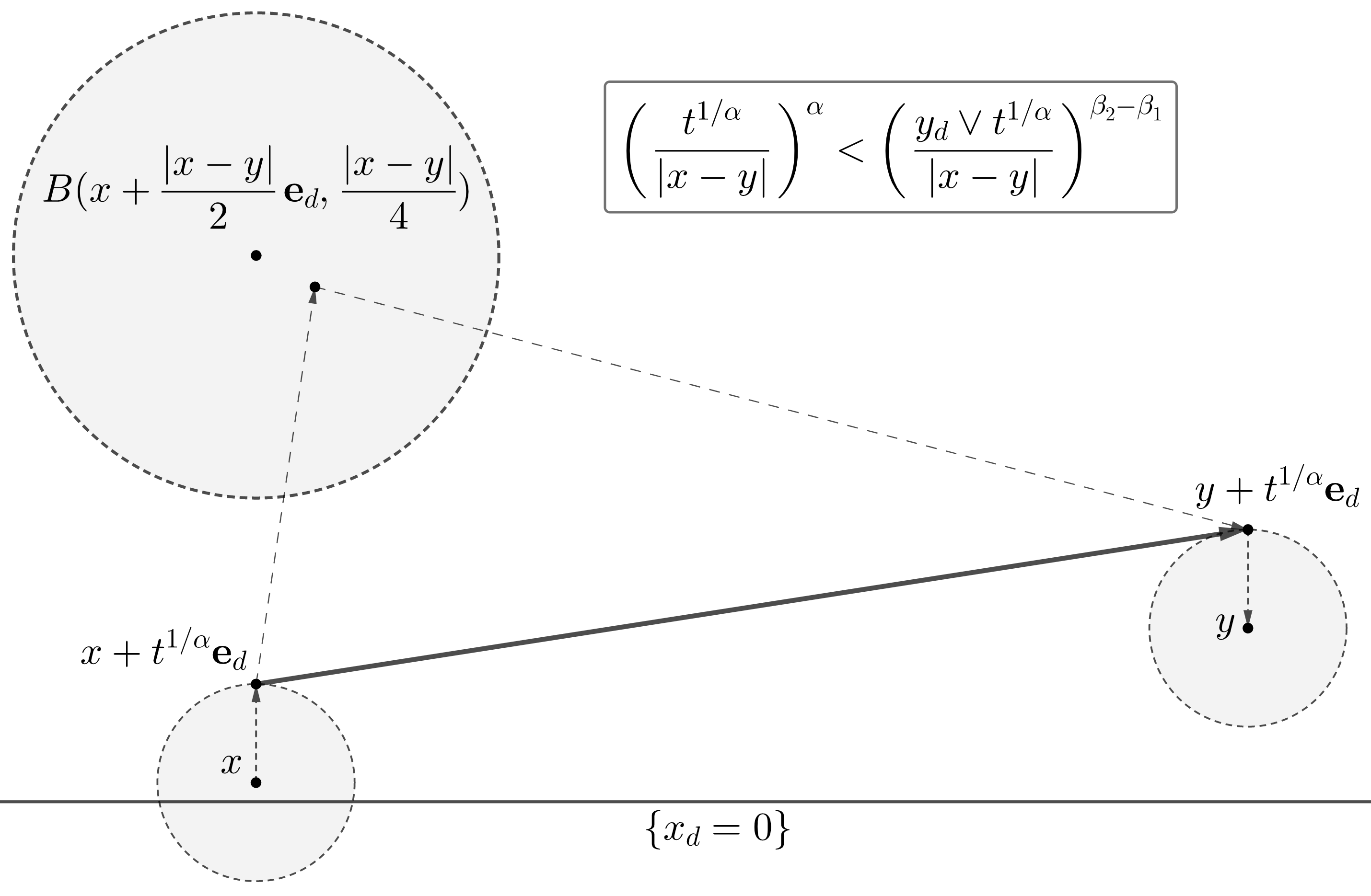}
	} 
	
	\vspace{0.1cm}
	
	\subfigure[Two jumps regime]{\includegraphics[width=
		9.2cm]{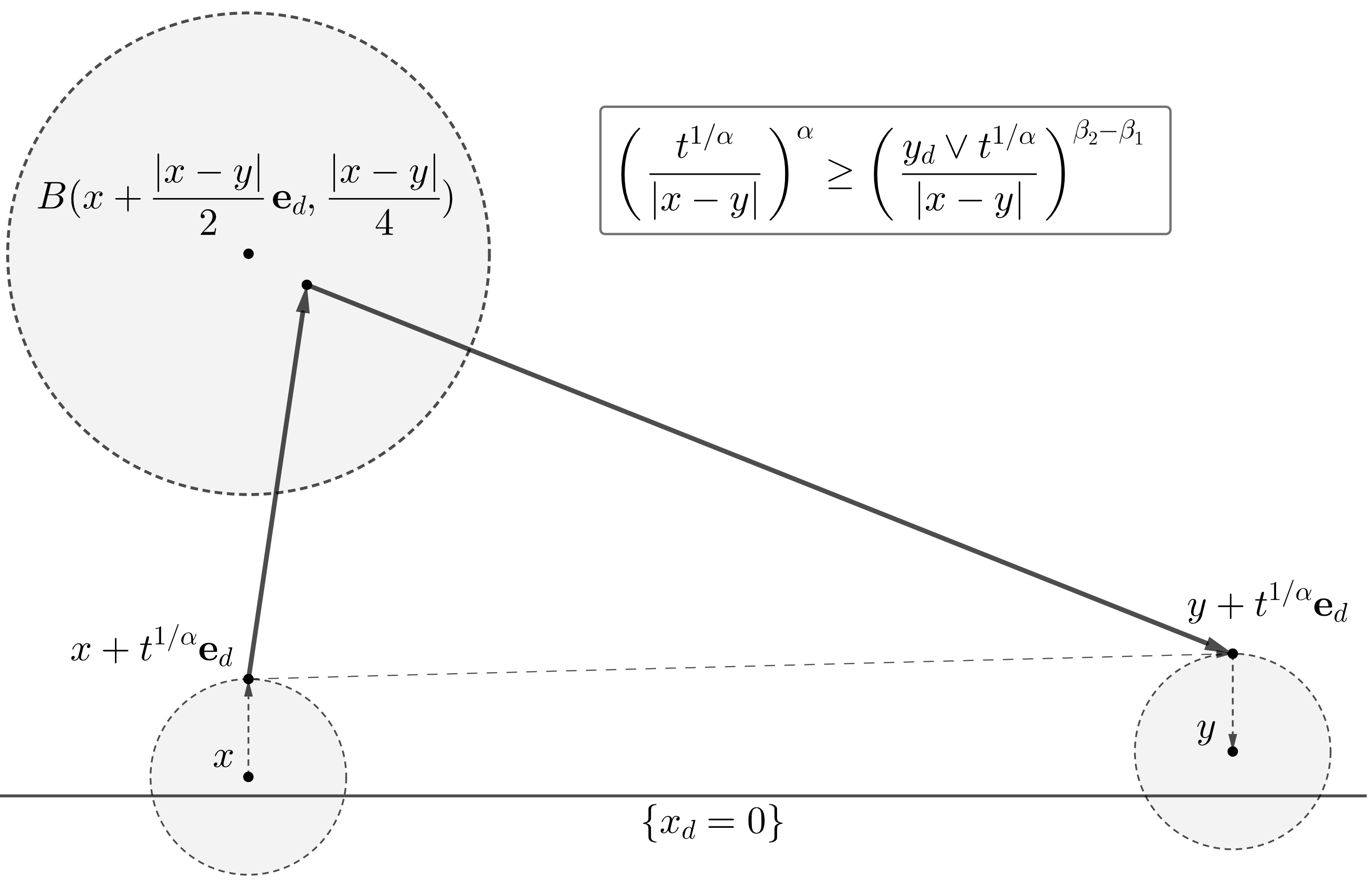}
	}
	\vspace{-0.2cm}
	\caption{\small Dominant path from 
		$x$ to $y$  (with $x_d\le y_d$)  at time 
		$t<(|x-y|/6)^\alpha$ 
		when $\beta_2>\alpha+\beta_1$, $\beta_3=\beta_4=0$.}\label{fig:1}
	\vspace{-0.4cm}	
\end{figure}

Let $d\ge 2$ and $ \beta_2>\alpha+\beta_1$. 
When $|x-y|>6t^{1/\alpha}$, the second term in the brackets in \eqref{e:Case2}, i.e., 
\begin{align}
\Big(1 \wedge  \frac{t}{|x-y|^\alpha}\Big)B_{\beta_1, \beta_1, 0, \beta_3}(x+t^{1/\alpha}\e_d,y+t^{1/\alpha}\e_d) \L^{\beta_3} \Big( \frac{|x-y|}{((x_d \wedge y_d)+ t^{1/\alpha}) \wedge |x-y|} \Big)
\label{e:Case21a}
\end{align}
is comparable to 
\begin{align}
	t |x-y|^{d+\alpha}  \int_{B(x+2^{-1}|x-y|\e_d, \,
	  	4^{-1}|x-y|)}  	J(x+t^{1/\alpha}\e_d,z)\, J(z,y+t^{1/\alpha}\e_d)dz ,\label{e:Case21}
	 	\end{align}
see Remark \ref{r:proofCase21} below.
In the special case when $\beta_3=\beta_4=0$ (and $ \beta_2>\alpha+\beta_1$), Figure \ref{fig:1},  
where the comparability of \eqref{e:Case21a} and \eqref{e:Case21} is used, 
illustrates the regions where one jump and two jumps dominate. 

The different forms of the heat kernel estimates in Theorem \ref{t:HKE}  are consequences of the relationship among $\beta_1$, $\beta_2$ and $\alpha$ only, the values of $\beta_3$ and $\beta_4$
do not play any role in this.
To reduce technicalities, the reader, on a first reading,  may assume  $\beta_3=\beta_4=0$,  without loosing
the essential features of this paper. Note that, even in this special case, a logarithmic term appears in the estimates when $\beta_2=\alpha+\beta_1$.  We need the logarithmic terms in the assumptions to cover the important examples in \cite{CKSV, KSV2}.

\smallskip

 Let $\sF^0$ be the closure of $C_c^\infty(\R^d_+)$ in $L^2(\R^d_+,dx)$ under the norm  
 $(\sE^0_1)^{1/2}$. 
  Then $(\sE^0,\sF^0)$ is also a regular Dirichlet form on $L^2(\R^d_+,dx)$ and it is the part form of $(\sE^0, \rF)$ on $\R^d_+$.
 Let $Y^0=(Y^0_t,t \ge 0; \P_x,x \in \R^d_+\setminus \sN_0)$ be the Hunt 
 process associated with $(\sE^0, \sF^0)$  where $\sN_0$ is an exceptional set.
 $Y^0$ can be regarded as the part process of $\rY$ killed upon exiting $\R^d_+$.

For $\kappa \in (0,\infty)$, define 
\begin{align*}
	\sE^\kappa(u,v)&:=\sE^0(u,v) + 
	\int_{\R^d_+} u(x)v(x)\kappa x_d^{-\alpha}dx, \\
	\sF^\kappa&:= 	\wt \sF^0 \cap L^2(\R^d_+, 
	\kappa x_d^{-\alpha}dx),\nn
\end{align*}
where $\wt \sF^0$ is the family of all $\sE^0_1$-quasi-continuous functions in $\sF^0$. Then $(\sE^\kappa,\sF^\kappa)$ is a regular Dirichlet form on $L^2(\R^d_+,dx)$ with $C_c^\infty(\R^d_+)$ as a special standard core,  see \cite[Theorems 6.1.1 and 6.1.2]{FOT}.
 Let  $Y^\kappa=(Y^\kappa_t,t \ge 0; \P_x,x \in \R^d_+\setminus \sN_\kappa)$ be the Hunt process associated with 
$(\sE^\kappa, \sF^\kappa)$ where $\sN_\kappa$ is an exceptional set. For $\kappa \in [0,\infty)$, we denote by $\zeta^\kappa$ the lifetime of $Y^\kappa$. Define $Y^\kappa_t=\partial$ for $t \ge \zeta^\kappa$, where $\partial$ is a cemetery point added to the state space $\R^d_+$. 

We now associate with the constant $\kappa$ a positive parameter $q$ which plays an important  role in the paper.  
For $q \in (-1-\beta_1, \alpha +\beta_1)$, 
we define $C(\alpha,q,\sB)$ to be
\begin{align*}
	\begin{cases}
		\displaystyle		\int_{\R^{d-1}}  \frac{1}{(|\wt u|^2+1)^{(d+\alpha)/2}}   \int_0^1 \frac{(s^q-1)(1-s^{\alpha-q-1})}{(1-s)^{1+\alpha}} 
	 \sB\big( ( (1-s)\wt u,1) ,s \e_d\big)  ds \,d\wt u 
	 &\mbox{if } d \ge 2,\\[0.5cm]
		\displaystyle \int_0^1 \frac{(s^q-1)(1-s^{\alpha-q-1})}{(1-s)^{1+\alpha}} \sB\big( 1 ,s \big)ds, 
	&\mbox{if } d =1.
	\end{cases}
\end{align*}
If we additionally assume that  {\bf (A3)}(II) holds, then the constant $C(\alpha,q,\sB)$ is 
well defined and finite 
for every 
$q \in (-1-\beta_1,\alpha +\beta_1)$,  
$C(\alpha,q,\sB)=0$ if and only if $q\in\{0,\alpha-1\}$,  and 
$\lim_{q\to-1-\beta_1} C(\alpha,q,\sB)= \lim_{q \to \alpha+\beta_1} C(\alpha,q,\sB)=\infty$
(see \cite[Lemma 5.4 and Remark 5.5]{KSV-jump}). Note that for every $s \in (0,1)$, $q \mapsto (s^q-1)(1-s^{\alpha-q-1})$ is strictly decreasing on 
$(-1-\beta_1, (\alpha-1)/2)$ 
and strictly increasing on $((\alpha-1)/2,\alpha+\beta_1)$. Thus, the shape of the map $q\mapsto C(\alpha,q,\sB)$ is given as follows.
\begin{table}[hbt!]
	\begin{center}
		\vspace{-0.2cm}
		\begin{tabular}{|c||c|c|c|c|c|c|c|c|c|} 
			\hline  
			\rule{0in}{2.5ex}$q$ & $ -1-\beta_1  $ & $\cdots$  &$(\alpha-1) \wedge 0$   & $\cdots$ & $\frac{1}{2}(\alpha-1)$ & $\cdots$ & $(\alpha-1) \vee 0$  &$\cdots$& $\alpha+\beta_1$\\[0.5ex]
			\hline 
			\rule{0in}{2.5ex}$C(\alpha,q, \sB)$& $\infty$ & $\searrow$ &0 &$\searrow$&$\text{minimum}
			 \le 0  $&$\nearrow$&0&$\nearrow$&$\infty$\\[3pt]
			\hline 
		\end{tabular}
	\vspace{-0.3cm}
	\end{center}
\end{table}\\
Consequently, for every $\kappa\ge0$,  there exists a unique 
$q_\kappa \in [(\alpha-1)_+, \alpha+\beta_1)$
such that
\begin{equation}\label{e:killing-potential}
	\kappa=C(\alpha, q_\kappa, \sB).
\end{equation}
\begin{thm}\label{t:HKE-kappa}
Suppose that 	{\bf (A1)}--{\bf (A4)} and \eqref{e:killing-potential} hold with $q_\kappa \in [(\alpha-1)_+, \alpha+\beta_1)$. Then the process  $Y^\kappa$  
can be refined to start from every point in $\R^d_+$ and 
	has a jointly continuous heat kernel $p^\kappa:(0,\infty) \times \R^d_+ \times \R^d_+ \to (0,\infty)$. Moreover, the following approximate factorization  holds  for all $(t,x,y) \in (0,\infty) \times \R^d_+ \times \R^d_+$:
	\begin{align}\label{e:factor}
		p^\kappa(t,x,y) \asymp  \Big(1 \wedge \frac{x_d}{t^{1/\alpha}}\Big)^{q_\kappa}\Big(1 \wedge \frac{ y_d}{t^{1/\alpha}}\Big)^{q_\kappa} \bar p(t,x,y)  \asymp \P_x(\zeta^\kappa>t) \P_y(\zeta^\kappa>t)\, \bar p(t,x,y),
	\end{align}
where $\bar p(t,x,y)$ is the heat kernel of $\rY$.
\end{thm}

As a consequence of Lemma \ref{l:alpha<1} below, 
we will see that, when $\alpha\le 1$, the process $\rY$ started from $\R^d_+$ will never hit $\partial \R^d_+$ and is equal to $Y^0$. Thus for $x,y\in \R^d_+$ we have that $p^0(t,x,y)=\bar p(t,x,y)$, implying that the non-trivial content of 
Theorem \ref{t:HKE-kappa} is for $\kappa>0$ when $\alpha\le 1$.

\medskip

Let 
\begin{equation}\label{e:green-def}
	\bar G(x,y)=\int_0^\infty \bar p(t,x,y)  dt 
	\quad \text{and}\quad
	G^\kappa(x,y)=\int_0^\infty p^\kappa(t,x,y) dt.
\end{equation} 
When $\overline{G} (\cdot, \cdot)$ is not identically infinite, it is called the Green function of $\rY$, 
and when $G^\kappa(\cdot,\cdot)$ is not identically infinite, it is called the Green function of $Y^\kappa$.

As a consequence of 
the heat kernel estimates, we get the  Green function estimates. 
The following theorem says  that 
the Green function of  $\rY$ is comparable to that of the isotropic $\alpha$-stable process in $\R^d$ even
though the jump kernel of $\rY$ may be degenerate.

\begin{thm}\label{t:green-1}
	Suppose that {\bf (A1)}, {\bf (A3)}{\rm (I)} and {\bf (A4)} hold. If $d >\alpha$, then
	\begin{equation}\label{e:green-reflect}
		\bar G(x,y) \asymp \frac{1}{|x-y|^{d-\alpha}}, \quad x, y \in\oR_+^d.
	\end{equation} 
	If $d \le \alpha$, then $\bar G(x,y)=\infty$ for all $x,y \in \oR_+^d$.
\end{thm}

When $d>(\alpha+\beta_1+\beta_2) \wedge 2$, sharp two-sided estimates on $G^\kappa(x,y)$ were obtained in \cite{KSV-Green, KSV-nokill}. In the following theorem, we extend those results by removing the restriction on  $d$ and give another proof 
using the heat kernel estimates. 
The advantage of the new proof is that it explains the reason for the phase transition in the 
Green  function estimates for $d \ge 2$.

Define  $H_q(x,y)$ by
	\begin{align*}
		\begin{cases}
			\displaystyle 1 &\mbox{if } q <\alpha + \frac12 (\beta_1+\beta_2),\\[4pt]
			\displaystyle		\L^{\beta_4 +1} \Big(\frac{|x-y|}{(x_d \vee y_d)\wedge |x-y|}\Big) &\mbox{if }
			q =\alpha + \frac12 (\beta_1+\beta_2)
			,\\[7pt]
			\displaystyle  
			\Big(\frac{x_d \vee y_d}{|x-y|} \wedge 1\Big)^{2\alpha+\beta_1+\beta_2-2q}
			\L^{\beta_4} \Big(\frac{|x-y|}{(x_d \vee y_d)\wedge |x-y|}\Big)
			&\mbox{if } 
			q >\alpha + \frac12 (\beta_1+\beta_2)
			.
		\end{cases}
	\end{align*}
\begin{thm}\label{t:green-2}
	Suppose that {\bf (A1)}--{\bf (A4)} and \eqref{e:killing-potential} hold with $q_\kappa \in [(\alpha-1)_+, \alpha+\beta_1)$. When $\alpha \le 1$, suppose also that $q_\kappa>0$ (or, equivalently, $\kappa>0$).  Then $G^\kappa$ has the following estimates:
	
	\smallskip
	
\noindent	(a)  If $d \ge 2$, then for all $x,y \in \R^d_+$,
	\begin{align*}
		G^\kappa(x,y) & \asymp \frac{H_{q_\kappa}(x,y)}{|x-y|^{d-\alpha}} 	\Big(\frac{x_d \wedge y_d}{|x-y|} \wedge 1\Big)^{q_\kappa}\Big(\frac{x_d \vee y_d}{|x-y|} \wedge 1\Big)^{q_\kappa}.
	\end{align*}

\noindent (b)   If $d=1$, then for all $x,y \in \R^d_+$,
\begin{align*}
	G^\kappa(x,y) \asymp  \begin{cases}
		\displaystyle \frac{1}{|x-y|^{1-\alpha}} 	\Big(\frac{x \wedge y}{|x-y|} \wedge 1\Big)^{q_\kappa} &\mbox{ if } \alpha<1,\\[9pt]
		\displaystyle			\Big(\frac{x \wedge y}{|x-y|} \wedge 1\Big)^{q_\kappa} \L \Big(\frac{(x \wedge y)\vee |x-y|}{|x-y|}\Big) &\mbox{ if }\alpha=1,\\[9pt]
		\displaystyle (x \wedge y)^{\alpha-1}	\Big(\frac{x \wedge y}{|x-y|} \wedge 1\Big)^{q_\kappa-\alpha+1}&\mbox{ if } \alpha>1.
	\end{cases}
\end{align*}	
\end{thm}
Note that when $d\ge 2$, at the threshold $q_{\kappa} = \alpha + \frac12 (\beta_1+\beta_2)$, there is a transition 
from the usual behavior of the Green function estimates to anomalous behavior. 

\medskip

It is natural and important to extend the results of this paper to smooth open sets other than the half-space. This is a pretty daunting task, one of the first challenges is the lack of scaling property for a general smooth open set. In the preprint \cite{CKSV24}, we have made significant progress in this direction.

We describe now the strategy for proving our main results and the organization of the paper.

It is well known that an appropriate Nash-type inequality implies the existence of the heat kernel (outside an exceptional set) and its $\alpha$-stable-type upper bound. So we start in Section \ref{s:Nash} with establishing 
a Nash-type inequality, see Proposition \ref{p:Y-Nash}. 
To this end, we consider a certain Feller process in $\R^d_+$ with continuous paths, subordinate it by an independent $\alpha/2$-stable subordinator, and show a Nash-type inequality for the Dirichlet form of the subordinate 
process. In case $\alpha<1$, one can estimate 
the Dirichlet form  of the subordinate process from above by $\sE^0$ and thus prove a Nash-type
inequality for $\sE^0$. In case $\alpha \ge 1$, we use  the scaling property of the process and  a truncation of the jump kernel together with 
the already obtained inequality for $\alpha <1$.

In Section \ref{s:regularity}, we 
prove parabolic H\"older regularity for both $\rY$ and $Y^{\kappa}$ and use this to remove the exceptional sets and extend $\overline{p}(t,x,y)$ and $p^{\kappa}(t,x,y)$ 
continuously
 to $(0,\infty)\times \overline{\R}^d_+\times \overline{\R}^d_+$, respectively $(0,\infty)\times \R^d_+\times \R^d_+$. 
To prove the parabolic H\"older regularity, we need an interior lower bound on the heat kernels. The proof of this lower bound is based on an argument that has already appeared in \cite[Section 3]{KSV-jump} and is given here in Proposition \ref{p:POTAl3.1}.
 Next, we prove several 
lower bounds on the heat kernels, mean exit times and exit distributions   
of the underlying process and its killed version that allow us to apply the standard arguments for establishing the parabolic H\"older regularity.
 We end Section \ref{s:regularity} with  the elementary Lemma \ref{l:general-upper-2} which will be used repeatedly
in deriving the heat kernel upper bounds.

Following the arguments from \cite{CKK09, CKSV}, in Section \ref{s:PHI} we establish the parabolic Harnack inequality, see Theorem \ref{t:phi}, and use it to establish the important preliminary off-diagonal lower bound $p^{\kappa}(t,x,y)\ge c tJ(x,y)$ for $t$ small compared to $|x-y|$, $x_d$ and $y_d$, see Proposition \ref{p:pld}.

In Section \ref{s:pupper} we prove the following preliminary upper bound:
$$
p^\kappa(t,x,y) \le C \Big(1 \wedge \frac{x_d}{t^{1/\alpha}}\Big)^{q_\kappa}\Big(1 \wedge \frac{y_d}{t^{1/\alpha}} \Big)^{q_\kappa} \Big( t^{-d/\alpha} \wedge \frac{t}{|x-y|^{d+\alpha}}\Big), \quad t>0, \ x,y\in \R^d_+,
$$
for all $q_{\kappa}\in [(\alpha-1)_+, \alpha+\beta_1)$, see Proposition \ref{p:UHK-rough}. 
Proposition \ref{p:UHK-rough} is proved in two steps:
Using Lemma \ref{l:general-upper-2} and several 
results from  \cite{KSV-jump, KSV-Green, KSV-nokill}
(see Lemmas \ref{l:(5.6)implies(5.1)}-\ref{l:UHKD_n}) 
we first prove Proposition \ref{p:UHK-rough} in case $q_{\kappa}<\alpha$.
The real challenge is to extend it to the full range of $q_{\kappa}$ by covering the 
case $q_{\kappa} \ge \alpha$.
For this we use the upper bound on Green potentials of powers of the distance to the boundary. In case $d\ge 2$, such bound was proved in \cite{KSV-Green, KSV-nokill}, and here we prove the corresponding result for $d=1$. The restriction $q_{\kappa}<\alpha$ is removed in Lemma \ref{l:UHKD} by using a bootstrap (induction) argument.

Section \ref{s:lowerb} is devoted to proving sharp heat kernel lower bounds.
The estimates are given in terms of a function $A_{b_1,b_2,b_3, b_4}(t,x,y)$, a space-time version of the function $B_{b_1,b_2,b_3, b_4}(x,y)$.
We first prove a preliminary lower bound, Proposition \ref{p:HKE-lower-1}, by using the semigroup property and some interior lower bound on the heat kernel. This bound turns out to be sufficient 
in the cases $d=1$ or $\beta_2<\alpha+\beta_1$. 
For the case $d\ge 2$ and $\beta_2\ge \alpha+\beta_1$,
we need a sharper lower bound obtained in the key Lemma \ref{l:lower}. There we apply the semigroup property 
(and the preliminary lower bound) 
on a carefully chosen interior set on which we have good control of 
the terms appearing in the preliminary lower bound. 
By combining Proposition \ref{p:HKE-lower-1} and Lemma \ref{l:lower}, we obtain in Proposition \ref{p:lower-1} the sharp lower bound in cases $\beta_2\neq \alpha+\beta_1$. The remaining case $\beta_2=\alpha+\beta_1$ is more delicate and is covered in Proposition \ref{p:lower-2}.

Sharp heat kernel upper bound is more difficult to establish and 
Section \ref{s:sharpupper}  is devoted to this task. The first step is to establish in Lemma \ref{l:UHK-case1-induction} an upper bound which includes the function $A_{\beta_1,0,\beta_3, 0}$. The proof uses Lemma \ref{l:general-upper-2} and an induction argument which consecutively increases the first parameter of the function $A$ until reaching $\beta_1$, thus making the decay successively sharper. The sharp upper bound
 in the case $d=1$ or $\beta_2<\alpha+\beta_1$, respectively $d=2$ and $\beta_2\ge \alpha+\beta_1$, 
 is given in Theorem \ref{t:UHK} and Corollary \ref{c:UHK}, respectively Theorem \ref{t:UHK2}. 
The proofs are based on Lemma \ref{l:general-upper-2}, Lemma \ref{l:UHK-case1-induction}, 
and a number of delicate technical lemmas involving multiple space-time integrals of 
the preliminary heat kernel estimates.

 In Section \ref{s:pmain},  
 we combine the upper bounds obtained 
in Section \ref{s:sharpupper} 
with the lower bounds from Section \ref{s:lowerb} and give the proofs of Theorems \ref{t:HKE}--\ref{t:HKE-kappa}.

It is well known that the Green function is the integral over time of the heat kernel. 
In Section \ref{s:Green} we first use this observation together with the estimates of 
$\overline{p}(t,x,y)$ obtained in Proposition \ref{p:upper-heatkernel} and Lemma \ref{l:NDL-bar-1} (see also Remark \ref{r:conti}) to prove Theorem \ref{t:green-1}. By using the same method of integrating 
the heat kernel estimates of $p^\kappa(t, x, y)$
over time we establish Theorem \ref{t:green-2}, thus reproving and extending the main result of \cite{KSV-Green}. This new proof sheds more light on 
the anomalous behavior of these Green function estimates.
Lemma \ref{l:int-smalltime} clearly shows that they are 
consequences of the different forms of the small time heat kernel estimates. 

The paper ends with an appendix which contains a number of technical results not depending
on the preliminary estimates of the heat kernel. 

Throughout this paper, the  constants $\beta_1$, $\beta_2$, $\beta_3$, $\beta_4$ 
will remain the same, and $\kappa$ always stands for a non-negative number.
The notation $C=C(a,b,\ldots)$ indicates that the constant $C$
depends on $a, b, \ldots$. The dependence on $\kappa, d$ and $\alpha$  may not be mentioned explicitly. 
Lower case letters 
$c_i, i=1,2,  \dots$ are used to denote the constants in the proofs
and the labeling of these constants starts anew in each proof. 
We denote by $m_d$   the Lebesgue measure on $\R^{d}$.
For Borel subset $D \subset \R^d$, $\delta_D(x)$ denotes the distance 
between $x$ and the boundary $\partial D$.

\medskip

\begin{center}
{\bf List of Notations}
\end{center}
$J(x,y):=|x-y|^{-d-\alpha}\sB(x,y)$,
\ \ 
$p_\alpha(t, x, y):= t^{-d/\alpha} \wedge \frac{t}{|x-y|^{d+\alpha}}$,\\
$\sE^0(u,v):= \frac{1}{2}\int_{\R^d_+} \int_{\R^d_+}  (u(x)-u(y))(v(x)-v(y)) J(x, y)dydx$,\\
$\sE^{\kappa}(u,v):=\sE^0(u,v)+\int_{\R^d_+}u(x)v(x)\kappa x_d^{-\alpha}\, dx$,\\
$\rF$:   the closure of  $C_c^\infty(\oR^d_+)$ in $L^2(\R^d_+,dx)$ under the norm  
$(\sE^0_1)^{1/2}$,\\
$\sF^0$:   the closure of  $C_c^\infty(\R^d_+)$ in $L^2(\R^d_+,dx)$ under the norm  
$(\sE^0_1)^{1/2}$,\\
 $\sF^\kappa:= 	\wt \sF^0 \cap L^2(\R^d_+, \kappa x_d^{-\alpha}dx)$,\\
$\rY$:  Hunt process associated with  $(\sE^0, \rF)$; \ \ $\bar p(t,x,y)$:    the heat kernel of $\rY$,\\
$Y^0$: Hunt process associated with  $(\sE^0, \sF^0)$; \ \ $p(t,x,y)$:    the heat kernel of $Y^0$,\\
$Y^{\kappa}:$   Hunt process associated with $(\sE^\kappa, \sF^\kappa)$; \ \ $p^{\kappa}(t,x,y)$:    the heat kernel of $Y^{\kappa}$,\\
$L_{\alpha}^{\sB}u(x):= \textrm{p.v.} \int_{\R^d_+}(u(y)-u(x)J(x,y)\, dy$; \ \  $L^{\sB}u(x):=L_{\alpha}^{\sB}u(x)-\kappa x_d^{-\alpha}u(x)$,
\begin{align*}
A_{b_1, b_2, b_3, b_4}(t,x,y)
&:= 	
\Big(\frac{(x_d \wedge y_d) \vee t^{1/\alpha}}{|x-y|} \wedge 1 \Big)^{b_1}  \Big(\frac{(x_d \vee y_d) \vee t^{1/\alpha}}{|x-y|} \wedge 1 \Big)^{b_2}  \\	
&  \times  \L^{b_3}\Big(\frac{ ((x_d \vee y_d) \vee t^{1/\alpha}) \wedge |x-y|}{((x_d \wedge y_d) \vee t^{1/\alpha}) \wedge |x-y|} \Big) \,\L^{b_4}\Big(\frac{ |x-y|}{( (x_d \vee y_d) \vee t^{1/\alpha}) \wedge |x-y|} \Big),
\end{align*}
$B_{b_1, b_2, b_3, b_4}(x,y):=A_{b_1, b_2, b_3, b_4}(0,x,y)$.

\section{Nash inequality and the existence of the heat kernel}\label{s:Nash}

The goal of this section is (i) to prove a Nash type inequality  
(Proposition \ref{p:Y-Nash}); (ii) to deduce the existence of the heat kernels of $\rY$ and $Y^\kappa$; and (iii) to establish their  preliminary upper bounds (Proposition \ref{p:upper-heatkernel}).

We begin the section 
by introducing the notation for the relevant semigroups and establishing their scale invariance and horizontal translation invariance.

For $\kappa \ge 0$, 
 we denote by $(\rP_t)_{t\ge0}$ and  $(P_t^\kappa)_{t\ge 0}$ the semigroups corresponding to $(\sE^0, \overline \sF)$ and $(\sE^\kappa, \sF^\kappa)$ respectively.   $(\rP_t)_{t\ge0}$ and  $(P_t^\kappa)_{t\ge 0}$ define 
 contraction semigroups
 on 
 $L^p(\oR^d_+,dx)=L^p(\R^d_+,dx)$   for every $p \in [1,\infty]$, and when $p\in [1, \infty)$, these semigroups are strongly continuous. 
 For $t>0$ and  $p,q \in [1,\infty]$, define
\begin{align*}
\lVert \rP_t\rVert_{p \to q}= \sup\Big\{ \lVert \rP_t f \rVert_{L^q(\R^d_+,dx)} :f \in L^p(\R^d_+,dx),\,\lVert f \rVert_{L^p(\R^d_+,dx)}\le 1 \Big\}. 
\end{align*}

For $f:\R^d \to \R$ and $r>0$, define $f^{(r)}(x)=f(rx)$. The  following scaling property of $(P_t^\kappa)_{t\ge0}$ comes from  \cite[Lemma 5.1]{KSV-jump} and \cite[Lemma 2.1]{KSV-nokill}. By the same proof, $(\rP_t)_{t\ge0}$ also has the 
same  scaling property. 

\begin{lemma}\label{l:scaling}
Let $p \in [1,\infty]$ and  $\kappa \ge 0$.	For any $f \in L^p(\R^d_+, dx)$, $t>0$ and $r>0$, we have 
$$ \overline P_t f(x)=\overline P_{r^{-\alpha }t}f^{(r)}(x/r) \quad \text{and} \quad  P^\kappa_t f(x)= P^\kappa_{r^{-\alpha }t}f^{(r)}(x/r)  \quad \text{in } \, L^p(\R^d_+,dx).
$$
In particular, we have
\begin{align}\label{e:semigroup-scaling}
	\lVert \rP_t\rVert_{1\to \infty } = t^{-d/\alpha}	\lVert \rP_1\rVert_{1 \to \infty } \quad \text{for all }\, t>0. 
\end{align}
\end{lemma}

Let $d\ge 2$. For $f:\R^d \to \R$ and $\wt z \in \R^{d-1}$, define $f_{z}(x)=f(x+(\wt z,0))$.
From {\bf (A4)}, we also get the following horizontal translation invariance property of the semigroups $(\rP_t)_{t\ge0}$ and  $(P_t^\kappa)_{t\ge 0}$.

\begin{lemma}\label{l:translation-invariance}
	Let $d\ge 2$, $p \in [1,\infty]$ and  $\kappa \ge 0$.	For any $f \in L^p(\R^d_+,dx)$, $t>0$ and $\wt z \in \R^{d-1}$, we have 
	$$ \overline P_t f(x)=\overline P_{t}f_{z}(x - (\wt z,0)) \quad \text{and} \quad  P^\kappa_t f(x)= P^\kappa_{t}f_z(x- (\wt z,0))  \quad \text{in } \, L^p(\R^d_+,dx).$$
\end{lemma}

A consequence of the next lemma is that, in case when $\alpha\le 1$, the process $\rY$ started from $\R^d_+$ will never hit $\partial \R^d_+$ and is equal to $Y^0$.

\begin{lemma}\label{l:alpha<1}
If $\alpha \le 1$, then   
$\sF^0=\rF$.
\end{lemma}
\begin{proof} Define
\begin{align*}
	\wt \sC(u,v)&:= \frac{1}{2}\int_{\R^d_+} \int_{\R^d_+}  \frac{(u(x)-u(y))(v(x)-v(y))}{|x-y|^{d+\alpha}}dydx,\\
	\sD(\wt \sC) &:=\text{closure of $C_c^\infty(\overline\R^d_+)$ in $L^2(\R^d_+,dx)$ under $\wt \sC+(\cdot, \cdot)_{L^2(\R^d_+,dx)}$}.
\end{align*}
Then $(\wt \sC, \sD(\wt \sC))$ is a regular Dirichlet form associated with the \textit{reflected} $\alpha$-stable process in $\overline\R^d_+$ in the sense of \cite{BBC}. 
By  {\bf (A3)}{(I)},   $\sE^0(u,u) \le C_2 \wt \sC(u,u)$ for all $u \in C_c^\infty(\overline \R^d_+)$ 
and hence $\sD(\wt \sC) \subset \rF$. By \cite[Theorem 2.5(i) and Remark 2.2(1)]{BBC}, since $\alpha \le 1$, $\overline \R^d_+ \setminus \R^d_+$  is 
$(\wt \sC, \sD(\wt \sC))$-polar 
and hence is also $(\sE^0, \rF)$-polar.  
Therefore, when starting from $\R^d_+$,  $\rY$ will never exit $\R^d_+$. Hence 
$\rY$ and $Y^0$ are the same when they  
start from $x \in \R^d_+$ and thus $\sF^0=\rF$. 
\end{proof}

In order to prove the Nash type inequality, we first consider a Brownian motion on $\R^d_+$ 
killed with a critical potential 
and a subordinate process obtained by time changing this killed Brownian motion with an independent $\alpha/2$-stable subordinator. Then using results 
from \cite{CKSV},   the Hardy inequality in 
\cite[Proposition 3.2]{KSV-nokill} and comparing the Dirichlet form corresponding to the subordinate process with $(\sE^0, \rF)$, we deduce the desired result.

For any $\gamma  \ge 0$, denote by $I_\gamma$ the modified Bessel function of the  first kind defined by
$$
I_\gamma(r)= \sum_{m=0}^\infty \frac{1}{m! \, \Gamma(\gamma+1+m)} \left(\frac{r}{2}\right)^{2m+\gamma},
$$
where $\Gamma(r):=\int_0^\infty u^{r-1}e^{-u}du$ is the Gamma function. It is known that  (see, e.g. \cite[(9.6.7) and (9.7.1)]{AS})
\begin{equation}\label{e:estimates-modifiedBessel}	I_\gamma(r) \asymp (1 \wedge r)^{\gamma+1/2} r^{-1/2}e^r, \quad  r>0.
\end{equation}

Define for $t>0$ and 
$x=(x_1, \dots, x_d), y=(y_1, \dots, y_d) \in \R^d_+$, 
\begin{equation*}
	q^\gamma(t,x,y)= \frac{\sqrt{x_dy_d}}{2t} I_\gamma \left( \frac{x_dy_d}{2t}\right) \exp \left( - \frac{x_d^2+y_d^2}{4t} \right) \prod_{i=1}^{d-1} \left[\frac{1}{\sqrt{4\pi t}} \exp \left( - \frac{|x_i-y_i|^2}{4t} \right)\right].
\end{equation*}
Note that by  \eqref{e:estimates-modifiedBessel},
\begin{equation}\label{e:W-HKE}	
	q^\gamma(t,x,y) \asymp 
	\left( 1 \wedge \frac{x_dy_d}{t}\right)^{\gamma+1/2}t^{-d/2}\exp \left( - \frac{|x-y|^2}{4t}	\right).
\end{equation}
By \cite[Lemma 4.1 and Theorem 4.9]{MNS},  $q^\gamma(t,x,y)$ is 
the transition density  of the Feller process 
(killed Brownian motion with critical potential) 
 $W^\gamma=(W^\gamma_t)_{t\ge0}$ 
on $\R^d_+$ associated with the following regular Dirichlet form $(\sQ^\gamma, \sD(\sQ^\gamma))$:
\begin{align*}
\sQ^\gamma(u,v)&:=
 \int_{\R^d_+}   \nabla u(x) \cdot \nabla v(x)dx + \Big(\gamma^2- \frac14\Big)  \int_{\R^d_+} u(x)v(x)x_d^{-2} dx,\\ 
\sD(\sQ^\gamma)&:=\text{closure of $C_c^\infty(\R^d_+)$ in $L^2(\R^d_+,dx)$ under $\sQ^\gamma_1=\sQ^\gamma+(\cdot,\cdot)_{L^2(\R^d_+,dx)}$}.
\end{align*}

Let $S=(S_t)_{t\ge0}$ be an $\alpha/2$-stable subordinator independent of $W^\gamma$, and let $X^\gamma=(X^\gamma_t)_{t\ge0}$ be the subordinate process $X^\gamma_t:=W^\gamma_{S_t}$. Then $X^\gamma$ is a Hunt process with no diffusion part. The transition density $p^\gamma(t,x,y)$ of $X^\gamma_t$ exists and is given by
\begin{equation*}
	p^\gamma(t,x,y)=\int_0^\infty q^\gamma(s,x,y) \frac{d}{ds}\P(S_t \le s).
\end{equation*}
Also, the jump kernel $J^\gamma(dx,dy)$ and the killing measure $\kappa^\gamma(dx)$ of  $X^\gamma$ have densities $J^\gamma(x,y)$ and $\kappa^\gamma(x)$ that are given by the following formulas 
(see, for instance, \cite[(2.8)--(2.9)]{Ok}):
\begin{align*}
	J^\gamma(x,y)&=\int_0^\infty q^\gamma(t,x,y) \, \nu_{\alpha/2}(t)dt, \quad\;\; \kappa^\gamma(x)=\int_0^\infty  \bigg( 1- \int_{\R^d_+} q^\gamma(t,x,y)dy \bigg)  \nu_{\alpha/2}(t)dt,
\end{align*}
where $\nu_{\alpha/2}(t)=\frac{\alpha/2}{\Gamma(1-\alpha/2)}t^{-1-\alpha/2}$ is the L\'evy density of the $\alpha/2$-stable subordinator $S$.
\begin{lemma}\label{l:Z-estimates}
	(i) There exists a constant $c_{\gamma, \alpha}>0$ such that $	\kappa^\gamma(x) = c_{\gamma, \alpha} x_d^{-\alpha}$ for every $x \in \R^d_+$.
		
	\noindent
	(ii) It holds that for any $x,y \in \R^d_+$,
	\begin{align*}
		J^\gamma(x,y) &\asymp  \Big( 1 \wedge \frac{x_d}{|x-y|}\Big)^{\gamma+1/2} \Big( 1 \wedge \frac{y_d}{|x-y|}\Big)^{\gamma+1/2}   \frac{1}{ |x-y|^{d+\alpha}}=  \frac{B_{\gamma+1/2,\gamma+1/2,0,0}(x,y)}{ |x-y|^{d+\alpha}}.
	\end{align*}
	
	\noindent (iii) There exists a constant $C>0$ such that 
	\begin{equation*}
		p^\gamma(t,x,y) \le C t^{-d/\alpha}, \quad \, t>0, \; x,y \in \R^d_+.
	\end{equation*}
\end{lemma}
\begin{proof} (i) 
By the scaling property and horizontal translation invariance of $W^\gamma$, using the change of variables $z=(y-(\wt x, 0))/x_d$, we get that for every $x \in \R^d_+$,
\begin{align*}
 \int_{\R^d_+} q^\gamma(t,(\wt 0, x_d),(\wt y - \wt x, y_d))dy &= x_d^{d}\int_{\R^d_+} q^\gamma(t,(\wt 0, x_d),  x_dz)dz=\int_{\R^d_+} q^\gamma(x_d^{-2}t,(\wt 0, 1), z)dz.\end{align*}
Therefore,  for every $x \in \R^d_+$, using the change of variables $s=x_d^{-2}t$, we get that
\begin{align*}
	\kappa^\gamma(x) &=\frac{\alpha/2}{\Gamma(1-\alpha/2)}\int_0^\infty  \bigg( 1- \int_{\R^d_+} q^\gamma(t,(\wt 0, x_d),(\wt y - \wt x, y_d))dy \bigg)  t^{-\alpha/2-1}dt\\
	&=\frac{\alpha/2}{\Gamma(1-\alpha/2)} x_d^{-\alpha}\int_0^\infty  \bigg( 1- \int_{\R^d_+} q^\gamma(s,(\wt 0, 1), z)dz \bigg)  s^{-\alpha/2-1}ds=  \kappa^\gamma((\wt 0,1)) x_d^{-\alpha}.
\end{align*}

(ii) Note that \eqref{e:W-HKE} can be written in the form
\begin{align*}
q^{\gamma}(t,x,y)\asymp h(t,x,y)\frac{1}{|B(x,\Phi^{-1}(t)|}\exp\left(-\frac{|x-y|^2}{\Phi^{-1}(t)^2}\right)\, \quad \text{for all }(t,x,y)\in (0,\infty)\times \R^d_+\times \R^d_+, 
\end{align*}
with $h(t,x,y) = (1\wedge (x_dy_d/t))^{\gamma+1/2}$, $\Phi(r)=r^2$, and $|B(x,\Phi^{-1}(t)|\asymp t^{-d/2}$ the volume of the ball centered at $x$ and radius $\Phi^{-1}(t)=t^{-d/2}$. This precisely means that the  condition $\mathbf{HK_U^h}$ in \cite{CKSV} holds with $C_0=0$,  $\Phi(r)=r^2$ and $h(t,x,y) = (1 \wedge (x_dy_d/t))^{\gamma+1/2}$. 

Note that  the tail of the L\'evy measure of $S$ is given by $ \nu_{\alpha/2}(r,\infty)=:\int_r^\infty \nu_{\alpha/2}(u)du = r^{-\alpha/2}/ \Gamma(1-\alpha/2)$ implying that $\nu_{\alpha/2}(r,\infty)/\nu_{\alpha/2}(R,\infty)=(R/r)^{\alpha/2}$ 
for all $0<r\le R<\infty$. Thus,  the scaling  condition  {\bf (Poly-$\infty$)} in \cite{CKSV} also holds and  we obtain from \cite[Theorem 4.1]{CKSV} that for $x,y \in \R^d_+$,
\begin{align*}
|x-y|^{d+\alpha}J^\gamma(x,y)\asymp  \Big( 1 \wedge \frac{x_dy_d}{|x-y|^2}\Big)^{\gamma+1/2} \asymp  \Big( 1 \wedge \frac{x_d}{|x-y|}\Big)^{\gamma+1/2}  \Big( 1 \wedge \frac{y_d}{|x-y|}\Big)^{\gamma+1/2}.
\end{align*}
In the second comparison, we used the fact that $(1 \wedge \frac{x_d}{|x-y|}) (1 \wedge \frac{y_d}{|x-y|}) \le (1 \wedge \frac{x_dy_d}{|x-y|^2}) \le 2(1 \wedge \frac{x_d}{|x-y|})(1 \wedge \frac{y_d}{|x-y|})$ for all $x,y \in \R^d_+$, which can be  proved 
 by using  $y_d \le x_d + |x-y|$. 

(iii) Since  the  conditions {\bf (Poly-$\infty$)} and $\mathbf{HK_U^h}$ in \cite{CKSV} hold, the result follows from \cite[Proposition 4.5(ii)]{CKSV}. \end{proof}

Denote by $(\sC^\gamma, \sD(\sC^\gamma))$ the regular Dirichlet form associated with the subordinate process $X^\gamma$. Then since $X^\gamma$ has no diffusion part, we get from Lemma \ref{l:Z-estimates}(i) that 
\begin{equation*}
	\sC^\gamma(u,v)=\int_{\R^d_+} \int_{\R^d_+} (u(x)-u(y))(v(x)-v(y))J^\gamma(x,y)dydx + c_{\gamma, \alpha}\int_{\R^d_+} u(x)v(x) x_d^{-\alpha}dx.
\end{equation*}
Also, we have $C_c^\infty(\R^d_+) \subset \sD(\sC^\gamma)$ since $\sD(\sQ^\gamma)\subset\sD(\sC^\gamma)$, see  \cite{Ok}.

\begin{lemma}\label{l:Z-Nash}
	There exists a constant $C>0$ such that 
	\begin{equation*}
		\lVert u \rVert_{L^2(\R^d_+,dx
			)}^{2(1+\alpha/d)} \le C \sC^\gamma(u,u) \quad \text{for every } \; u \in C_c^\infty(\R^d_+) \text{ with } \lVert u \rVert_{L^1(\R^d_+,dx)
		} \le 1.
	\end{equation*}
\end{lemma}
\begin{proof} By \cite[Theorem 2.1]{CKS} (see also \cite[Theorem 3.4]{CKKW} and \cite[Theorem II.5]{Co}), the result follows from  Lemma \ref{l:Z-estimates}(iii). \end{proof}

\begin{prop}\label{p:Y-Nash}
	There exists a constant $C>0$ such that 
	\begin{equation}\label{e:Y-Nash}
		\lVert u \rVert_{L^2(\R^d_+,dx)}^{2(1+\alpha/d)} \le C \sE^0(u,u) \quad \text{for every } \; u \in \rF \text{ with } \lVert u \rVert_{L^1(\R^d_+,dx)} \le 1.
	\end{equation}
\end{prop}
\begin{proof}  
 We first assume that $\alpha < 1$. 
Let $\gamma=\beta_1\vee \beta_2$. 
 Using Lemmas \ref{l:Z-Nash} and 
\ref{l:Z-estimates}(i)--(ii),  the Hardy inequality in \cite[Proposition 3.2]{KSV-nokill} 
and {\bf (A3)}(I), 
we get that for any $u \in C_c^\infty(\R^d_+)$ with $\lVert u \rVert_{L^1(\R^d_+,dx)} \le 1$,
\begin{align*}
		\lVert u \rVert_{L^2(\R^d_+,dx)}^{2(1+\alpha/d)} &\le c_1 \int_{\R^d_+}\int_{\R^d_+} (u(x)-u(y))^2 \frac{B_{\gamma+1/2, \gamma+1/2, 0, 0}(x,y)}{|x-y|^{d+\alpha}}dydx 
		+ c_{\gamma, \alpha} \int_{\R^d_+}u(x)^2 x_d^{-\alpha}dx\\
 & \le c_2 \sE^0(u,u), 
\end{align*}
where $B_{\gamma+1/2,\gamma+1/2,0,0}$ is defined in \eqref{def:wtB}.  By Lemma \ref{l:alpha<1}, $\rF$ is the closure of $C_c^\infty(\R^d_+)$ under $\sE^0_1$. Therefore,  we conclude that  \eqref{e:Y-Nash} is true when $\alpha<1$.

Now, we assume that $\alpha \ge 1$.  Since \eqref{e:Y-Nash} is valid when $\alpha < 1$ and $C_c^\infty(\overline\R^d_+) \subset \rF$, 
we get that for every $u \in C_c^\infty(\overline\R^d_+)$ with $ \lVert u \rVert_{L^1(\R^d_+,dx)} \le 1$,
\begin{align*}
	\sE^{0}(u,u) 	&
	\ge\frac12 \int_{\R^d_+} \int_{\R^d_+}    (u(x)-u(y))^2 \frac{\sB(x,y)}{|x-y|^{d+1/2}} (1 -\1_{\{|x-y|> 1\}}  )  dydx  \\
	& \ge 
	c_3	\lVert u \rVert_{L^2(\R^d_+,dx)}^{2(1+1/(2d))} - 
	2 C_2\lVert u \rVert^2_{L^2(\R^d_+,dx)} 
	 \sup_{x \in \R^d_+} \int_{\R^d_+, |x-y| >1} \frac{dy}{ |x-y|^{d+1/2}} \\
		& \ge 
		c_3	\lVert u \rVert_{L^2(\R^d_+,dx)}^{2(1+1/(2d))} -
		c_4\lVert u \rVert^2_{L^2(\R^d_+,dx)}.
\end{align*}
Since $\rF$ is the closure of $C_c^\infty(\overline\R^d_+)$ under $\sE^0_1$,  \cite[Theorem 2.1]{CKS} yields that 
 $\lVert \rP_t\rVert_{1\to \infty}  \le c_5 t^{-2d} e^{c_4t}$    for all $t>0$. 
By \eqref{e:semigroup-scaling}, it follows that 
$\lVert \rP_t\rVert_{1\to \infty }  =t^{-d/\alpha} 	\lVert \rP_1\rVert_{1\to \infty}  \le c_5e^{c_4} t^{-d/\alpha}$  for all $t>0$. Using  \cite[Theorem 2.1]{CKS}  again, we conclude that \eqref{e:Y-Nash} holds for $\alpha\ge1$ and finish the proof. \end{proof}

As a consequence of the Nash-type inequality \eqref{e:Y-Nash}, we get the existence and a priori upper bounds of the heat kernels $\bar p(t,x,y)$ and $p^\kappa(t,x,y)$ of $\rY$ and $Y^\kappa$ respectively. 
Recall that
$p_\alpha(t, x, y):= t^{-d/\alpha} \wedge \frac{t}{|x-y|^{d+\alpha}}$

\begin{prop}\label{p:upper-heatkernel}
Let $\kappa\ge0$. The processes $\rY$ and $Y^\kappa$ have heat kernels $\bar p(t,x,y)$ and $p^\kappa(t,x,y)$ defined on $(0,\infty) \times (\overline \R^d_+ \setminus \sN) \times (\overline \R^d_+ \setminus \sN)$ 
and $(0,\infty) \times (\R^d_+ \setminus \sN) \times (\R^d_+ \setminus \sN)$ 
respectively, where $\sN \subset \overline\R^d_+$ is a properly exceptional set for
 $\rY$.
Moreover, there exists a constant $C>0$ such that
\begin{equation}\label{e:rel-heatkernel}
p^\kappa(t,x,y) \le \bar	p(t,x,y), \quad t>0, x,y \in \R^d_+ \setminus \sN
\end{equation}
and
\begin{equation}\label{e:upper-heatkernel}
\bar	p(t,x,y) \le C 
p_\alpha(t, x, y),
\quad t>0, \; x,y \in \overline\R^d_+ \setminus \sN.
\end{equation}
\end{prop}
\begin{proof} 
By our Nash-type inequality  \eqref{e:Y-Nash} and \cite{BBCK}, 
$\rY$ has a heat kernel $\bar p(t,x,y)$ on $(0,\infty) \times (\overline \R^d_+ \setminus \sN) \times (\overline \R^d_+ \setminus \sN)$ for a properly exceptional set $\sN$ and 
\begin{equation}\label{e:upper-heatkernel1}
\bar	p(t,x,y) \le c  t^{-d/\alpha}, \quad t>0, \; x,y \in \overline\R^d_+ \setminus \sN.
\end{equation}
Since $\sB(x,y)$ is bounded above, using \eqref{e:upper-heatkernel1}, we can follow the arguments given in \cite[Example 5.5]{CKKW} line by line and conclude that 
\eqref{e:upper-heatkernel} holds. 
According to \cite{KSV-nokill} (see the discussion before Lemma 2.1 there),
$Y^\kappa$ can be realized as a subprocess of $\rY$. 
Thus, 
 $Y^\kappa$ has 
a heat kernel $p^\kappa(t,x,y)$ defined on
$(0,\infty) \times (\R^d_+ \setminus \sN) \times (\R^d_+ \setminus \sN)$ and 
we also obtain 
\eqref{e:rel-heatkernel}. 
\end{proof}

 For notational convenience, we extend 
the domain of $p^\kappa(t,x,y)$ 
 to $(0,\infty) \times (\overline \R^d_+ \setminus \sN) \times (\overline \R^d_+ \setminus \sN)$ by letting $p^\kappa(t,x,y)=0$ if 
$x \in \partial \R^d_+ \setminus \sN$ or $y \in \partial \R^d_+ \setminus \sN$.

\section{Parabolic  H\"older regularity and consequences}\label{s:regularity}

For $\kappa \ge 0$ and an open set $D \subset \oR^d_+$ relative to the  topology on $\oR^d_+$,
 we  denote by $\bar{Y}^{D}$ and $\bar{P}^{D}_t$ the part of the process $\bar{Y}$ killed upon exiting $D$ and
  its semigroup, respectively,
  and by  $Y^{\kappa,D}$ and $P_t^{\kappa,D}$ the part of the process $Y^\kappa$ killed upon exiting $D \cap \R^d_+$ and its semigroup, respectively.    
The Dirichlet forms of $\bar{Y}^{D}$ and  $Y^{\kappa,D}$  are  $(\sE^0,\bar{\sF}^D)$ and $(\sE^\kappa,\sF^{\kappa, D})$, where
 $\bar{\sF}^D=\{u\in \bar{\sF}: u=0 \mbox{ q.\/e. on } \oR^d_+\setminus D\}$ and $\sF^{\kappa, D}=\{u\in \sF^\kappa: u=0 \mbox{ q.\/e. on } \R^d_+\setminus D\}$ respectively. 
For $u, v\in \sF^{\kappa, D}$, 
\begin{equation}\label{e:POTA3.1}
\sE^\kappa(u, v)=\frac12\int_D\int_D((u(x)-u(y))(v(x)-v(y))J(x, y)dydx+\int_Du(x)v(x)\kappa_D(x)dx,
\end{equation}
where
\begin{equation}\label{e:POTA3.2}
\kappa_D(x)=\int_{\R^d_+\setminus D}J(x, y)dy+\kappa x_d^{-\alpha}.
\end{equation}
Let $\bar \tau_D= \inf\{t>0: \rY_t \notin D\}$ and  $\tau^\kappa_D=\inf\{t>0: Y^\kappa_t \notin D \cap \R^d_+\}$. 
For $x,y\notin \mathcal N$, let 
\begin{align}
	\bar p^D(t,x,y) &= \bar p(t,x,y) - \E_x\big[\bar p(t- \bar \tau_D, \rY_{\overline\tau_D},y); \, \bar\tau_D<t\big],\nn\\
	p^{\kappa, D}(t,x,y) &= p^\kappa (t,x,y) - \E_x\big[ p^\kappa (t-  \tau_D^\kappa, Y^\kappa_{\tau^\kappa_D},y); \, \tau_D^\kappa<t\big].\label{e:Dirichlet-kernel}
\end{align}
By the strong Markov property, $\bar p^D(t,x,y)$ and $p^{\kappa,D}(t,x,y)$ are  the transition densities of $\rY^D$ and $Y^{\kappa, D}$ respectively.

In case when $D$ is a relatively compact open subset of $\R^d_+$,  one can show that the process $Y^{\kappa,D}$ can start from every point in $D$ and provide an interior lower bound for its transition density. We accomplish this by identifying the semigroup $P^{\kappa, D}_t$ with the Feynman-Kac semigroup  of the part process on $D$ of an auxiliary process.
This idea has already been used in \cite[Subsection 3.1]{KSV-jump}. For the benefit of the reader, we repeat some of the details here. 
 Note that, unlike  \cite[Subsection 3.1]{KSV-jump}, {\bf (A3)}{(II)} is not assumed.

 Recall that we denote by $m_d$   the Lebesgue measure on $\R^{d}$. 
\begin{prop}\label{p:POTAl3.1}
If $D$ is a relatively compact open subset of $\R^d_+$, then $Y^{\kappa,D}$ 
has a  transition density 
$p^{\kappa, D}(t, x, y)$ defined for any $(t, x, y)\in (0, \infty)\times D\times D$.
Furthermore, for any $T>0$ and $b\in (0, 1]$, there exists 
a constant $C=C(b, T, D)>0$ such that
$$
p^{\kappa  ,D}(t,x,y)
\ge C
p_\alpha(t, x, y),
\quad t\in (0,T), \; x,y \in D 
\textrm{ with } 
\delta_D(x)\wedge \delta_D(y)>bt^{1/\alpha}.
$$
\end{prop}
\begin{proof}
For $\gamma>0$ let $J_{\gamma}:\R^d\times \R^d\to [0,\infty]$ be defined by $J_{\gamma}(x,y)=J(x,y)$ if $x,y\in D$, and $J_{\gamma}(x,y)=\gamma|x-y|^{-\alpha-d}$ otherwise. It follows from \textbf{(A3)}(I) and the relative compactness of $D$ that $J_{\gamma}(x,y)\asymp |x-y|^{-\alpha-d}$. Hence  by \cite[Theorem 1.2]{CK08}, there exists a Feller and strongly Feller process $Z$ (that can start from every point in $\R^d)$ with a continuous transition density $\wt{q}(t,x,y)$ on $(0,\infty)\times \R^d \times \R^d$ 
such that for all $t>0, \; x,y \in \R^d$,
\begin{equation}\label{e:alpha-stable-estimate}
		c_1^{-1} 
		p_\alpha(t, x, y)
		\le \wt q(t,x,y) \le c_1 
		p_\alpha(t, x, y)
	\end{equation} 
for some constant $c_1\ge 1$.

Denote the part of the process $Z$ killed upon exiting $D$ by $Z^D$. The Dirichlet form of $Z^D$ is $(\sC, \mathcal{D}_D(\sC))$, where for $u,v\in \sD_D(\sC)$,
\begin{eqnarray*}
\sC(u,v)&=&\frac{1}{2}\int_D \int_D (u(x)-u(y))(v(x)-v(y))J_{\gamma}(x,y)\, dy\, dx +\int_D u(x)v(x)\kappa^{Z}_D(x)\, dx\\
&=&\frac{1}{2}\int_D \int_D (u(x)-u(y))(v(x)-v(y))J(x,y)\, dy\, dx +\int_D u(x)v(x)\kappa^{Z}_D(x)\, dx
\end{eqnarray*}
with
\begin{equation}\label{e:kappa-Z-new}
\kappa^{Z}_D(x)=\int_{\R^d\setminus D}J_{\gamma}(x,y)\, dy=\gamma \int_{\R^d\setminus D}
|x-y|^{-d-\alpha}\, dy\, ,\quad x\in D\, ,
\end{equation}
and $\mathcal{D}_D(\sC)=\{u\in \mathcal{D}(\sC):\, u=0 \textrm{ q.e.~on  } \R^d\setminus D\}$.

Let $\delta_D=\mathrm{dist}(D, \partial \R^d_+)$ and let $V$ be the $\delta_D/2$-neighborhood of $D$, that is $V:=\{x\in \R^d_+:\, \mathrm{dist}(x, D)<\delta_D/2\}$. Then
$$
\kappa_D(x)=\int_{\R^d_+\setminus V}J(x,y)\, dy+\int_{V\setminus D}J(x,y)\, dy +\kappa x_d^{-\alpha}\, .
$$
It follows from \textbf{(A3)}(I) and the relative compactness of $D$ that $c_2 |x-y|^{-d-\alpha}\le J(x,y) \le C_2 |x-y|^{-d-\alpha}$ for all $x,y\in V$ with $c_2:=c_2(D)>0$.
It is quite  easy to see that 
$\sup_{x\in  D}\int_{\R^d_+\setminus V}J(x,y)\, dy =:c_3\ < \infty$. Therefore 
$$
c_2\int_{V\setminus D}|x-y|^{-d-\alpha}\, dy \le \kappa_D(x)\le c_3 +C_2 \int_{V\setminus D}|x-y|^{-d-\alpha}dy+c_{4}\, ,\quad x\in U\, ,
$$
where $c_4:=\kappa\sup_{x\in D}x_d^{-\alpha}$. Since 
$$
\inf_{x\in D}\int_{V\setminus D}|x-y|^{-d-\alpha}\, dy \ge  m_d(V \setminus D) \,\mathrm{diam}(V)^{-d-\alpha} =:c_{5}>0\, ,
$$
we conclude that
$$
c_2\int_{V\setminus D}|x-y|^{-d-\alpha}\, dy \le \kappa_D(x)\le c_{6} \int_{V\setminus D}|x-y|^{-d-\alpha}\, dy\, ,\quad x\in D\, .
$$
Further, since
$$
\kappa_D^Z(x)=\gamma \bigg(\int_{\R^d\setminus V}  |x-y|^{-d-\alpha}\, dy+\int_{V\setminus D} |x-y|^{-d-\alpha}\, dy\bigg)\, ,\quad x\in D\, 
$$
and $\sup_{x\in D}\int_{\R^d\setminus V}|x-y|^{-d-\alpha}\, dy=:c_{7}<\infty$, we see that there is a constant $c_{8}>0$ such 
$$
\int_{V\setminus D}|x-y|^{-d-\alpha}\, dy \le \gamma^{-1}\kappa_D^Z(x)\le c_{8} \int_{V\setminus D}|x-y|^{-d-\alpha}\, dy\, ,\quad x\in 
D.
$$
It follows that
$c_{6}^{-1}\kappa_D(x)\le \gamma^{-1}\kappa_D^Z(x)\le  c_{8}c_2^{-1}\kappa_{D}(x)$ 
for all $x\in D$
with positive constants $c_2, c_{6}, c_{8}$ not depending on $\gamma$. Now we choose $\gamma>0$ so small that $\gamma c_{8}c_2^{-1}\le 1$. With this choice we get that $\kappa_D^Z(x)\le \kappa_D(x)$ for all $x\in D$. In particular, with $c_{9}:=\gamma c_{6}^{-1}$ we see that
\begin{equation}\label{e:kappa-U-kappa-Z-comparable}
c_{9}\kappa_D(x)\le \kappa_D^Z(x)\le \kappa_D(x)\, , \qquad x\in D\, .
\end{equation}
It follows that for $u\in C_c^{\infty}(D)$,
\begin{eqnarray*}
\sE^{\kappa, D}_1(u,u)
&=&\frac{1}{2}\int_D \int_D (u(x)-u(y))^2 J(x,y)\, dy\, dx +\int_D u(x)^2 \kappa_D(x)\, dx+\int_{D} u(x)^2 dx \\
&\asymp & \frac{1}{2}\int_D \int_D (u(x)-u(y))^2 J_{\gamma}(x,y)\, dy\, dx +\int_D u(x)^2 \kappa_D^Z(x)\, dx+\int_D u(x)^2 dx\\
& =&\sC(u,u) +\int_D u(x)^2 dx=\sC_1(u,u)\, .
\end{eqnarray*}
Since $C_c^{\infty}(D)$ is a core of both $(\sE^{\kappa}, \sF^{\kappa, D})$ and $(\sC, \mathcal{D}_D(\sC))$, we have that $\sF^{\kappa, D}=\mathcal{D}_D(\sC)$.

Define $\wt{\kappa}:D\to \R$ by 
$
\wt{\kappa}(x):=\kappa_D(x)-\kappa_D^Z(x)$,  $x\in D.
$
By the choice of $\gamma$ we have that $\wt{\kappa}\ge 0$. On the other hand, 
it follows from \eqref{e:kappa-Z-new}  and \eqref{e:kappa-U-kappa-Z-comparable} that there is a constant $c_{10}>0$ such that
\begin{equation}\label{e:estimate-kappaU}
\wt{\kappa}(x)\le \kappa_{D}(x)\le  c_{10}\delta_D(x)^{-\alpha}\, , \quad x\in D.
\end{equation}
Let $\mu(dx)=\wt{\kappa}(x)\, dx$ be a measure on $D$. Using \eqref{e:alpha-stable-estimate} and \eqref{e:estimate-kappaU} one can  check that 

$$
N^{D,\mu}_a(t):=\sup_{x\in \R^d_+} \int_0^t \int_{z\in D, \delta_D(z)>at^{1/\alpha}} \wt{q}(s,x,z)\, \mu(dz)\, ds, \quad t>0, a\ge 0,
$$
satisfies that $\sup_{t<1} N^{D, \mu}_a(t)<\infty$ for all $a\in (0,1)$, and $\lim_{t\to 0}N^{U,\mu}_0(t)=0$ for every relatively compact open subset $U$ of $D$. This means that 
$\mu\in {\bf K}_1(D)$, where the class ${\bf K}_1(D)$ is defined in \cite[Definition 2.12]{CKSV20}. 

For any Borel function $f:D\to [0,\infty)$ let  
\begin{equation}\label{e:Feynman-Kac}
T^{D, \wt{\kappa}}_t f(x) = \E_x \left[\exp \left( - \int_0^t \wt \kappa(Z^D_s) ds \right)  f(Z^D_t) \right], 
 \quad t>0, \; x \in D,
\end{equation}
be the Feynman-Kac transform of the semigroup corresponding to the killed process $Z^D$. 
By \cite[Theorem 2.15]{CKSV20}, the semigroup $(T^{D, \wt{\kappa}}_t)_{t\ge 0}$ has a transition density $p^{Z,D}(t,x,y)$ (with respect to the Lebesgue measure) such that for every $T>0$ and $b\in (0,1]$ there exists a constant $c_2=c_2
 (b, T, D) 
>0$ such that
\begin{equation}\label{e:pZD}
p^{Z,D}(t,x,y)\ge c_2\wt{q}(t,x,y), \quad t\in (0,T), \; x,y \in D \textrm{ such that }\delta_D(x)\wedge \delta_D(y)>bt^{1/\alpha}.
\end{equation}

Finally, by computing the Dirichlet forms of $Y^{\kappa, D}$ and $Z^D$ (for the latter use \cite[Theorem 6.1.2]{FOT}), we conclude that they coincide. This implies that $P^{\kappa,  D}_t=T^{D, \wt{\kappa}}_t$.   
Combining the lower bound in  \eqref{e:alpha-stable-estimate} with \eqref{e:pZD}, the proof is complete.   
\end{proof}

Recall that $\e_d=(\wt 0, 1)$. 
Define for $a\in (0,1]$, 
$$
S(a)=\{(\wt z, z_d) \in \R^d_+: |\wt z| < 2/a, \, a/(2a+2) < z_d < (2a+2)/a\}.
$$ 
Note that $B(\e_d, 1/8) \subset S(a)$ for all $a\in (0,1]$. Moreover, for any $a \in (0,1]$, $x\in B(\e_d,1/8)$ and  $y=(\wt y, y_d) \in \R^d_+$ with $x_d \wedge y_d > a|x-y|$, we have $|\wt y| \le |\wt x| + |y- x| <1/8 +  x_d/a<2/a$, $y_d + y_d/a \ge x_d -|x-y| + |x-y| \ge 7/8$ and $y_d \le x_d + |x-y| \le x_d + x_d/a <(2a+2)/a$.  Thus,  
\begin{align}\label{e:ncl}	
	&\big\{(t, (\wt x, x_d), (\wt y, y_d) ) \in \R^1_+ \times  B(\e_d, 1/8) \times  \R^d_+: x_d \wedge y_d > a(t^{1/\alpha} \vee |x-y|)\big\}\\
	& \subset (0, 2/a^\alpha)  \times S(a) \times  S(a) \qquad \text{for every} \;\, a\in (0,1].\nn
\end{align}

\begin{prop}\label{p:lower-heatkernel}
For any $a \in (0,1]$, there exists a constant $C =C(a)>0$ 
such that the following estimates hold:  For any $t>0$ and $x \in \R^d_+ \setminus \sN$, there is a  measurable set $N_{t,x} \subset \R^d_+$ 
of zero Lebesgue measure
such that for  all $y \in \R^d_+ \setminus N_{t,x}$ with $x_d \wedge y_d > a (t^{1/\alpha} \vee |x-y|)$,	\begin{equation}\label{e:lower-heatkernel}		\bar p(t,x,y) \ge p^\kappa(t,x,y) \ge C	
 p_\alpha(t, x, y).
 \end{equation}
\end{prop}
\begin{proof} 
By \eqref{e:rel-heatkernel},
it suffices to prove the second inequality in  \eqref{e:lower-heatkernel}. 
Since $\sN$ has 
zero Lebesgue measure, 
there exist $r>0$ and $\wt z \in \R^{d-1}$ such that $r^{-1}(x- (\wt z,0)) \in B(\e_d, 1/8) \setminus \sN$. By Lemmas \ref{l:scaling} and \ref{l:translation-invariance}, we have
\begin{align}\label{e:heatkernel-scaling}
	p^\kappa(t,x,y) = r^{-d}p^\kappa (r^{-\alpha}t, r^{-1}(x- (\wt z,0)), r^{-1}(y- (\wt z,0)) ) \qquad \text{for a.e.} \;\, y \in \R^d_+.
\end{align}
Therefore, by \eqref{e:ncl}, it is enough to prove that  there exists  a constant $c=c(a)>0$  such that  for all $0<t<2/a^\alpha$, 
$x \in B(\e_d, 1/8) \setminus \mathcal N$ 
and a.e. $y \in S(a) \setminus \sN$,
\begin{align}\label{e:nclaim}
p^{\,\kappa, S(a/2)}(t,x,y) \ge c 
p_\alpha(t, x, y),
\end{align}
which implies 
the second inequality in  \eqref{e:lower-heatkernel}
by \eqref{e:Dirichlet-kernel} and \eqref{e:heatkernel-scaling}.

Note that  $\delta_{S(a/2)}(y) \ge  a/(a^2+3a+2)$ for every $y \in S(a)$. Hence for all $0<t<2/a^\alpha$ and  $z,y \in S(a/2)$ it holds that $\delta_{S(a/2)}(z) \wedge \delta_{S(a/2)}(y) \ge (a^2/(2a^2+6a+4))t^{1/\alpha}$. 
Now it follows  from Proposition \ref{p:POTAl3.1} (with $b=(a^2/(2a^2+6a+4))$) that \eqref{e:nclaim} holds.
\end{proof}

As a direct consequence of Proposition \ref{p:lower-heatkernel},  for any $a \in (0,1]$, there exists $C>0$ such that 
\begin{equation}\label{e:ndl-1}
	\bar p(t,x,y) \ge p^\kappa(t,x,y) \ge C t^{-d/\alpha}
\end{equation}
for all $t>0$, $x\in \R^d_+ \setminus \sN $ and a.e. $y \in \R^d_+\setminus \sN$ with $x_d \wedge y_d > at^{1/\alpha}$ and $|x-y| \le t^{1/\alpha}$.

 By repeating the proofs of \cite[Lemmas 6.1 and 6.3]{CKSV}, we obtain the following two results from \eqref{e:Dirichlet-kernel},  Proposition \ref{p:upper-heatkernel} and  \eqref{e:ndl-1}.
\begin{lemma}\label{l:NDL}
	There exist constants $C>0$ and $\eta \in (0,1/4)$ such that 
	for all $x \in \R^d_+$, $r \in (0, x_d)$, $t \in (0, (\eta r)^\alpha]$ and $z \in B(x, \eta t^{1/\alpha}) \setminus \sN$,  
	\begin{equation*}
		\bar	p^{B(x,r)}(t,z,y)  \ge 	p^{\kappa, B(x,r)}(t,z,y) \ge C t^{-d/\alpha} \quad \text{for a.e. }  y \in B(x, \eta t^{1/\alpha}).
	\end{equation*}
\end{lemma}

\begin{lemma}\label{l:mean-exit-time}
	There exists a constant $C>1$ such that for all $x \in \R^d_+ \setminus \sN$ and $r \in (0, x_d)$,
	\begin{equation*}
		C^{-1}r^\alpha \le  \E_x[\tau^\kappa_{B(x,r)}] \le 
		\sup_{z \in B(x,r)\setminus \sN  } 
		\E_z[\bar\tau_{B(x,r)}] \le C r^\alpha.
	\end{equation*}
\end{lemma}

The L\'evy system formula (see \cite[Theorem 5.3.1]{FOT} and the arguments in \cite[p.40]{CK03}) 
states
that for any 
non-negative Borel function $F$ on $\R^d_+\times \R^d_+$ vanishing on the diagonal and any stopping time $T$ for $Y^\kappa$, it holds that
\begin{align}
\label{e:levysk}
\E_x\sum_{s\le T} F (Y^\kappa_{s-}, Y^\kappa_s)=
\E_x\int^T_0\int_{\R^d_+} F (Y^\kappa_s, y)J(Y^\kappa_s, y) dyds, \quad x\in \R^d_+\setminus \sN.
\end{align} 
Here $Y^\kappa_{s-}=\lim_{t\uparrow s}Y^\kappa_t$ denotes the left limit of the process $Y$ at time $s>0$.
Similarly, for any 
non-negative Borel function $F$ on $\oR^d_+\times \oR^d_+$ vanishing on the diagonal and any stopping time $T$ for $\rY$, it holds that
\begin{align}
\label{e:levysb}
\E_x\sum_{s\le T} F (\rY_{s-}, \rY_s)=
\E_x\int^T_0\int_{\R^d_+} F (\rY_s, y)J(\rY_s, y) dyds, \quad x\in \oR^d_+\setminus \sN.
\end{align} 
See \cite[(3.3) and (3.4)]{KSV-Green} for a simper form following from \eqref{e:levysk}, which will be used in this paper too.

For $x=(\wt x, x_d) \in \oR^d_+$ and $t>0$,  we define
\begin{align}\label{e:def-VW}
	 V_x(t)&=\{z=(\wt{z}, z_d)\in \R^d:\, |\wt{z}-\wt{x}|<2t^{1/\alpha},\, z_d \in [0, 2t^{1/\alpha})\},\\
	  W_x(t)&=\{z=(\wt{z}, z_d)\in \R^d:\, |\wt{z}-\wt{x}|<2t^{1/\alpha},\, z_d \in [x_d+5t^{1/\alpha},x_d + 8t^{1/\alpha})\}.\nn
\end{align}
In dimension 1, we abuse notation and use  $V_x(t)=[0,2t^{1/\alpha})$ and 
 $W_x(t)=[x+5t^{1/\alpha},x + 8t^{1/\alpha})$. 

\begin{lemma}\label{l:barY}
	(i) 
	There exists $C>0$ such that for all $t>0$ and  $x \in \oR^d_+ \setminus \sN$
	with $x_d\le t^{1/\alpha}$, 
	$$
	\P_x\big(\rY_{\overline\tau_{V_x(t)}} \in W_x(t)\big) \ge C.
	$$
	
\noindent	(ii) There exists $C>0$ such that for all $n \ge 1$, $t>0$ and $x \in \oR^d_+ \setminus \sN$ with $x_d\le t^{1/\alpha}$,
	$$
	\P_x(\overline\tau_{V_x(t)}>nt) \le 2e^{-Cn}.
	$$ 
\end{lemma}
\begin{proof} (i) Define $V_x(t,r)=\{z \in V_x(t):\delta_{V_x(t)}(z)>r t^{1/\alpha}\}$ for $r>0$. For any 
$r\in(0,1)$, $z \in V_x(t,r)$, $u \in B(z, r t^{1/\alpha}/2)$ and $w \in W_x(t)$, we have
\begin{align}\label{e:ndl-bar-1}
   r t^{1/\alpha}/2 \le u_d \le  5t^{1/\alpha} \le w_d    \,\,\, \text{and} \,\,\,	|u-w| \le |u-z| + |z-x| + |x-w| < 16t^{1/\alpha}.  
\end{align}
Using the  L\'evy system formula in \eqref{e:levysb},  {\bf (A3)}(I), \eqref{e:ndl-bar-1} and Lemma \ref{l:mean-exit-time},  we get that for any  $r\in(0,1)$ and  $z \in V_x(t,r) \setminus \sN$,
\begin{align}\label{e:barY-1}
	&\P_z\big(\rY_{\overline\tau_{V_x(t)}} \in W_x(t)\big)\ge \P_z\big(\rY_{\overline\tau_{B(z,r t^{1/\alpha}/2)}} \in W_x(t)\big)  \\
	&= \E_z \bigg[\int_0^{ \overline\tau_{B(z,r t^{1/\alpha}/2)}} \int_{W_x(t)} J(\rY_s, w)dw ds \bigg] \nn \\
	&\ge c_1 \E_z \bigg[\int_0^{ \overline\tau_{B(z,r t^{1/\alpha}/2)}} \int_{W_x(t)} \left( \frac{r/2}{16} \right)^{\beta_1}  \left( \frac{5}{16} \right)^{\beta_2} (16t^{1/\alpha})^{-d-\alpha}  dw ds \bigg] \nn\\
	&\ge c_2 r^{\beta_1} t^{-1-d/\alpha}\E_z \big[ \overline\tau_{B(z,r t^{1/\alpha}/2)} \big]\int_{W_x(t)} dw \ge c_3 
r^{\beta_1+\alpha}. \nn
\end{align}

By  Proposition \ref{p:upper-heatkernel} and \cite[Remark 3.3]{GH}, the condition (i) in \cite[Theorem 3.1]{GH} holds true with $\rho(s)=s^{1/\alpha}$. Then, since $\rY$ is conservative, by the implication (i) $\Rightarrow$ (ii) of \cite[Theorem 3.1]{GH}, there exists a constant $\eps_0>0$ independent of $t$ and $x$ such that
\begin{align}\label{e:barY-2}
	\P_x(\overline\tau_{V_x}(t)>\eps_0 t) \ge 1/2. 
\end{align} On the other hand, for all $r>0$, we get from  Proposition \ref{p:upper-heatkernel}  that
\begin{align}\label{e:barY-3}
         &\P_x(\overline\tau_{V_x(t) \setminus V_x(t,r)}>\eps_0t)  \le \P_x(\rY_{\eps_0 t} \in V_x(t) \setminus V_x(t,r)) \\
	&\le \int_{V_x(t) \setminus V_x(t,r)} p(\eps_0 t,x,y)dy \le  c_4t^{-d/\alpha} \int_{V_x(t) \setminus V_x(t,r)}  dy \le c_5 r. \nn
\end{align}
Set $r_0:=1/(4c_5+1)$. Using the strong Markov property and   \eqref{e:barY-1}--\eqref{e:barY-3}, we obtain
\begin{align*}
         &\P_x\big(\rY_{\overline\tau_{V_x(t)}} \in W_x(t)\big) \ge 	\P_x\big(\rY_{\overline\tau_{V_x(t)}} \in W_x(t), \, \overline\tau_{V_x(t) \setminus V_x(t,r_0)}\le \eps_0 t<\overline\tau_{V_x(t)}\big)\\
	&\ge \P_x\Big( \P_{\rY_{\overline\tau_{V_x(t) \setminus V_x(t,r_0)}}} \big(\rY_{\overline\tau_{V_x(t)}} \in W_x(t)\big) \, : \, \overline\tau_{V_x(t) \setminus V_x(t,r_0)}\le \eps_0 t<\overline\tau_{V_x(t)}\Big)\\
	&\ge c_3
 r_0^{\beta_1+\alpha}	\P_x\left( \overline\tau_{V_x(t) \setminus V_x(t,r_0)}\le \eps_0 t<\overline\tau_{V_x(t)} \right)\\
	&\ge c_3
	r_0^{\beta_1+\alpha} \left(\P_x\left(\overline\tau_{V_x(t)}>\eps_0 t \right)-\P_x\left( \overline\tau_{V_x(t) \setminus V_x(t,r_0)}> \eps_0 t\right) \right) \ge c_3 	r_0^{\beta_1+\alpha}/4. 
\end{align*}

(ii) By Proposition \ref{p:upper-heatkernel}, there exists a constant $k_0>0$ independent of $t$ and $x$ such that  for any $z \in V_x(t) \setminus \sN$,  
\begin{align}\label{e:SP-1}
	\P_z\big(\overline\tau_{V_x(t)}>k_0t\big)  \le \P_z\big(\rY_{k_0t} \in V_x(t)\big) \le  \int_{V_x(t)} \bar p(k_0t,z,y)dy \le
	\frac{c_1t^{d/\alpha}}{(k_0t)^{d/\alpha}} \le \frac12.
\end{align}
For $r>0$, let $\lfloor r \rfloor:=\sup\{m \in \N: m \le r\}$. Now for any $n \ge1$, using the Markov property and \eqref{e:SP-1}, we get that
\begin{align*}
         &\P_x\left(\overline\tau_{V_x(t)}>nt\right) \le 	\P_x\left(\overline\tau_{V_x(t)}> \lfloor n/k_0 \rfloor k_0t\right) \\
	&\le \P_x\left(\P_{X_{(\lfloor n/k_0 \rfloor-1)k_0t}}(\overline\tau_{V_x(t)}>k_0t) : \overline\tau_{V_x(t)}>(\lfloor n/k_0 \rfloor-1)k_0t\right)  \\
	&\le 2^{-1}\P_x\big(\overline\tau_{V_x(t)}>(\lfloor n/k_0 \rfloor-1)k_0t\big) \le \cdots \le 2^{-\lfloor n/k_0 \rfloor}  \le 2e^{- (\log 2)n/k_0}.
\end{align*}
\end{proof}
\begin{lemma}\label{l:bar-lower}
 	There exist constants $M>1$ and $C>0$ such that 
	$$	
	\essinf_{z \in W_x(t) } \,	\bar p(Mt,x,z) \ge Ct^{-d/\alpha} \quad \text{for all	$t>0$ and $x \in \oR^d_+ \setminus \sN$}.
	$$
\end{lemma}
\begin{proof} 
 Suppose $x_d>t^{1/\alpha}$. For any $M>1$ and  $z \in W_x(t)$, since $|x-z|<10t^{1/\alpha}$, we see that
 \begin{align*}
 	x_d \wedge z_d \ge t^{1/\alpha} \ge \frac{((Mt)^{1/\alpha}\vee |x-z|)}{M^{1/\alpha} \vee 10} \quad \text{and} \quad (Mt)^{-d/\alpha} \wedge \frac{t}{|x-z|^{d+\alpha}} \ge  \frac{t^{-d/\alpha}}{M^{d/\alpha} \vee 10^{d+\alpha}}.
 \end{align*}
 Thus, the result follows from Proposition \ref{p:lower-heatkernel}. 
 
 Suppose $x_d\le t^{1/\alpha}$. For all $w,z \in W_x(t)$, we have $w_d \wedge z_d \ge 5t^{1/\alpha}$ and  $|w-z|<7t^{1/\alpha}$. Therefore, by Proposition \ref{p:lower-heatkernel}, for any $M>1$, there is $c_1=c_1(M)>0$ independent of $t$ and $x$ such that 
 for all $w \in W_x(t) \setminus \sN$, $z \in W_x(t)$ and $\eps \in (0,t^{1/\alpha})$,
 \begin{align}\label{e:lbd1}
 \inf_{t \le s \le Mt}\,\P_w\big(\rY_{s} \in B(z,\eps)\big)  = \inf_{t \le s \le Mt} \int_{B(z,\eps)} \bar p(s,w,y)dy \ge c_1 t^{-d/\alpha} m_d(B(z,\eps)).
 \end{align}
 By the strong Markov property and \eqref{e:lbd1}, for all $M>1$,  $z \in W_x(t)$ and $\eps \in (0, t^{1/\alpha})$, 
 \begin{align}\label{e:interior-lower-1}
 	&\P_x\left(\rY_{Mt} \in B(z,\eps)\right)	 \\
 	&\ge  	\E_x\left[  \P_{\rY_{\overline\tau_{V_x(t)}}}\big(\rY_{Mt-\overline \tau_{V_x(t)}} \in B(z,\eps)\big)  : \overline\tau_{V_x(t)} \le (M-1)t,\, \rY_{\overline\tau_{V_x(t)}} \in W_x(t)\right] \nn \\
 	&\ge \Big(\inf_{t \le s \le Mt}\,\inf_{w \in W_x(t) \setminus \sN}\P_w\big(\rY_{s} \in B(z,\eps)\big)  \Big)
 	\,	\P_x \left(\overline\tau_{V_x(t)} \le (M-1)t, \,\rY_{\overline\tau_{V_x(t)}} \in W_x(t)\right) \nn\\
 	&\ge  c_1 t^{-d/\alpha} m_d(B(z,\eps))
 	\left(\P_x \big( \rY_{\overline\tau_{V_x(t)}} \in W_x(t)\big) -\P_x \big(\overline\tau_{V_x(t)} > (M-1)t\big) \right). \nn
 \end{align}
By Lemma \ref{l:barY}(i)-(ii), there are constants $c_2,c_3,c_4>0$ such that for all $M>1$,
\begin{align*}
	\P_x \big( \rY_{\overline\tau_{V_x(t)}} \in W_x(t)\big) -\P_x \big(\overline\tau_{V_x(t)} > (M-1)t\big) \ge c_2 - c_3 e^{-c_4M}.
\end{align*} 
Choosing $M=c_4^{-1}\log(2c_3/c_2)+1$,   we arrive at the result by \eqref{e:interior-lower-1} 
and the Lebesgue differentiation theorem. \end{proof}

\begin{lemma}\label{l:NDL-bar-1}
   There exists a constant $C>0$ such that
$\bar p(t,x,y) \ge C t^{-d/\alpha}$ for all $t>0$    and $x,y \in \oR^d_+ \setminus \sN$ with $|x-y|\le t^{1/\alpha}$.
\end{lemma}
\begin{proof}  Let $M> 1$ be the constant in  
Lemma \ref{l:bar-lower}. 
Note that for all $(z,w) \in W_x(t/M) \times W_y(t/M)$, we have $z_d \wedge w_d \ge 5(t/M)^{1/\alpha}$ and  $|z-w| \le |z-x| + |x-y| + |y-w|  (20+M^{1/\alpha})(t/M)^{1/\alpha}$. 
Therefore, by Proposition \ref{p:lower-heatkernel}, there is $c_1>0$ independent of $t,x,y$ such that 
\begin{align}\label{e:NDL-bar}
\essinf_{z \in W_x(t/(3M))} \essinf_{w\in  W_y(t/(3M))} \bar p(t/3,z,w) \ge c_1t^{-d/\alpha}.
\end{align}
By the semigroup property, \eqref{e:NDL-bar} and Lemma \ref{l:bar-lower},
\begin{align*}
	\bar p(t,x,y)&\ge
	\int_{W_x(t/(3M))}\int_{W_y(t/(3M))}\bar p(t/3,x,z)
\bar	p(t/3,z,w) \bar p(t/3,w,y)
	dzdw\nn\\
	& \ge   
	 \left(\essinf_{z \in W_x(t/(3M)) } \bar p(t/3,x,z) \right) 
	\left(\essinf_{w \in W_y(t/(3M)) } \bar p(t/3,y,w) \right)
	 \nn\\
	&\quad\;\; \times 
	\left( \essinf_{z \in W_x(t/(3M))} \essinf_{w\in  W_y(t/(3M))} \bar p(t/3,z,w) \right)
	\int_{W_x(t/(3M))}\int_{W_y(t/(3M))}
	dzdw \nn\\ 
	&\ge c_2 t^{-3d/\alpha+2d/\alpha} = c_2t^{-d/\alpha}.
\end{align*}
The proof is complete.
\end{proof}

For $x \in \oR^d_+$ and $r>0$, we denote $\uB(x,r):=B(x,r) \cap \oR^d_+$.
 We observe that
\begin{align}\label{e:d-set}
	m_d(\uB(x,r)) \asymp r^d \qquad \text{for all} \;\, x \in \oR^d_+, \, r>0.
\end{align}

Now, using  \eqref{e:upper-heatkernel} and Lemma \ref{l:NDL-bar-1} (instead of \eqref{e:ndl-1}), we extend the results in  Lemmas \ref{l:NDL} and \ref{l:mean-exit-time} to $\overline Y$ removing the restrictions on $x$ and $r$.
\begin{lemma}\label{l:NDL-bar}
	There exist constants $C>0$ and $\eta \in (0,1/4)$ such that for all $x \in \oR^d_+ \setminus \sN$, $r>0$, $t \in (0, (\eta r)^\alpha]$ and $z \in \uB(x, \eta t^{1/\alpha}) \setminus \sN$,
	\begin{equation*}
		\bar	p^{\uB(x,r)}(t,z,y)  \ge  C t^{-d/\alpha} \quad \text{for a.e. }  y \in \uB(x, \eta t^{1/\alpha}) \setminus \sN.
	\end{equation*}
\end{lemma}

\begin{lemma}\label{l:mean-bar}
	There exists a constant $C>1$ such that for all 
	$x \in \oR^d_+ \setminus \sN$
	and $r>0$,
	\begin{equation}\label{e:mean-bar}
		C^{-1}r^\alpha \le  \E_x\left[\overline\tau_{\uB(x,r)}\right] \le 
		\sup_{z \in \uB(x,r)\setminus \sN} 
		\E_z\left[\overline\tau_{\uB(x,r)}\right] \le C r^\alpha.
	\end{equation}
\end{lemma}

Let $\bar X := (T_s, \rY_s)_{s\ge 0}$ and  $X^\kappa := (T_s, Y^\kappa_s)_{s\ge 0}$ be  time-space processes  where $T_s = T_0 -s$. 
The law of the time-space process $s\mapsto \bar X_s$ or $s\mapsto X^\kappa_s$ starting from $(t, x)$ will be denoted by 
$\P_{(t,x)}$. 
For
every open subset $U$ of $[0,\infty) \times \R^d$, define $\bar\tau_{U} = \inf\{s > 0:\, \bar X_s\notin U\}$ and $\tau^{\kappa}_U = \inf\{s > 0:\,  X_s^\kappa \notin U\}$. 
We also define $\bar\sigma_{U} = \inf\{s > 0:\, \bar X_s\in U\}$ and $\sigma^{\kappa}_U = \inf\{s > 0:\,  X_s^\kappa \in U\}$.

A Borel  function $u:[0,\infty)\times \overline\R^d_+\to \R$ is 
said to be parabolic in $
 (a,b] 
\times \uB(x,r)\subset (0, \infty)\times \oR^d_+$ 
with respect to $\rY$ if for every relatively compact open set $U\subset 
 (a,b] \times \uB(x,r)$ with respect to the  topology on $[0,\infty) \times \oR^d_+$,
 it holds that 
 $u(t,z)=\E_{(t,z)}u(\bar X_{\bar\tau_U})$ for all 
$(t,z)\in U$ with $z\notin \sN$.
Similarly, a Borel  function $u:[0,\infty)\times \R^d_+\to \R$ is 
said to be parabolic in $
 (a,b] \times B(x,r)\subset (0, \infty)\times \R^d_+$ 
with respect to $Y^\kappa$ if for every relatively compact open set $U\subset 
 (a,b] \times B(x,r)$, it holds that 
$u(t,z)=\E_{(t,z)}u(X^\kappa_{\tau^\kappa_U})$ for all 
$(t,z)\in U$ with $z\notin \sN$.

\begin{lemma}\label{l:ckw-3-7} 
\noindent (i) Let $\eta\in (0,1/4)$ be the constant from Lemma \ref{l:NDL}. 
For every $\delta\in (0,\eta]$, there exists a constant $C>0$  such that for all $x \in \R^d_+\setminus \sN$, $r \in (0, x_d)$,
 $t\ge \delta r^\alpha$, and any compact set $A\subset [t-\delta r^\alpha, t-\delta r^\alpha/2]\times B(x, (\eta \delta/2)^{1/\alpha} r)$, 
\begin{equation*}
 \P_{(t,x)}\big(
\sigma^\kappa_A < \tau^\kappa_{[t-\delta r^\alpha, t]\times B(x,r)}\big)
\ge C \frac{m_{d+1}(A)}{ r^{d+\alpha}}. 
\end{equation*}
(ii)	Let $\eta\in (0,1/4)$ be the constant from Lemma \ref{l:NDL-bar}. 
For every $\delta\in (0,\eta]$, there exists a constant $C>0$  such that for all $x \in \oR^d_+\setminus \sN$, $r>0$,
$t\ge \delta r^\alpha$, and any compact set $A\subset [t-\delta r^\alpha, t-\delta r^\alpha/2]\times \uB(x, (\eta \delta/2)^{1/\alpha} r)$, 
\begin{equation*}
 	\P_{(t,x)}\left(
	\bar\sigma_A < \bar\tau_{[t-\delta r^\alpha, t]\times \uB(x,r)}\right) 
	\ge C \frac{m_{d+1}(A)}{ r^{d+\alpha}}. 
\end{equation*}
\end{lemma}
\begin{proof}   By repeating the proofs of \cite[Lemma 6.5]{CKSV}
(using the L\'evy system formulas in \eqref{e:levysk} and \eqref{e:levysb}), 
we deduce the results from  Lemmas  \ref{l:NDL} and  \ref{l:NDL-bar} respectively.
\end{proof}

\begin{thm}\label{t:phr}
(i) For any $\delta \in (0,1)$, 
there exist $ \lambda \in (0,1]$ and $C>0$ such that
for all  $x \in \R^d_+$,  
$r \in (0, x_d)$,  $t_0\ge 0$, and any  function $u$ on $(0,\infty)\times \R^d_+$ which is parabolic in $(t_0, 
 t_0+ r^\alpha]\times B(x, r)$ with respect to $Y^\kappa$ and bounded in $(t_0, 
 t_0+ r^\alpha]
\times \R^d_+$, we have
\begin{equation}\label{e:phr}
	|u(s,y) - u(t,z)| \le C \bigg(\frac{|s-t|^{1/\alpha} + |y-z|}{r}\bigg)^{\lambda } \esssup_{[t_0, t_0 + r^\alpha] \times \R^d_+}|u|,
\end{equation}
for every $s,t \in (t_0 + (1-\delta^\alpha)r^\alpha, 
 t_0+ r^\alpha]$ and $y,z \in B(x, \delta r)  \setminus \sN $.

\noindent (ii) For any $\delta \in (0,1)$, 
there exist $\lambda   \in (0,1]$ and $C>0$ such that
for all 
$x \in \oR^d_+$, $r>0$,  $t_0\ge 0$, and any  function $u$ on $(0,\infty)\times \oR^d_+$ which is parabolic in $(t_0, 
 t_0+ r^\alpha]\times \uB(x, r)$ with respect to $\rY$ and bounded in $(t_0, 
 t_0+ r^\alpha]\times \overline\R^d_+$, \eqref{e:phr} holds true 
for every $s,t \in (t_0 + (1-\delta^\alpha)r^\alpha, t_0+ r^\alpha]$ and 
$y,z \in \uB(x, \delta r)\setminus \sN$.
\end{thm}
\begin{proof}  Using Lemmas  
\ref{l:NDL},  \ref{l:mean-exit-time} and \ref{l:ckw-3-7}(i) for (i), and \eqref{e:d-set} and  Lemmas  
\ref{l:NDL-bar},  \ref{l:mean-bar} and \ref{l:ckw-3-7}(ii) for (ii), we get the desired results
using the same argument as in the proof of \cite[Theorem 4.14]{CK03} (see also the proof of \cite[Proposition 3.8]{CKW20}). We omit details here. 
\end{proof}

\begin{remark} \label{r:conti}
 {\rm By Theorem  \ref{t:phr}, since 
the heat kernels $\bar p(t,x,y)$ and $p^\kappa(t,x,y)$ are 
parabolic with respect to  $\rY$ and $Y^\kappa$ respectively, they can be extended continuously to
 $(0,\infty) \times \oR^d_+ \times \oR^d_+$  and  $(0,\infty) \times \R^d_+ \times \R^d_+$ respectively. 
As consequences, by Proposition \ref{p:upper-heatkernel},  $\rY$ and $Y^\kappa$ can be refined to be 
a strongly Feller processes starting from every point in $\oR^d_+$ and $\R^d_+$ respectively, and  the exceptional set $\sN$ in Proposition \ref{p:upper-heatkernel} can be taken to be  the empty set.

In the remainder of this paper, 
we  take the jointly continuous version 
of $\bar p:(0,\infty) \times \oR^d_+ \times \oR^d_+ \to [0,\infty)$
and $p^\kappa:(0,\infty) \times \R^d_+ \times \R^d_+ \to [0,\infty)$, 
take the exceptional set $\sN$ in Proposition \ref{p:upper-heatkernel} to be empty set, 
and replace the ${\rm essinf}$ in Lemma \ref{l:bar-lower} by $\inf$.
Again, for notational convenience, we extend
the domain of $p^\kappa(t,x,y)$ 
 to  $(0,\infty) \times \oR^d_+ \times \oR^d_+ $  by letting $p^\kappa(t,x,y)=0$ if 
 $x \in \partial  \R^d_+$ or $y\in \partial  \R^d_+$.

The following scaling and horizontal translation invariance properties of the heat kernels  (the latter for $d\ge 2$), which come from   Lemmas \ref{l:scaling} and \ref{l:translation-invariance},  will be used throughout the paper:	For any  $(t,x,y) \in (0,\infty) \times \oR^d_+ \times \oR^d_+$, $r>0$ and $\wt z \in \R^{d-1}$,
 \begin{equation*}
		\bar p(t,x,y)  =r^{-d} \bar p(t/r^\alpha, x/r, y/r) = \bar p(t, x + (\wt z, 0), y + (\wt z,0)),
\end{equation*}
\begin{equation}\label{e:kernel-scaling}
	 p^\kappa(t,x,y)  =r^{-d} p^\kappa(t/r^\alpha, x/r, y/r) = p^\kappa(t, x + (\wt z, 0), y + (\wt z,0)).
\end{equation}
}
	\end{remark}

\medskip

From Proposition \ref{p:upper-heatkernel}, since the exceptional set $\sN$ is removed, we obtain

\begin{cor}\label{c:green-upper}	If $d >\alpha$, 	\label{c:Green}	$$ G^\kappa(x,y) \le \bar G(x,y) \le \frac{c}{|x-y|^{d-\alpha}}, \quad x, y \in\oR_+^d.$$\end{cor}	

\begin{remark} \label{r:dalpha}	
	{\rm The assumption $d>(\alpha+\beta_1+\beta_2) \wedge 2$ in \cite{KSV-Green, KSV-nokill} is only used to show   		$G^\kappa(x,y) \le c|x-y|^{-d+\alpha}$. Thus, by Corollary \ref{c:Green}, all results in \cite{KSV-Green, KSV-nokill}	with the assumption $d>(\alpha+\beta_1+\beta_2) \wedge 2$		hold under the weaker 		assumption $d>\alpha$.
	Note that there is a typo in 
	the statement of \cite[Theorem 1.3]{KSV-Green}: The assumption  $d>\alpha+\beta_1+\beta_2$   in \cite[Theorem 1.3]{KSV-Green} should be $d>(\alpha+\beta_1+\beta_2)\wedge 2$. 
		}
\end{remark}

Since the heat kernels $\bar p(t,x,y)$ and $p^\kappa(t,x,y)$ are jointly continuous, we get the next lemma from the strong Markov properties of $\rY$ and $Y^\kappa$. 
 This lemma is a refined version of  
\cite[Lemma 4.2]{CKS14} which was inspired by \cite{siu}.  
In this paper, this lemma will play an important role in the bootstrap method to prove sharp upper estimates on the heat kernels. 
 Although the proof of next lemma is standard, we give it for  the reader's  convenience.

\begin{lemma}\label{l:general-upper-2}	
	(i)	Let  $V_1$ and $V_3$ be open subsets of $\oR^d_+$ with ${\rm dist}(V_1,V_3)>0$. Set $V_2:=\oR^d_+ \setminus (V_1 \cup V_3)$. For any $x\in V_1$, $y \in V_3$ and $t>0$, it holds that 	\begin{align*}	\bar p(t,x,y) &\le  \P_x(\overline \tau_{V_1}<t) \sup_{s \le t, \, z \in V_2} \bar p(s,z,y)\\		&\quad  + {\rm dist}(V_1,V_3)^{-d-\alpha}\int_0^t\int_{V_3} \int_{V_1} \bar p^{V_1}(t-s, x, u) \sB(u,w) \bar p(s, y,w) du dw ds.	\end{align*}
	(ii)	Let  $V_1$ and $V_3$ be open subsets of $\R^d_+$ with ${\rm dist}(V_1,V_3)>0$. Set $V_2:=\R^d_+ \setminus (V_1 \cup V_3)$. For any $x\in V_1$, $y \in V_3$ and $t>0$, it holds that 	\begin{align*}	p^\kappa (t,x,y) &\le  \P_x(\tau^\kappa_{V_1}<t< \zeta^\kappa) \sup_{s \le t, \, z \in V_2} p^\kappa(s,z,y)\\		&\quad  + {\rm dist}(V_1,V_3)^{-d-\alpha}\int_0^t\int_{V_3} \int_{V_1} p^{\kappa, V_1}(t-s, x, u) \sB(u,w) p^\kappa(s, y,w) du dw ds.	\end{align*}
\end{lemma} 
\begin{proof} Since the proofs are the same, we only give the proof of (ii).

By the strong Markov property,  the L\'evy system formula in \eqref{e:levysk} and symmetry, 
we get that for any $x \in V_1$ and $y \in V_3$, 
\begin{align*}	
p^\kappa(t,x,y)	
&=\E_x \left[p^\kappa(t-\tau^\kappa_{V_1}, Y^\kappa_{\tau^\kappa_{V_1}}, y) : \tau^\kappa_{V_1}<t<\zeta^\kappa,  Y^\kappa_{\tau^\kappa_{V_1}} \in V_2 \right]\\	&\quad +\E_x \left[p^\kappa(t-\tau^\kappa_{V_1}, Y^\kappa_{\tau^\kappa_{V_1}}, y) : \tau^\kappa_{V_1}<t<\zeta^\kappa,  Y^\kappa_{\tau^\kappa_{V_1}} \in V_3 \right]\\	& \le  \P_x ( \tau^\kappa_{V_1}<t<\zeta^\kappa) \sup_{s \le t, \, z \in V_2} p^\kappa(s,z,y)\\	&\quad  + \int_0^t \int_{V_3}\int_{V_1} p^{\kappa, V_1}(t-s, x, u) \frac{\sB(u,w)}{|u-w|^{d+\alpha}} p^\kappa(s, w,y)  dudwds\\	& \le \P_x ( \tau^\kappa_{V_1}<t<\zeta^\kappa) \sup_{s \le t, \, z \in V_2} p^\kappa(s,z,y)\\	& \quad  +{\rm dist}(V_1,V_3)^{-d-\alpha}\int_0^t\int_{V_3} \int_{V_1} 	p^{\kappa, V_1}(t-s, x, u) 	\sB(u,w) p^\kappa(s, y,w) du dw ds.
\end{align*}
\end{proof}

\section{Parabolic Harnack inequality and preliminary lower bound of heat kernels}\label{s:PHI}

In this section, we prove that the parabolic Harnack inequality holds for 
$\rY$ and $Y^\kappa$, and get some preliminary lower bounds for 
the heat kernels  $\bar p(t,x,y)$ and $p^\kappa(t,x,y)$.

Recall that $\uB(x,r)=B(x,r) \cap \oR^d_+$.
\begin{lemma}\label{l:ckw-4-2}
(i)	Let $\eta \in (0,1/4)$ be the constant from Lemma \ref{l:NDL} and let $\delta\in (0,2^{-\alpha-2}\eta)$. There exists a constant $C>0$ such that  for all $y \in \R^d_+$, $R \in (0, y_d)$, $r \in  (0, (\eta\delta /2)^{1/\alpha}R/2]$,
$\delta R^\alpha/2\le t-s \le 4\delta (2R)^\alpha$, 
$x\in B(y, (\eta\delta /2)^{1/\alpha} R/2  )$  and $z\in B(x_0, (\eta\delta /2)^{1/\alpha} R )$, 
	\begin{equation*}
	\P_{(t,z)}
	\big(\sigma^\kappa_{ \{s\}\times B(x,r)}\le \tau^\kappa_{[s,t]\times B(y,R)}\big)\ge C(r/R)^d. 
	\end{equation*}
	(ii)	Let $\eta \in (0,1/4)$ be the constant from Lemma \ref{l:NDL-bar} and let $\delta\in (0, 2^{-\alpha-2}\eta)$. There exists a constant $C>0$ such that  for all $y \in \oR^d_+$, $R>0$, $r \in  (0, (\eta\delta /2)^{1/\alpha}R/2]$,
$\delta R^\alpha/2\le t-s \le 4\delta (2R)^\alpha$, 
$x\in \uB(y, (\eta\delta /2)^{1/\alpha} R/2  )$  and $z\in \uB(x_0, (\eta\delta /2)^{1/\alpha} R )$, 
\begin{equation*}
	\P_{(t,z)}\big(\bar\sigma_{ \{s\}\times \uB(x,r)}\le \bar\tau_{[s,t]\times \uB(y,R)}\big) \ge C(r/R)^d. 
\end{equation*}
\end{lemma}
\begin{proof}  Using  Lemmas \ref{l:NDL} and  \ref{l:NDL-bar}, and \eqref{e:d-set}, the result can be proved by the  same argument as that of \cite[Lemma 6.7]{CKSV}. We omit details here. \end{proof} 

In order to obtain 
the parabolic Harnack inequality for $\rY$ and $Y^\kappa$, we introduce two conditions:

\smallskip

{\bf (UBS)} There exists a constant $C>0$ such that  for  all $x,y\in \oR_+^d$ and $0<r \le  |x-y|/2$,
\begin{align}\label{e:UBS}
\sB(x,y)\le \frac{C}{r^d}\int_{\uB(x,r)} \sB(z,y)dz.
\end{align}

{\bf (IUBS)}  
There exists a constant $C>0$ such that  \eqref{e:UBS} holds  for all $x,y\in \R_+^d$ and $0<r \le  (|x-y| \wedge x_d)/2$. (Note that, $\uB(x,r)=B(x,r)$ for this range of $r$.)

\smallskip

\begin{lemma}\label{l:check-UJS}
If  {\bf (A3)}{\rm (II)} also holds, then {\bf (UBS)}  is satisfied -- and thus {\bf (IUBS)} as well.
\end{lemma}
\begin{proof} Let $x,y\in \oR_+^d$ and $0<r \le  |x-y|/2$.
	For all $z \in \uB(x,r)$, $|x-y|/2 \le |z-y| \le 2|x-y|$ by the triangle inequality. Thus,  by {\bf (A3)}, there is $c_1>0$ independent of $x,y$ and $r$ such that for all $z \in \uB(x,r)$ with $z_d \ge x_d$, $\sB(x,y) \le c_1 \sB(z,y)$. Using this, we get
\begin{align*}
	\frac{1}{r^d}\int_{\uB(x,r)} \sB(z,y)dz \ge 
	\frac{\sB(x,y)}{c_1r^d}  \int_{B(x,r):z_d \ge x_d} dz \ge c_2 \sB(x,y).
\end{align*}
\end{proof}

We now show  that  the following parabolic Harnack inequalities  hold.
\begin{thm}\label{t:phi}
(i) Suppose that $\sB(x,y)$ 
satisfies {\bf (IUBS)}. 
Then there exist constants $\delta>0$ and $C,M\ge1$  such that for all $t_0\ge 0$, $x\in \R^d_+$ and $R\in (0, x_d)$, and any non-negative function $u$ on $(0,\infty)\times \R^d_+$ which is parabolic on 
$Q:=(t_0, t_0+4\delta  R^\alpha] \times B(x, R)$ 
with respect to $\bar Y$ or $Y^\kappa$, we have	
\begin{equation}\label{e:phi-1}\sup_{(t_1,y_1)\in Q_{-}}u(t_1,y_1)\le C \inf_{(t_2, y_2)\in Q_{+}}u(t_2, y_2),	
\end{equation}	
where $Q_-=[t_0+\delta R^\alpha, t_0+2\delta  R^\alpha]\times B(x, R/M)$ and $Q_+=[t_0+3\delta  R^\alpha, t_0+4\delta  R^\alpha]\times B(x, R/M)$.

\noindent(ii) Suppose that $\sB(x,y)$ satisfies 
{\bf (UBS)}. Then there exist constants $\delta>0$ and $C,M\ge1$  such that for all $t_0\ge 0$, $x\in \oR^d_+$ and $R>0$, and any non-negative function $u$ on $(0,\infty)\times \oR^d_+$ which is parabolic on 
$Q^0:=(t_0, t_0+4\delta  R^\alpha] \times \uB(x, R)$ with respect to $\bar Y$, we have	
\begin{equation}\label{e:phi-0}
	\sup_{(t_1,y_1)\in Q^0_{-}}u(t_1,y_1)\le C \inf_{(t_2, y_2)\in Q^0_{+}}u(t_2, y_2),	
\end{equation}	
where $Q^0_-=[t_0+\delta R^\alpha, t_0+2\delta  R^\alpha]\times \uB(x, R/M)$ and $Q^0_+=[t_0+3\delta  R^\alpha, t_0+4\delta  R^\alpha]\times \uB(x, R/M)$.
\end{thm}
\begin{proof}  (i) By  {\bf (IUBS)}, 
 there is a constant $C>0$ such that for all $x,y \in \R^d_+$ and $0<r \le (|x-y| \wedge x_d)/2$, 
\begin{align}\label{e:UJS}
 J(x, y)=	\frac{\sB(x,y)}{|x-y|^{d+\alpha}}\le \frac{C}{r^d}\int_{\uB(x,r)}\frac{\sB(z,y)}{|z-y|^{d+\alpha}}dz=
 \frac{C}{r^d}\int_{\uB(x,r)}J(z,y)dz. 
\end{align}
 Using Proposition \ref{p:upper-heatkernel}  and  \eqref{e:UJS}, one can follow the proof of \cite[Lemma 6.10]{CKSV} and see that 
\cite[Lemma 6.10]{CKSV}  is also 
 valid  for our case. (Note that, in the proof of  \cite[Lemma 6.10]{CKSV}, a pointwise comparison for the jump kernel from \cite[Proposition 6.8]{CKSV} was used to bound the term $I_2$ therein which can be replaced by \eqref{e:UJS}.) 
Using this and
 Lemmas \ref{l:NDL},   \ref{l:mean-exit-time},  \ref{l:ckw-3-7}, \ref{l:ckw-4-2},  the result can be proved 
using the same argument as in the proof of \cite[Lemma 5.3]{CKK09} (see also the proof of \cite[Lemma 4.1]{CKW20}). We omit details here. 

(ii) Since {\bf (UBS)} implies that \eqref{e:UJS} is satisfied for all $x,y \in \R^d_+$ and $0<r \le |x-y|/2$, using Proposition \ref{p:upper-heatkernel} and  \eqref{e:d-set}, one can also deduce that 
\cite[Lemma 6.10]{CKSV}  is 
 valid  for this case with $\uB(x_0,\cdot)$ instead of  $B(x_0, \cdot)$ in the definitions of $Q_i$, $1\le i\le 4$ therein. Then using Lemmas \ref{l:NDL-bar},   \ref{l:mean-bar},  \ref{l:ckw-3-7}(ii), \ref{l:ckw-4-2}(ii), one can follow the arguments in the  proof of \cite[Lemma 5.3]{CKK09} and conclude the result. \end{proof}

Using Lemma \ref{l:NDL}, and Theorem \ref{t:phi}(i), we obtain the following lemma.
\begin{lemma}\label{L:3.3}
Suppose that $\sB(x,y)$ 
satisfies {\bf (IUBS)}. For any  positive
 constants  $a, b$, there exists a constant
$C=C(a,b, \kappa)>0$ such that for all $z \in \bR_+^d$ and $r>0$ 
with $B(z, 2b r) \subset \R^d_+$,
$$
\inf_{y\in B(z, b r/2)}\P_y \big(\tau^\kappa_{B(z,b r)}> ar^\alpha  \big)
\, \ge\, C.
$$
\end{lemma}
\begin{proof}  By Lemma \ref{l:NDL}, there exist constants
$c_1,c_2,\eps_1>0$ such that for all $z\in \R^d_+$ and $r>0$ with $B(z,2b r) \subset \R^d_+$,
\begin{align}\label{e:L3.3}
\P_z \big( \tau^\kappa_{B(z, br/2 )} >\eps_1 r^\alpha \big) \ge
\int_{B(z, c_1r )} p^{\kappa, B(z, br/2 )}  (\eps_1 r^\alpha, z,w)dw \ge c_2. 
\end{align}
Thus it suffices 
to prove the lemma for $a>\eps_1$.  Applying the parabolic Harnack inequality (Theorem \ref{t:phi}) repeatedly, 
we conclude that there  
exists $c_3>0$ such that for any
 $w, y \in B(z, b r/2 )$,
 $$p^{\kappa, B(z, br )}(a r^\alpha ,y,w)\ge c_3\,p^{ \kappa, B(z, br )}(\eps_1 r^\alpha ,z,w). $$
Thus, using \eqref{e:L3.3}, we deduce that  for any  $y \in B(z, b r/2 )$,
\begin{align*}
\P_y \big( \tau^\kappa_{B(z, br)} > a r^\alpha \big)\ge c_3 \int_{B(z, br/2 )} p^{\kappa,B(z, br )}(\eps_1 r^\alpha ,z,w) dw 
\ge c_3 \P_z \big(\tau^\kappa_{B(z, br/2)}>\eps_1 r^\alpha\big) \ge  c_2c_3.
\end{align*}
This proves the lemma.
\end{proof}

Now, we follow the proof of \cite[Proposition 3.5]{CK} to get the following preliminary lower bound.
\begin{prop}\label{p:pld}
Suppose that $\sB(x,y)$ 
satisfies {\bf (IUBS)}. For any $a>0$,
 there exists a constant $C=C(a, \kappa)>0$ such 
 that for any $(t, x, y)\in (0, \infty)\times \R_+^d\times \R_+^d$ with $a t^{1/\alpha} \leq (4|x-y|)\wedge
 x_d\wedge y_d $, 
$$
\bar p (t,x,y) \ge p^\kappa (t, x, y)\ge C t
 J(x, y).
 $$ 
\end{prop}
\begin{proof}   The first inequality holds true by \eqref{e:rel-heatkernel} and Remark \ref{r:conti}.

For the second inequality, it suffices to consider the case $a\in (0,1]$. Let $a\in (0,1]$. Note  that for 
$ y_d \ge  a t^{1/\alpha}$ and $z \in B(y,5^{-1} a(t/2)^{1/\alpha})$, 
\begin{align}
	z_d \ge y_d-|z-y| \ge a (1-5^{-1}(1/2)^{1/\alpha})t^{1 /\alpha} \ge (4a/5)t^{1 /\alpha}.
	\label{e:nv61}
\end{align}
 By Lemma~\ref{L:3.3} and \eqref{e:nv61}, starting at $z\in B(y, 10^{-1} a (t/2)^{1/\alpha})$, with probability at least $c_1=c_1(a)>0$, the process $Y^\kappa$ does not move more than $10^{-1} a (t/2)^{1/\alpha} $ by time $t/2$.
Thus, using the strong Markov property and  the L\'evy system formula in \eqref{e:levysk}, we obtain
that for 
$a t^{1/\alpha} \leq (4|x-y|)\wedge
x_d\wedge y_d $, 
\begin{align}\label{e:nv1}
	\P_x \left( Y^\kappa_{t/2} \in B ( y, 5^{-1} a (t/2)^{1/\alpha}) \right)&\ge  c_1\P_x\big(\,Y^\kappa_{(t/2)\wedge 
		\tau^{\kappa}_{B(x,  20^{-1} a (t/2)^{1/\alpha})}}\in B(y,  10^{-1} a (t/2)^{1/\alpha})\\
		&\qquad \qquad \hbox{ and } (t/2) \wedge 
	\tau^{\kappa}_{B(x,  20^{-1} a  (t/2)^{1/\alpha})} \hbox{ is a jumping time\,}\big)\nonumber \\
	&= c_1 \E_x \bigg[\int_0^{(t/2)\wedge  \tau^{\kappa}_{B(x,  20^{-1} a (t/2)^{1/\alpha})}} \int_{B(y, 10^{-1} a (t/2)^{1/\alpha})}
	J(Y^\kappa_s, u) duds \bigg]. \nonumber
\end{align}
Observe that for any $a t^{1/\alpha} \leq (4|x-y|)\wedge
x_d\wedge y_d $ and $z\in B(x,20^{-1} a(t/2)^{1/\alpha})$, we have
\begin{align*}
 \frac{|y-z| \wedge y_d}2  \ge \frac{(|x-y|- 20^{-1} a(t/2)^{1/\alpha}) \wedge (at^{1/\alpha})}{2} \ge 10^{-1} a(t/2)^{1/\alpha}.
\end{align*}
Thus, by \eqref{e:UJS}, we obtain   that for  $a t^{1/\alpha} \leq (4|x-y|)\wedge
x_d\wedge y_d $,
\begin{align}\label{e:nv2}
	& \E_x \bigg[\int_0^{(t/2)\wedge
		\tau^{\kappa}_{B(x, 20^{-1} a (t/2)^{1/\alpha})}} \int_{B(y, 10^{-1} a (t/2)^{1/\alpha})}
	J(Y^\kappa_s,u) duds \bigg] \\
	&=  \E_x \bigg[\int_0^{t/2} \int_{B(y, 10^{-1} a (t/2)^{1/\alpha})}
	J(Y^{\kappa, B(x, 20^{-1} a (t/2)^{1/\alpha})}_s,u) duds \bigg] \nn\\
	&\ge  c_2 (10^{-1} a(t/2)^{1/\alpha})^d \int_0^{t/2} \E_x \left[
	J(Y^{\kappa, B(x, 20^{-1} a (t/2)^{1/\alpha})}_s,y) \right]ds \nn\\
	& \ge  c_3t^{d/\alpha} \int_{( a\eta /40)^\alpha t/2}^{(a\eta/20)^\alpha t/2} \int_{B(x, (a\eta^2/40) (t/2)^{1/\alpha})}
	J(w,y)\, p^{\kappa, B(x, 20^{-1} a (t/2)^{1/\alpha})}(s, x,w)     dw ds,
\end{align}
where $\eta\in (0,1/4)$ is the constant in Lemma \ref{l:NDL}. By Lemma \ref{l:NDL}, we see that for
$( a\eta /40)^\alpha t/2 <s <( a\eta /20)^\alpha t/2$ and $w\in B(x, (a\eta^2/40) (t/2)^{1/\alpha})$,
\begin{eqnarray}
p^{\kappa, B(x, 20^{-1} a (t/2)^{1/\alpha})}(s, x,w)  \ge  c_4 s^{-d/\alpha} \ge c_4 ((a \eta /20)^\alpha t/2)^{-d/\alpha}=:c_5t^{-d/\alpha}.
	\label{e:nv4}
\end{eqnarray}
Combining \eqref{e:nv1}, \eqref{e:nv2} with \eqref{e:nv4} and applying \eqref{e:UJS} again, we get   that for 
$a t^{1/\alpha} \leq (4|x-y|)\wedge
x_d\wedge y_d $, 
\begin{align}
	\label{e:nv6}
	\P_x \left( Y^\kappa_{t/2} \in B \big( y,   5^{-1} a (t/2)^{1/\alpha} \big) \right)&\ge  c_1c_3c_5\int_{( a\eta /40)^\alpha t/2}^{(a\eta/20)^\alpha t/2} \int_{B(x, (a\eta^2/40) (t/2)^{1/\alpha})}
	J(w,y)    dw ds\\
&	=  c_6 t  \int_{B(x, (a\eta^2/40) (t/2)^{1/\alpha})}
	J(w,y)    dw
	\ge   c_7   t^{1+d/\alpha}
	J(x,y).\nn
\end{align}
The proposition now follows from the Chapman-Kolmogorov equation along with  \eqref{e:nv6} and  Proposition \ref{p:lower-heatkernel} (using \eqref{e:nv61}):  for 
$a t^{1/\alpha} \leq (4|x-y|)\wedge
x_d\wedge y_d $, 
\begin{align*}
	p^\kappa(t, x, y) &= \int_{\R^d_+ } p^\kappa(t/2, x, z)p^\kappa(t/2, z, y)dz\ge \int_{B(y, \, 5^{-1}a (t/2)^{1/\alpha})}p^\kappa(t/2, x, z) p^\kappa(t/2, z, y) dz\\&\ge c_8 t^{-d/\alpha} \,\bP_x \left( Y^\kappa_{t/2} \in B(y, 5^{-1}a(t/2)^{1/\alpha} ) \right)\ge  c_9{ t} J(x,y).
\end{align*}
\end{proof}

 \section{Preliminary heat kernel upper bounds}\label{s:pupper}

The goal of this section is to prove the following  proposition. 
\begin{prop}\label{p:UHK-rough} Suppose that {\bf (A1)}--{\bf (A4)} and \eqref{e:killing-potential} hold 
with $q_\kappa\in [(\alpha-1)_+, \alpha+\beta_1)$. 
Then there exists a constant $C>0$ 	
such that for all $t>0, \, x,y \in \R^d_+$, 
	\begin{equation}\label{s:UHK-rough}		
	p^\kappa(t,x,y) \le C \left(1 \wedge \frac{x_d}{t^{1/\alpha}}\right)^{q_\kappa}\left(1 \wedge \frac{y_d}{t^{1/\alpha}} \right)^{q_\kappa} 
	p_\alpha(t, x, y).
	\end{equation}
\end{prop}

Note that, since 
$p_\alpha(t, x, y)$
is comparable to the transition density of 
the isotropic $\alpha$-stable process in $\R^d$, 
there exists a constant $C>0$ such that for all $t,s>0$ and  $x,y \in \R^d_+$,
\begin{align}\label{e:stableu1} 
	 \int_{\R^d_+} 
	 p_\alpha(t, x, z)dz
	 \le C, \end{align}
and
\begin{align} \label{e:stableu2} 
	\int_{\R^d_+} 
	p_\alpha(t, x, z)p_\alpha(s, z, y)
	dz \le C 
	 p_\alpha(t+s, x, y),
	 \end{align}
where in \eqref{e:stableu2} we used the semigroup property.

Before giving the proof of Proposition \ref{p:UHK-rough}, we record its simple consequence which directly follows from Proposition \ref{p:UHK-rough} and \eqref{e:stableu1}.

\begin{cor}\label{c:life}
	There exists a constant $C>0$ such that $$\P_x (\zeta^\kappa>t) \le  C \left(1 \wedge \frac{x_d}{t^{1/\alpha}}\right)^{q_\kappa}, \quad t>0, \; x \in \R^d_+.$$\end{cor}

When $q_\kappa=0$, Proposition \ref{p:UHK-rough} follows from Proposition \ref{p:upper-heatkernel}. Hence, 
in the remainder of this section,  
we  assume that {\bf (A1)}--{\bf (A4)}  hold, fix 
$\kappa \in [0,\infty)$  such that  $q_\kappa>0$ and denote by $q$ the constant $q_\kappa$ in \eqref{e:killing-potential}.
For the proof we will need several results from \cite{KSV-jump, KSV-Green, KSV-nokill} that we now recall for the convenience of the reader.

For $r>0$, we define a subset $U(r)$ of $\oR^d_+$, $d\ge 2$, by
$$
U(r):=\{x=(\wt{x}, x_d)\in \R^d:\, |\wt{x}|<r/2,\, 0\le x_d<r/2\}.
$$ 
When $d=1$, we abuse notation and use $U(r)=[0,r/2)$. For $t>0$ and an open set $V \subset \R^d_+$, denote by $Y_t^{\kappa,d}$ and $Y_t^{\kappa, V,d}$ the last coordinates of $Y_t^\kappa$  and $Y_t^{\kappa, V}$ respectively.

\begin{lemma}\label{l:l.5.2-gen}
(i) There exists 
	 $C>0$ such that for all $x,y \in \R^d_+$ satisfying $|x-y| \ge x_d$, 
	\begin{align*}
		\sB(x,y) \le C  x_d^{\beta_1}(|\log x_d|^{\beta_3}\vee 1)\big(1+{\bf 1}_{|y|\ge1}(\log|y|)^{\beta_3}\big)|x-y|^{-\beta_1}.
	\end{align*}
	
\noindent	(ii)  For any $R>0$, there exists 
$C>0$ such that for all $r\in (0, R]$ and $x \in U(2^{-4}r) \cap \R^d_+$,
	\begin{align*}
		\E_x \int_0^{\tau^\kappa_{U(r)}} (Y^{\kappa,d}_t)^{\beta_1} |\log Y^{\kappa, d}_t|^{\beta_3} dt \le 
		C x_d^q.
	\end{align*}

\noindent 	(iii) If  $q\in [(\alpha-1)_+, \alpha)$,
 then there exists  
  $C>0$ such that for all $r>0$ and $x \in U(2^{-4}r) \cap \R^d_+$,
$
	\P_x\big(\tau^\kappa_{U(r)} <\zeta^\kappa \big)= \P_x\big(Y^\kappa_{\tau^\kappa_{U(r)}}  \in \R^d_+\big) \le 
		C \left( {x_d}/{r} \right)^q.
$

\noindent 	(iv) If $q\in [(\alpha-1)_+, \alpha)$,  then there exists 
$C>0$ such that for all $r>0$ and $x \in U(2^{-4}r) \cap \R^d_+$,	
$	
\E_x[\tau^\kappa_{U(r)}] \le C \left( {x_d}/{r} \right)^q.
$
\end {lemma}

\begin{proof}
(i) See  \cite[Lemma 5.2(a)]{KSV-jump}. (ii) The result follows from scaling (Lemma \ref{l:scaling}), and \cite[Lemma 5.7(a)]{KSV-jump} if $\kappa>0$ and  \cite[Lemma 5.3]{KSV-nokill} if $\kappa=0$.
(iii) When $\kappa>0$, the  result follows from \cite[Lemma 3.4]{KSV-Green}. When $\kappa=0$, using \cite[Lemma 5.5]{KSV-nokill} and  
part (ii),
one can follow the proof of \cite[Lemma 3.4]{KSV-Green} and deduce the result. (iv) The result follows from scaling (Lemma \ref{l:scaling}), and \cite[Lemma 5.13]{KSV-jump} if $\kappa>0$ and  \cite[Lemma 4.5]{KSV-nokill} if $\kappa=0$.
\end{proof} 

\begin{prop}\label{p:Carleson}	
There exists a constant $C>0$ such that for any $w \in \partial \R^d_+$, $r>0$ and any non-negative function $f$ in $\R^d_+$ that is harmonic in $\R^d_+ \cap B(w,r)$ with respect to $Y^\kappa$ and vanishes continuously on $\partial \R^d_+ \cap B(w,r)$,  we have
	\begin{align*}
		f(x) \le C f(\wh x) \quad \text{for all} \;\, x \in \R^d_+ \cap B(w,r/2),
	\end{align*}
	where $\wh x \in \R^d_+ \cap B(w,r)$ with $\wh x_d \ge r/4$.
\end{prop}
\begin{proof}  The result follows from \cite[Theorem 1.2]{KSV-jump} if $\kappa>0$ and  \cite[Theorem 5.6]{KSV-nokill} if $\kappa=0$. \end{proof}

After recalling the known results above, we now continue with several auxiliary lemmas leading to the proof of Proposition \ref{p:UHK-rough}.

\begin{lemma}\label{l:lemma-upper-3}
	For all $\gamma\ge 0$,  $t>0$ and $x \in U(1) \cap \R^d_+$, it holds that
	\begin{equation}\label{e:lemma-upper-3-1}
		\int_{\R^d_+} p^\kappa(t,x,z)(1 \wedge z_d)^\gamma dz \le \E_x \left[ (1 \wedge Y^{\kappa, U(1),d}_{t})^\gamma : \tau^\kappa_{U(1)}>t\right] + \P_x\big(Y^\kappa_{\tau^\kappa_{U(1)}} \in \R^d_+\big).
	\end{equation}
	In particular, it holds that
	\begin{equation}\label{e:lemma-upper-3-2}
		\P_x(\zeta^\kappa>t) \le t^{-1}\E_x\big[\tau^\kappa_{U(1)}\big] + \P_x\big(Y^\kappa_{\tau^\kappa_{U(1)}}\in \R^d_+\big).
	\end{equation}
\end{lemma}
\begin{proof} Since $Y^{\kappa,U(1)}_t=Y^\kappa_t$ for $t<\tau^\kappa_{U(1)}$, we have
\begin{align*}
         &\int_{\R^d_+} p^\kappa(t,x,z)(1 \wedge z_d)^\gamma dz 
         = \E_x \left[ (1 \wedge Y^{\kappa, d}_{t})^\gamma : t<\zeta^\kappa \right]\\
	&=  \E_x \left[ (1 \wedge Y^{\kappa, d}_{t})^\gamma : \tau^\kappa_{U(1)}>t \right] 
	+ \E_x \left[ (1 \wedge Y^{\kappa, d}_{t})^\gamma : \tau^\kappa_{U(1)}\le t<\zeta^\kappa \right]\\
& \le  \E_x \left[ (1 \wedge Y^{\kappa, U(1), d}_{t})^\gamma : \tau^\kappa_{U(1)}>t \right] 
	+ \P_x \big(Y^\kappa_{\tau^\kappa_{U(1)}}\in \R^d_+\big).
\end{align*}
By taking $\gamma=0$ in \eqref{e:lemma-upper-3-1} and using Markov's inequality, we get 
\begin{equation*}
	\P_x(\zeta^\kappa>t)  \le \P_x(\tau^\kappa_{U(1)}>t) + \P_x\big(Y^\kappa_{\tau^\kappa_{U(1)}}\in \R^d_+\big)\le t^{-1}\E_x\big[\tau^\kappa_{U(1)}\big]+ \P_x \big(Y^\kappa_{\tau^\kappa_{U(1)}} \in \R^d_+\big).
\end{equation*}
\end{proof}

\begin{lemma}\label{l:general-upper-1}	There exists  $C>0$ such that 	
\begin{equation*}		p^\kappa(t,x,y) \le C t^{-d/\alpha} \P_x(\zeta^\kappa>t/3) \P_y(\zeta^\kappa>t/3), \quad t>0, \, x,y \in \R^d_+.	
\end{equation*}
\end{lemma}
\begin{proof} By the semigroup property, the symmetry of $p^\kappa(t, \cdot, \cdot)$ and Proposition \ref{p:upper-heatkernel}, we obtain
\begin{align*}		
&p^\kappa(t,x,y) = \int_{\R^d_+}\int_{\R^d_+} p^\kappa(t/3,x,z)p^\kappa(t/3,z,w) p^\kappa(t/3,y,w) dzdw \\
&\le c_1t^{-d/\alpha} \int_{\R^d_+} p^\kappa(t/3,x,z)dz \int_{\R^d_+} p^\kappa(t/3,y,w)dw  = c_1t^{-d/\alpha} \P_x(\zeta^\kappa>t/3) \P_y(\zeta^\kappa>t/3).
\end{align*}
\end{proof}

The next lemma shows that  \eqref{s:UHK-rough} 
(hence Proposition \ref{p:UHK-rough}) is a consequence of the following, seemingly weaker, inequality: There exists $C>0$ such that
\begin{equation}\label{e:UHKDn}
		p^\kappa(t,x,y) \le C \left(1 \wedge \frac{x_d}{t^{1/\alpha}}\right)^{q}\left(1 \wedge \frac{y_d}{t^{1/\alpha}}\right)^{q} t^{-d/\alpha}, \quad  t>0, \, x,y \in \R^d_+.
	\end{equation}
\begin{lemma}\label{l:(5.6)implies(5.1)}
If \eqref{e:UHKDn} holds true,  then  \eqref{s:UHK-rough} also holds.
\end{lemma}
\begin{proof} We claim that there exists a constant $c_1>0$ such that 
\begin{equation}\label{e:UHK-claim}
	p^\kappa(t,x,y) \le c_1 \left(1 \wedge \frac{x_d}{t^{1/\alpha}}\right)^{q} 
	p_\alpha(t, x, y).
\end{equation} 
By \eqref{e:kernel-scaling}, 
we can assume $\wt x=\wt 0$ and $t=t_0=(1/2)^\alpha$. If $x_d \ge 2^{-4}t_0^{1/\alpha}$ or $|x-y| \le 4t_0^{1/\alpha}$, then \eqref{e:UHK-claim} follows from 
Proposition \ref{p:upper-heatkernel} or  the assumption \eqref{e:UHKDn}  respectively.
Hence, we assume $x_d<2^{-4}t_0^{1/\alpha}$ and $|x-y|>4t_0^{1/\alpha}$, and  will show that
\begin{equation}\label{e:UHK-claim1}
	p^\kappa(t_0,x,y) \le c_1(t_0) \frac{x_d^q}{|x-y|^{d+\alpha}}.
\end{equation}

Let $V_1=U(t_0^{1/\alpha})$, $V_3=\{w \in \R^d_+:|w-y| < |x-y|/2\}$ and $V_2=\R^d_+ \setminus (V_1 \cup V_3)$. By 
Lemma \ref{l:l.5.2-gen}(iii)
we have 
\begin{align}\label{l:lemma-upper-1-1}
	 \P_x(\tau^\kappa_{V_1}<t_0< \zeta^\kappa)  \le \P_x(Y^\kappa_{\tau^\kappa_{V_1}} \in \R^d_+)  \le c_2(t_0^{-1/\alpha}x_d)^q.
\end{align}
Also, we get from Proposition \ref{p:upper-heatkernel} that 
\begin{equation}\label{e:UHK-rough-0}
	\sup_{s \le t_0, \, z \in V_2} p^\kappa(s,z,y) \le c_3 \sup_{s \le t_0, \, z \in \R^d_+, |z-y|>|x-y|/2} \frac{s}{|z-y|^{d+\alpha}} = 2^{d+\alpha}c_3 \frac{t_0}{|x-y|^{d+\alpha}}.
\end{equation}
Next, we note that by the triangle inequality, for any $u \in V_1$ and $w \in V_3$,
\begin{equation}\label{e:UHK-rough-1}
	|u-w| \ge |x-y| - |x-u| - |y-w| \ge |x-y| - t_0^{1/\alpha} - \frac{|x-y|}{2} \ge \frac{|x-y|}{4} \ge t_0^{1/\alpha} \ge u_d.
\end{equation}
In particular, recalling  that $\beta_1>0$ if $\beta_3>0$, 
we see that  $\big(1+{\bf 1}_{|w|\ge1}(\log|w|)^{\beta_3}\big)|u-w|^{-\beta_1} \le c_4$ for $u \in V_1$ and $w \in V_3$, so by Lemma \ref{l:l.5.2-gen}(i), we have that for any $u \in V_1$ and $w \in V_3$,
\begin{equation}\label{e:UHK-rough-11}
	\sB(u,w) \le 
	c_5 u_d^{\beta_1}(|\log u_d|^{\beta_3}\vee 1)\big(1+{\bf 1}_{|w|\ge1}(\log|w|)^{\beta_3}\big)|u-w|^{-\beta_1}
	\le c_6 u_d^{\beta_1} |\log u_d|^{\beta_3}.
\end{equation}
Thus, by {\bf (A3)}(II), 
Lemma \ref{l:l.5.2-gen}(ii) and
\eqref{e:UHK-rough-1} and \eqref{e:UHK-rough-11} we get that
\begin{align} \label{e:UHK-rough111}
	&\int_0^{ t_0 }\int_{V_3} \int_{V_1} p^{\kappa, V_1}(t_0 -s, x, u) \sB(u,w) p^\kappa(s, y,w) du dw ds \\
	&\le c_7
	\int_0^{t_0}\left(\int_{V_1} p^{\kappa, V_1}(t_0-s, x, u) u_d^{\beta_1} |\log u_d|^{\beta_3} du \right) \left( \int_{V_3}p^\kappa(s, y,w) dw\right) ds \nn\\
	&\le c_7
	\int_0^\infty\left(\int_{V_1} p^{\kappa, V_1}(s, x, u) u_d^{\beta_1} |\log u_d|^{\beta_3}  du \right)  ds \nn\\
	&=
	c_7
	\E_x
	\int_0^{\tau^\kappa_{V_1} } (Y_s^{\kappa,d})^{\beta_1} |\log Y_s^{\kappa,d}|^{\beta_3}   ds\le c_8 x_d^q. \nn
\end{align}
Now  \eqref{e:UHK-claim1} (and so \eqref{e:UHK-claim}) follows
from \eqref{l:lemma-upper-1-1}--\eqref{e:UHK-rough-1}, \eqref{e:UHK-rough111}  and  Lemma \ref{l:general-upper-2}. 

Finally,  by the semigroup property, symmetry, \eqref{e:stableu2} and \eqref{e:UHK-claim}, 
\begin{align*}
	&p^\kappa(t,x,y) = \int_{\R^d_+} p^\kappa(t/2,x,z)p^\kappa(t/2,y,z)dz \\
	&\le c_1^2 \left(1 \wedge \frac{x_d}{t^{1/\alpha}}\right)^{q}  \left(1 \wedge \frac{y_d}{t^{1/\alpha}}\right)^{q} \int_{\R^d_+}
	 p_\alpha(t/2, x, z)p_\alpha(t/2, z, y)
	 dz\\
	&\le c_9 \left(1 \wedge \frac{x_d}{t^{1/\alpha}}\right)^{q}  \left(1 \wedge \frac{y_d}{t^{1/\alpha}}\right)^{q} 
	p_\alpha(t, x, y)
\end{align*}
\end{proof}
Now we prove \eqref{e:UHKDn} holds.  We first consider the case $q<\alpha$.
\begin{lemma}\label{l:UHKD_n}
	If $q<\alpha$, then  \eqref{e:UHKDn} holds true. 
\end{lemma}
\begin{proof} By \eqref{e:kernel-scaling}, it suffices to prove \eqref{e:UHKDn} when $t=1$.
By Lemma  \ref{l:general-upper-1}, it suffices to prove  that there exists a constant $c_1>0$ such that  $\P_x(\zeta^\kappa>1/3) \le c_1(1 \wedge x_d)^q$ for all $x \in \R^d_+$.
By Lemma \ref{l:translation-invariance} and the fact that $\P_x(\zeta>1/3) \le 1$, 
without loss of generality, we can assume $\wt x=0$ and $x_d<2^{-5}$.
Then, since $q<\alpha$, by \eqref{e:lemma-upper-3-2} and 
Lemma \ref{l:l.5.2-gen}(iii)--(iv) (with $r=1$),
we get  $\P_x(\zeta^\kappa>1/3)\le c_2 x_d^q$.
\end{proof}

To remove the assumption $q<\alpha$ in Lemma \ref{l:UHKD_n},  we will make use of 
a result from \cite{KSV-Green, KSV-nokill}. Recall from  Remark \ref{r:dalpha} that all results in \cite{KSV-Green, KSV-nokill} are valid when $d>\alpha$. 
Now we  state, and will prove later,  that  the desired result  still holds true for $d=1 \le \alpha$ (thus for all $d$).

\begin{lemma}\label{l:6.10}
	Let $\gamma>q-\alpha$. There exists  $C>0$ such that	for any $R>0$, $U(R) \subset D \subset U(2R)$ and any $x=(\wt 0, x_d) \in \R^d_+$ with $x_d \le R/10$, it holds that 
	\begin{equation*}
\E_x \int_0^{\tau^\kappa_{D}} (Y^{\kappa, D,d}_{s})^{\gamma} ds = \int_0^{\infty} \int_D p^{\kappa, D}(t,x,z)z_d^\gamma dz dt \le C R^{\gamma+\alpha-q} x_d^q.
	\end{equation*}
\end{lemma}
\begin{proof} 
When $d>\alpha$, by Remark \ref{r:dalpha}, the result follows from \cite[Proposition 6.10]{KSV-Green} if $\kappa>0$, and from   \cite[Proposition 6.8]{KSV-nokill} if $\kappa=0$.  The proof of the case $d=1\le \alpha$ is postponed to  
the end of this section. 
\end{proof}

\begin{lemma}\label{l:UHKD}
	Inequality \eqref{e:UHKDn} holds true. 
\end{lemma}
\begin{proof} Again using  \eqref{e:kernel-scaling}, it suffices to prove \eqref{e:UHKDn} when $t=1$.

 Inequality \eqref{e:UHKDn}  holds when  $q<\alpha$ by Lemma \ref{l:UHKD_n}.
Now,
assume that \eqref{e:UHKDn} holds  for all $q<k\alpha$
for some $k\in \mathbb{N}$. We now show \eqref{e:UHKDn} also holds for $q\in [k\alpha,
(k+1)\alpha)$ and hence \eqref{e:UHKDn} always holds by induction.

Fix $\eps \in (0,\alpha)$ such that $q-\alpha + \eps<k\alpha$.   Then $q-\alpha + \eps<q$ and 
$\wt \kappa:=C(\alpha, q-\alpha+\eps , \sB)<C(\alpha, q, \sB)=\kappa$.
By \eqref{e:kernel-scaling} and the induction hypothesis, it holds that for any $s,u\in[0,1/4]$ and $z,w \in \R^d_+$, 
\begin{align*}	
&p^\kappa(1-s-u,z,w) =(1-s-u)^{-d/\alpha} p^\kappa(1, (1-s-u)^{-1/\alpha}z,(1-s-u)^{-1/\alpha}w)\\	
&\le 2^{d/\alpha}p^{\wt \kappa}(1, (1-s-u)^{-1/\alpha}z,(1-s-u)^{-1/\alpha}w)  \le c_3 (1 \wedge z_d)^{q-\alpha+\eps}  (1 \wedge w_d)^{q-\alpha+\eps}.
\end{align*}
Here in the first inequality above we used the 
fact that 
$p^\kappa(1,x,y)\le p^{\wt \kappa}(1,x,y)$,
which is a consequence of 
$\kappa> \wt \kappa$.
Therefore, by the semigroup property and symmetry, we get
\begin{align*}
	&p^\kappa(1,x,y) = 16\int_0^{1/4} \int_0^{1/4} p^\kappa(1,x,y) dsdu\\
	&= 16\int_0^{1/4} \int_0^{1/4} \int_{\R^d_+}\int_{\R^d_+} p^\kappa(s,x,z) p^\kappa(1-s-u,z,w) p^\kappa(u,y,w)dzdw dsdu\\
	& \le 16c_3  \int_0^{1/4} \int_{\R^d_+} p^\kappa(s,x,z) (1 \wedge z_d)^{q-\alpha+\eps} dz ds \int_0^{1/4} \int_{\R^d_+} p^\kappa(u,y,w) (1 \wedge w_d)^{q-\alpha+\eps} dw du.
\end{align*}
Thus, to conclude \eqref{e:UHKDn} by induction, it suffices to show that there exists 
$c_4>0$ such that 
\begin{equation}\label{e:UHKD-claim}
	\int_0^{1/4} \int_{\R^d_+} p^\kappa(s,v,z) (1 \wedge z_d)^{q-\alpha+\eps} dz ds \le c_4 (1 \wedge v_d)^q, \quad v \in \R^d_+.
\end{equation}
By Lemma \ref{l:translation-invariance}, we can assume $\wt v=0$.  If $v \notin U(2^{-4})$, then we get
\begin{align*}
(1 \wedge v_d)^q &\ge 2^{-5q} \ge 2^{-5q} \int_0^{1/4} \int_{\R^d_+} p^\kappa(s,v,z) (1 \wedge z_d)^{q-\alpha+\eps} dz ds.
\end{align*}
Otherwise, if $v \in U(2^{-4})$, then by \eqref{e:lemma-upper-3-1},  Fubini's theorem, 
and Lemmas  \ref{l:l.5.2-gen}(iii)  and \ref{l:6.10},
\begin{align*}
	&\int_0^{1/4} \int_{\R^d_+} p^\kappa(s,v,z) (1 \wedge z_d)^{q-\alpha+\eps} dz ds\\
	& \le \int_0^{1/4} \E_v \left[ (1 \wedge Y^{\kappa, U(1),d}_{s})^{q-\alpha+\eps} : \tau^\kappa_{U(1)}>s\right]  ds+ \int_0^{1/4} \P_v(Y^\kappa_{\tau^\kappa_{U(1)}} \in \R^d_+) ds\\
	& \le \E_v \int_0^{\tau^\kappa_{U(1)}} (1 \wedge Y^{\kappa, U(1),d}_{s})^{q-\alpha+\eps} ds   +  \frac14\P_v(Y^\kappa_{\tau_{U(1)}} \in \R^d_+) \le  c_5v_d^q.
\end{align*}
This completes the proof. 
\end{proof}

\noindent \textsc{Proof of Proposition \ref{p:UHK-rough}:}  
The assertion is a direct consequence of  Lemmas \ref{l:UHKD} and \ref{l:(5.6)implies(5.1)}.
\qed

\begin{remark} \label{r:dalpha2}	
{\rm 
Note that, when $d=1$,  \cite[Lemma 7.1]{KSV-Green} holds for $p\in ((\alpha-1)_+, \alpha+\beta_1)$ since 
the parameter $\beta_2$ is irrelevant. See Remark \ref{r:d=1}. Moreover, in the proof of \cite[Theorem 1.2]{KSV-Green}, 
only the upper bound in \cite[Proposition 6.10]{KSV-Green} is used. 
Thus, by  Remark \ref{r:dalpha}	 and using our Lemma \ref{l:6.10} 
instead of \cite[Proposition 6.10]{KSV-Green}, we see that 
\cite[Theorem 1.2]{KSV-Green}
 is valid in case and $d\ge 2$ and $p\in ((\alpha-1)_+,\alpha+(\beta_1\wedge \beta_2))$, and in case $d=1$ for all $p\in ((\alpha-1)_+, \alpha+\beta_1)$. 
}
\end{remark}

In the remainder of this section we complete the proof of Lemma \ref{l:6.10} in case $d=1\le \alpha$. 
\begin{lemma}\label{l:G1}
	If $d=1 \le \alpha$, then there exists $C>0$ such that for all $x,y \in (0, \infty)$,
\begin{equation}\label{e:G1-new}
		G^\kappa(x,y) \le C   \Big(1 \wedge \frac{x \wedge y}{|x-y|}\Big)^{q \wedge (\alpha- \frac12)} (x \vee y)^{\alpha-1} \L\Big(\frac{x \vee y}{|x-y|}\Big).
\end{equation}
\end{lemma}
\begin{proof}  Let  $\wt q:=q \wedge (\alpha-\frac12) \in (0,q] \cap  [\alpha-1, \alpha)$. By symmetry and scaling
\eqref{e:kernel-scaling},  we can assume $|x-y|=1$ and $x \le y$ without loss of generality. Then  $y = x +1 > x \vee 1$. It suffices to show that
\begin{align*}
	G^\kappa(x,y) \le C x^{\wt q} y^{\alpha-1}\L(y).
\end{align*}
 Note that $1+2\wt q >1+ \wt q  \ge \alpha$. Since $C(\alpha,\wt q, \sB) \le C(\alpha, q, \sB)=\kappa$ implies that 
$p^\kappa(t,x,y) \le p^{C(\alpha, \wt q, \sB)}(t,x,y) $, 
using the fact $y>x \vee 1$ and 
Lemmas \ref{l:UHKD_n} and  \ref{l:(5.6)implies(5.1)}, 
 we get
\begin{align*}
&G^\kappa(x,y)\le \int_0^\infty p^{C(\alpha, \wt q, \sB)}(t,x,y) dt \\
&\le 	c_1\Big(  x^{\wt q}\int_0^{1}  t^{(\alpha - \wt q)/\alpha}dt
	+ x^{\wt q}\int_1^{y^\alpha} t^{-(1+\wt q)/\alpha}dt 
	+x^{\wt q} y^{\wt q}\int_{y^\alpha}^{\infty} t^{-(1+2\wt q)/\alpha}
	dt\Big)\\
	&\le c_1 x^{\wt q} \Big(1 + \log( y^\alpha) +\frac{\alpha}{1+2\wt q -\alpha} y^{\alpha-\wt q-1} \Big)\le c_2 x^{\wt q} y^{\alpha-1} \L (y).
\end{align*}
The proof is complete.
\end{proof}

We now improve \eqref{e:G1-new} by removing the term $\alpha- \frac12$ from the power of the first factor. 

\begin{lemma}\label{l:G2}
	If $d=1 \le \alpha$, then there exists $C>0$ such that for all $x,y \in (0, \infty)$,
	\begin{align*}
		G^\kappa(x,y) \le C  \Big(1 \wedge \frac{x \wedge y}{|x-y|}\Big)^{q} (x \vee y)^{\alpha-1} 
		\L\Big(\frac{x \vee y}{|x-y|}\Big).
	\end{align*}
\end{lemma}
\begin{proof} 
Let $r=2^{-6}$. 
By symmetry, \eqref{e:kernel-scaling} and Lemma \ref{l:G1},  without loss of generality, we can assume $|x-y|=1$ and $x  <y \wedge (r/2)$. Note that $y=x+1 \in (1, 1+r/2)$. It suffices to show that $
	G^\kappa(x,y) \le C x^{q}.$

For $z\in (0, r/2)$ and $w\in (r, \infty)$, we have $z<w/2 \le w-z <w$ and $w \vee y \asymp w$. Thus, for any $z \in (0,r/2)$, using 
Lemmas \ref{l:l.5.2-gen}(i)
and  \ref{l:G1}, we obtain
\begin{align*}
	&\int_{r}^{\infty}G^\kappa(w,y) \sB(z,w)(w-z)^{-1-\alpha} dw\nonumber\\
	&\le c_1z^{\beta_1}(|\log z|^{\beta_3}\vee 1)\int_{r}^{\infty}\frac{G^\kappa (w,y)}{w^{1+\alpha+\beta_1}}
	\big(1+{\bf 1}_{\{w\ge1\}}(\log w)^{\beta_3}\big)
	dw\nn\\
	&\le c_2 z^{\beta_1} |\log z|^{\beta_3}\int_{r}^{\infty} \L \Big(\frac{w }{|w-y|}\Big)
	\frac{ \big(1+{\bf 1}_{\{w\ge1\}}(\log w )^{\beta_3}\big)}{w^{2+\beta_1}} 
	dw.
\end{align*}
Hence, by the L\'evy system formula and 
Lemma \ref{l:l.5.2-gen}(ii),
we get 
\begin{align}\label{e:green-d=1}
	&\E_x\left[G^\kappa(Y^\kappa_{\tau^\kappa_{(0, r/2)}}, y); Y^\kappa_{\tau^\kappa_{(0, r/2)}}\notin (0, r)\right]\\
	&\le c_2 \E_x\int^{\tau^\kappa_{(0, r/2)}}_0(Y^{\kappa,d}_t)^{\beta_1}  |\log (Y^{\kappa, d}_t)|^{\beta_3}dt\int_{r}^{\infty} 
	\L\Big(\frac{w}{|w-y|}\Big)
	\frac{ \big(1+{\bf 1}_{\{w\ge1\}}(\log w )^{\beta_3}\big)}{w^{2+\beta_1}}
	dw\nn\\
	&\le c_3 x^q \int_{r}^{\infty} \L \Big(\frac{w }{|w-y|}\Big)
	\frac{ \big(1+{\bf 1}_{\{w\ge1\}}(\log w )^{\beta_3}\big)}{w^{2+\beta_1}} dw. \nn
\end{align}
Since $w-y \ge w/2$ for $w \in (2y,\infty)$ and $y \in (1,1+r/2)$, using  a change of the variables,  we obtain
\begin{align*}
	&\int_{r}^{\infty} \L \Big(\frac{w }{|w-y|}\Big)
	\frac{ \big(1+{\bf 1}_{\{w\ge1\}}(\log w )^{\beta_3}\big)}{w^{2+\beta_1}}
	dw\\
	& \le 	\frac{ 1+(\log (2+r) )^{\beta_3}}{r^{2+\beta_1}}\int_{r}^{2y} \L \Big(\frac{2+r }{|w-y|}\Big)
	dw +\L(2) \int_{2y}^{\infty} 
	\frac{ 1+(\log w )^{\beta_3}}{w^{2+\beta_1}}
	dw \\
		& \le  \frac{2( 1+(\log (2+r) )^{\beta_3})}{r^{2+\beta_1}}\int_{0}^{1+r/2} \L \Big( \frac{2+r }{v}\Big)
	dv +c_4\le c_5.
\end{align*}
Thus, by \eqref{e:green-d=1}, it holds that
\begin{align}\label{e:TAMSe6.41}
	\E_x\left[G^\kappa(Y^\kappa_{\tau^\kappa_{(0, r/2)}}, y); Y^\kappa_{\tau^\kappa_{(0, r/2)}}\notin(0, r)\right] \le 
	c_3c_5 x^q.
\end{align}

Since  $z \mapsto G^\kappa(z, y)$ is harmonic in $(0,2r)$ with respect to $Y^\kappa$ and vanishes continuously  as $z \to 0$ by Lemma \ref{l:G1}, we get from Proposition \ref{p:Carleson} that $G^\kappa(z,y) \le c_6 G^\kappa(r,y)$ for all $z \in (0,r)$. Therefore,  
 using 
 Lemmas \ref{l:l.5.2-gen}(iii)
 and \ref{l:G1}, since $y \in (1,1+r/2)$, we obtain
\begin{align}\label{e:TAMSe6.42}
	&\E_x\left[G^\kappa(Y^\kappa_{\tau^\kappa_{(0, r/2)}}, y); Y^\kappa_{\tau^\kappa_{(0, r/2)}}\in (0, r)\right]\le c_6 G^\kappa(r, y)\P_x(Y^\kappa_{\tau^\kappa_{(0, r/2)}}\in (0, r))\\
	&\le c_7 \Big(\frac{x}{r}\Big)^q\Big( \frac{r}{y-r} \Big)^{q \wedge (\alpha - \frac12)} y^{\alpha-1} 
	\L\Big(\frac{y}{y-r}\Big)  \le c_8x^q. \nn
\end{align}
Combining \eqref{e:TAMSe6.41} and \eqref{e:TAMSe6.42} and using the   harmonicity of $z \mapsto G^\kappa(z, y)$ in $(0,2r)$,
we arrive at
\begin{align*}
	G^\kappa(x, y)&=\E_x\left[G^\kappa(Y^\kappa_{\tau^\kappa_{(0, r/2)}}, y); Y^\kappa_{\tau^\kappa_{(0, r/2)}}\notin (0, r)\right]+ \E_x\left[G^\kappa(Y^\kappa_{\tau^\kappa_{(0, r/2)}}, y); Y^\kappa_{\tau^\kappa_{(0, r/2)}}\in (0,r)\right]\\
	&\le c_9x^q.
\end{align*}
The proof is complete.
\end{proof}

\noindent \textsc{Proof of Lemma \ref{l:6.10} for $d=1\le \alpha$:} 
Assume that $d=1 \le \alpha$. Since  
$D \subset U(2R)=[0,R)$,
using Fubini's theorem and Lemma \ref{l:G2},  we have 
\begin{align*} 
	&\int_0^{\infty} \int_D p^{\kappa,D}(t,x,z)z^\gamma dz dt \le  \int_0^R \int_0^{\infty}  p^{\kappa,D}(t,x,z)
	  dt\,  z^\gamma dz 
	  \le   \int_0^R G^\kappa(x,z)z^\gamma dz\nn\\ 
	&\le c_1 \int_0^{x/2}   \frac{z^{q+\gamma}}{|x-z|^q} x^{\alpha-1} \L\Big(\frac{x}{|x-z|}\Big)dz + c_1 
	 x^{\gamma+\alpha-1} \int_{x/2}^{2x} 
	 \L\Big(\frac{2x}{|x-z|}\Big)dz \\
	&\quad + c_1 \int_{2x}^R  \frac{ x^q}{|x-z|^q} z^{\gamma+\alpha-1}\L\Big(\frac{z}{|x-z|}\Big)dz=:c_1(I+II+III).
\end{align*}
Since $\gamma+\alpha-q>0$ and $x\le R/10$, we see that
\begin{align*}
	 I\le x^{q+\gamma+\alpha-1}\int_0^{x/2}   \frac{1}{(x/2)^q}  \L\Big(\frac{x}{(x/2)}\Big)dz  
	 \le 2^{q-1} \L(2) x^{\gamma+\alpha} \le c_2 R^{\gamma+\alpha-q} x^q.
\end{align*}
Next, using the change of the variables $z=xy$, we also get
\begin{align*}
	II = x^{\gamma+\alpha} \int_{1/2}^{2}\L\Big(\frac{2}{|1-y|}\Big)dy\le c_3 x^{\gamma+\alpha} \le c_3R^{\gamma+\alpha-q} x^q.
\end{align*}
Lastly, since $|x-z| \ge z/2$ for all $z \ge 2x$, we obtain
\begin{align*}
III \le 2^q\L(2) x^q\int_{2x}^R  z^{\gamma+\alpha-q-1}dz \le c_4 R^{\gamma+\alpha-q} x^q.
\end{align*}
The proof is complete.
\qed

\section{Sharp heat kernel lower bounds}\label{s:lowerb}

\textit{In this and the next section, 
whenever we consider $\rY$, we assume  that {\bf (A1)}, {\bf (A3)} and {\bf (A4)} hold, and whenever we consider $Y^\kappa$, we assume that all {\bf (A1)}--{\bf (A4)} and \eqref{e:killing-potential} hold with $q_\kappa\in [(\alpha-1)_+, \beta_1)$.} Note that {\bf (IUBS)} holds  under this setting by Lemma \ref{l:check-UJS}. 

We fix $\kappa \in [0,\infty)$. The following notational convention
will be used throughout this and the next section.
When we consider $Y^{\kappa}$, we assume that $\kappa>0$ if $\alpha \in (0,1]$ (see Lemma \ref{l:alpha<1}), and denote by $q$ the strictly positive constant $q_k \in ((\alpha-1)_+, \alpha+\beta_1)$ from \eqref{e:killing-potential}. Additionally, we write $Y$,  $p(t,x,y)$, $p^D(t,x,y)$, $\tau_D$ and $\zeta$ instead of $Y^{\kappa}$, $p^{\kappa}(t,x,y)$, $p^{\kappa, D}(t,x,y)$, $\tau_D^{\kappa}$ and $\zeta^\kappa$. When we consider 
$\overline{Y}$, the letter $q$ denotes 0, and we write $Y$, $p(t,x,y)$, $p^D(t,x,y)$, $\tau_D$ and $\zeta$ instead of $\overline{Y}$, 
$\bar{p}(t,x.y)$, $\bar{p}^D(t,x,y)$, $\overline{\tau}_D$ and $\infty$.

\medskip
Recall the definitions of $V_x(t)$ and $W_x(t)$ from \eqref{e:def-VW}. We let $V_x:=V_x(1)$ and $W_x:=W_x(1)$.

\begin{lemma}\label{l:M}
 	There exist constants $M>1$ and $C>0$ such that 
	for all $t>0$ and  $x \in \R^d_+$, 
$$ 
\inf_{z \in W_x(t)} p(Mt,x,z) \ge C  \left(1 \wedge \frac{x_d}{t^{1/\alpha}}\right)^{q}t^{-d/\alpha}.
$$
\end{lemma}
\begin{proof}  When $q=0$, the result 
is given in  Lemma \ref{l:bar-lower}. 

Suppose $q>0$. 
By the scaling property \eqref{e:kernel-scaling}, 
it suffices to prove the lemma for $t=1$. If $x_d  \ge 1$, then the result follows from Proposition \ref{p:lower-heatkernel}. 
Assume $x_d <1$. By Proposition \ref{p:lower-heatkernel},  for any $M>1$, there is $c_1=c_1(M)>0$  such that 
\begin{align}
 \inf_{w, z \in W_x, \, 1 \le s \le M} p (s, w,z) \ge  c_1  \quad \text{for all $x \in \R^d_+$}.
	\label{e:lbd1-}
\end{align} 
Using the strong Markov property and \eqref{e:lbd1-},
we see that for all $M >1$ and $z \in W_x$,
\begin{align*}
	p(M,x,z) &\ge 	\E_x\left[\,p(M- \tau_{V_x}, Y_{\tau_{V_x}}, z): \tau_{V_x} \le M-1,\, Y_{\tau_{V_x}} \in W_x\right]\nn\\
	&\ge \Big(
 \inf_{w \in W_x, \, 1\le s \le M } p (s, w,z) 	\Big)	\P_x \left(\tau_{V_x} \le M-1, \,Y_{\tau_{V_x}} \in W_x\right) \nn\\
 	&\ge c_1 \left(\P_x ( Y_{\tau_{V_x}} \in W_x) -\P_x (\tau_{V_x(t)} > M-1) \right).
\end{align*}
Note that, by \cite[Lemma 5.10]{KSV-jump} for $\kappa>0$ and \cite[Theorem 1.1]{KSV-nokill} for $\kappa=0$, we have 
\begin{equation*}
	\P_x ( Y_{\tau_{V_x}} \in W_x) \ge 2 c_2x_d^q
\end{equation*}
and,   by Corollary \ref{c:life}, we also have 
\begin{equation*}
	\P_x (\tau_{V_x} > M-1) \le \P_x (\zeta > M-1) \le c_3  (x_d/(M-1)^{1/\alpha})^q.
\end{equation*}
Thus, we can choose $M =1+ (c_3/c_2)^{\alpha/q}$  so that 
$2c_2-c_3  (M-1)^{-q/\alpha} = c_2$, which implies  
$p(M,x,z)  \ge c_2 x_d^q$.  \end{proof}

Recall that $\L^{b}(a)=\log^{b}(e+a)$. 
For any $b_1, b_2, b_3, b_4\ge 0$, we define for $t\ge0$ and $x,y \in \oR^d_+$,
\begin{align}\label{def:wtB2}	
A_{b_1, b_2, b_3, b_4}(t,x,y)&
:= 	
\Big(\frac{(x_d \wedge y_d) \vee t^{1/\alpha}}{|x-y|} \wedge 1 \Big)^{b_1}  \Big(\frac{(x_d \vee y_d) \vee t^{1/\alpha}}{|x-y|} \wedge 1 \Big)^{b_2}  \\	
& \,\times  \L^{b_3}\Big(\frac{ ((x_d \vee y_d) \vee t^{1/\alpha}) \wedge |x-y|}{((x_d \wedge y_d) \vee t^{1/\alpha}) \wedge |x-y|} \Big) \,\L^{b_4}\Big(\frac{ |x-y|}{( (x_d \vee y_d) \vee t^{1/\alpha}) \wedge |x-y|} \Big). \nn
\end{align}

\begin{remark}\label{r:Ad=1}
	{\rm When $d=1$,  the parameters $b_2$ and $b_4$ in the above definition \eqref{def:wtB2}	  are irrelevant (see  Remark \ref{r:d=1}).}
\end{remark}

We first note that 
$A_{b_1, b_2, b_3, b_4}(0, x,y)=B_{b_1, b_2, b_3, b_4}(x,y)$
for $x,y \in \R^d_+$ and 
\begin{align}\label{e:comp-AB}
	A_{b_1, b_2, b_3, b_4}(t, x,y)  \asymp 
B_{b_1, b_2, b_3, b_4}(x + t^{1/\alpha}\e_d,y + t^{1/\alpha}\e_d) 
\quad \text{ for }  t \ge 0, x,y \in \R^d_+,
\end{align}
since $a \vee b \asymp a+b$ for all $a,b \ge 0$. 
Note also that for any $a>0$, 
there exists $c>0$ such that
\begin{equation}\label{e:wtB-interior}
	A_{b_1, b_2, b_3, b_4}(t,x,y) \ge 		c(a\wedge 1)^{b_1 +b_2},
\end{equation}
for all $t>0$ and $x,y \in \R^d_+$ with $(x_d \wedge y_d)+t^{1/\alpha} \ge a |x-y|$.

\begin{prop}\label{p:HKE-lower-1}
 	There exists a constant $C>0$ such that for all $t>0$ and $x,y \in \R^d_+$,
	\begin{align*}
		p(t,x,y) \ge C \left(1 \wedge \frac{x_d}{t^{1/\alpha}}\right)^{q} 
		\left(1 \wedge \frac{ y_d}{t^{1/\alpha}}\right)^{q} A_{\beta_1, \beta_2, \beta_3, \beta_4}(t,x,y) 
		p_\alpha(t, x, y).
	\end{align*}
\end{prop}
\begin{proof}
Without loss of generality, we assume $x_d \le y_d$.
Let $M> 1$ be the constant in Lemma \ref{l:M}.
By the semigroup property and Lemma \ref{l:M}, 
\begin{align}\label{e:lb11}
	p(t,x,y)&\ge
	\int_{ W_x(t/(3M))}\int_{W_y(t/(3M))}
	p(t/3,x,z) \,p(t/3,z,w) \,p(t/3,w,y)
	dzdw\\
	& \ge c_1 t^{2d/\alpha}  \left(\inf_{z \in W_x(t/(3M)) } 
	p(t/3,x,z) \right)  \left(\inf_{w \in W_y(t/(3M)) } p(t/3,y,w) \right)\nn\\
	&\quad \times 
	\left(\inf_{(z,w) \in W_x(t/(3M))\times W_y(t/(3M))} p(t/3,z,w) \right) \nn\\
	& \ge c_2 \left(1 \wedge \frac{x_d}{t^{1/\alpha}}\right)^{q} \left(1 \wedge \frac{y_d}{t^{1/\alpha}}\right)^{q}	
	\inf_{(z,w) \in W_x(t/(3M))\times W_y(t/(3M))} p(t/3,z,w). \nn
\end{align}
For all $(z,w) \in W_x(t/(3M))\times W_y(t/(3M))$, we have $z_d \ge x_d + 5(t/(3M))^{1/\alpha}$, $w_d \ge y_d + 5(t/(3M))^{1/\alpha}$,
$|z-w| \le |z-x| + |x-y| + |y-w| \le |x-y| + 20(t/(3M))^{1/\alpha}$ and $|z-w| \ge  |x-y|-|z-x|- |y-w| \ge |x-y| - 20(t/(3M))^{1/\alpha}$.
 Hence, if $|x-y|\le 21(t/(3M))^{1/\alpha}$, then 
 for all $(z,w) \in  W_x(t/(3M))\times W_y(t/(3M))$,
  \begin{align}\label{e:LHK-diag}
 z_d \wedge w_d \ge 5(\frac{t}{3M})^{1/\alpha} \ge \frac{5}{41} \big((\frac{t}{3M})^{1/\alpha} \vee |z-w|\big).
  \end{align}
If  $|x-y|>21(t/(3M))^{1/\alpha}$, then using  {\bf (A3)}(I)  in the first inequality below, 
 $|z-w|\le 2|x-y|$ and Lemma \ref{cal:basic} in the second  and \eqref{e:comp-AB}  in the third,  we get 
\begin{align}\label{e:LHK-jump}
& \inf_{(z,w) \in W_x(t/(3M))\times W_y(t/(3M))}J(z,w) \ge C_2^{-1} 	
\inf_{(z,w) \in W_x(t/(3M))\times W_y(t/(3M))} \frac{B_{\beta_1,\beta_2,\beta_3,\beta_4} (z,w)}{|z-w|^{d+\alpha}}\\
&\ge c_3\frac{B_{\beta_1,\beta_2,\beta_3,\beta_4} 
(x + 5 (t/(3M))^{1/\alpha} \e_d,y+ 5 (t/(3M))^{1/\alpha} \e_d) }{|x-y|^{d+\alpha}}
\ge c_4\frac{A_{\beta_1,\beta_2,\beta_3,\beta_4} (t,x,y) }{|x-y|^{d+\alpha}}. \nn
\end{align}

Now, by Propositions \ref{p:lower-heatkernel} and \ref{p:pld}, \eqref{e:LHK-diag}, \eqref{e:LHK-jump} and \eqref{e:wtB-interior}, we have 
\begin{align}\label{e:lb12}
&\inf_{(z,w) \in W_x(t/(3M))\times W_y(t/(3M))} p(t/3,z,w)\\
& \ge c_5 \begin{cases}		t^{-d/\alpha} 	
& \text{ if } |x-y| \le 21(t/(3M))^{1/\alpha},\\[4pt] 	
tA_{\beta_1,\beta_2,\beta_3,\beta_4} (t,x,y) |x-y|^{-d-\alpha } 	
& \text{ if } |x-y| > 21(t/(3M))^{1/\alpha}\end{cases}\nn\\
& \ge c_6 A_{\beta_1, \beta_2, \beta_3, \beta_4}(t,x,y) 
p_\alpha(t, x, y).
\nn
\end{align}
Combining \eqref{e:lb11} and \eqref{e:lb12}, 
we get the desired result.
\end{proof}

\begin{cor}\label{c:life-2}
	It holds that for any $t>0$ and $x \in \R^d_+$,
	\begin{align*}
		\P_x (\zeta>t) \asymp   \left(1 \wedge \frac{x_d}{t^{1/\alpha}}\right)^{q}.
	\end{align*}
\end{cor}
\begin{proof} Using Proposition \ref{p:HKE-lower-1} and \eqref{e:wtB-interior}, we get
\begin{align*}
		\P_x (\zeta>t) &\ge \int_{B(x,2t^{1/\alpha}): y_d \ge t^{1/\alpha}} p(t,x,y)dy \ge c_1  \left(1 \wedge \frac{x_d}{t^{1/\alpha}}\right)^{q}t^{-d/\alpha} \int_{B(x,2t^{1/\alpha}): y_d \ge t^{1/\alpha}} dy\\
		& \ge c_2\left(1 \wedge \frac{x_d}{t^{1/\alpha}}\right)^{q}.
\end{align*}
Combining the above with Corollary \ref{c:life}, we arrive at the result. 
\end{proof}

\begin{lemma}\label{l:lower}
Suppose $d \ge 2$. There exists 
$C>0$ such that for all $t>0$ and $x,y \in \R^d_+$,
	\begin{align*}
	&p(t,x,y) \ge C \left(1 \wedge \frac{x_d}{t^{1/\alpha}}\right)^{q} \left(1 \wedge \frac{ y_d}{t^{1/\alpha}}\right)^{q}
		p_\alpha(t, x, y) \big( |x-y|^{\alpha}\wedge t \big)\\
		&\qquad \times 
		\int_{(y_d \vee t^{1/\alpha}) \wedge (|x-y|/4)}^{|x-y|/2}A_{\beta_1, \beta_2,  \beta_3, \beta_4}(t,x,x+r\e_d)\,A_{\beta_1, \beta_2,  \beta_3, \beta_4}(t,x+r\e_d,y)\frac{dr}{r^{\alpha+1}}.
	\end{align*}
\end{lemma}
\begin{proof}    
Case (i): $(x_d \wedge y_d) \vee t^{1/\alpha} \ge |x-y|/4$. In this  case, we 
see from  \eqref{e:wtB-interior} that $A_{\beta_1, \beta_2, \beta_3, \beta_4}(t,x,y)$ is bounded  below by a  positive constant. Moreover,   since we have that $\sup_{s>0,\, z,w \in \R^d_+} A_{\beta_1,\beta_2,\beta_3,\beta_4}(s,z,w) \le c_1<\infty$,  we get  
\begin{align*}
	&|x-y|^\alpha\int_{(y_d \vee t^{1/\alpha}) \wedge( |x-y|/4)}^{|x-y|/2}A_{\beta_1,
		\beta_2,  \beta_3, \beta_4}(t,x,x+r\e_d)\,A_{\beta_1, \beta_2, \beta_3, \beta_4}(t,x+r\e_d,y)\frac{dr}{r^{\alpha+1}}\\
	& \le 
	c_1^2|x-y|^\alpha \int_{|x-y|/4}^{\infty} \frac{dr}{r^{\alpha+1}} =  
	4^\alpha c_1^2/\alpha.
\end{align*}
Therefore, we get the result  from Proposition \ref{p:HKE-lower-1} in this case.

Case (ii): $(x_d \wedge y_d) \vee t^{1/\alpha} < |x-y|/4$. 
For $r>0$, set $x(r):=x+r\e_d$ and
$$
K(r):= \left\{z=(\wt z, z_d) \in \R^d_+: |\wt z - \wt x| <\frac{r}2, \, z_d = x_d +r\right\}.
$$
Let
$$
K:=\left\{z \in \R^d_+: z \in K(r) \text{ for some } (y_d \vee t^{1/\alpha})\wedge \frac{|x-y|}4<r<\frac{|x-y|}2 \right\}.
$$
For any $z=(\wt z, x_d+r)\in K$, since 
\begin{equation}\label{e:lower-1}
	r\le |x-z|\le \frac{\sqrt{5}}2r\le \frac{\sqrt{5}}4|x-y|,
	\qquad (1-\frac{\sqrt{5}}4)|x-y|\le  |y-z|\le (1+\frac{\sqrt{5}}4)|x-y|,
\end{equation} 
we see that
\begin{align}\label{e:crit-1}
	A_{\beta_1,\beta_2,\beta_3,\beta_4}(t,x,z)  &\asymp  \Big(1  \wedge  \frac{x_d \vee t^{1/\alpha}}{r}\Big)^{\beta_1}  \Big(1 \wedge  \frac{(x_d+r) \vee t^{1/\alpha}}{r}\Big)^{\beta_2} \\  
	&\quad \times
	\L^{\beta_3} \Big( \frac{((x_d+r) \vee t^{1/\alpha} ) \wedge r}{(x_d \vee t^{1/\alpha}) \wedge r} \Big)\L^{\beta_4} \Big( \frac{r}{((x_d+r) \vee t^{1/\alpha} ) \wedge r} \Big) \nn\\
	& \asymp  \Big(1 \wedge  \frac{x_d \vee t^{1/\alpha}}{r}\Big)^{\beta_1}   \L^{\beta_3} \Big( \frac{r}{(x_d \vee t^{1/\alpha}) \wedge r} \Big)\nn\\
	& \asymp 	A_{\beta_1,\beta_2,\beta_3,\beta_4}(t,x,x(r)), \nn
\end{align}
and
\begin{align}\label{e:crit-2}
	A_{\beta_1,\beta_2,\beta_3,\beta_4}(t,z,y)	& \asymp  \Big(1 \wedge  \frac{y_d \vee t^{1/\alpha}}{|x-y|}\Big)^{\beta_1}\Big(1 \wedge  \frac{x_d+r}{|x-y|}\Big)^{\beta_2} \L^{\beta_3} \Big( \frac{x_d+r}{y_d \vee t^{1/\alpha}} \Big) \L^{\beta_4} \Big( \frac{|x-y|}{x_d+r} \Big)\\
	&\asymp 	A_{\beta_1,\beta_2,\beta_3,\beta_4}(t,x(r),y).\nn
\end{align}
Thus, using  the semigroup property, Proposition  \ref{p:HKE-lower-1} and \eqref{e:lower-1}, we get
\begin{align*}
         &p(t,x,y) \ge \int_{K} p(t/2,x,z) p(t/2,z,y)dz\\
	&\ge c_2^2  \left(1 \wedge \frac{x_d}{t^{1/\alpha}}\right)^{q}   \left(1 \wedge \frac{y_d}{t^{1/\alpha}}\right)^{q}\int_K \frac{tA_{\beta_1,\beta_2,\beta_3,\beta_4}(t,x,z)}{|x-z|^{d+\alpha}}\, \frac{tA_{\beta_1,\beta_2,\beta_3,\beta_4}(t,z,y)}{|y-z|^{d+\alpha}}  dz \\
	&\ge c_3t^2\left(1 \wedge \frac{x_d}{t^{1/\alpha}}\right)^{q} \left(1 \wedge \frac{y_d}{t^{1/\alpha}}\right)^{q}\\
	&\quad \times \int_{(y_d \vee t^{1/\alpha}) \wedge (|x-y|/4)}^{|x-y|/2} \frac{A_{\beta_1,\beta_2,\beta_3,\beta_4}(t,x,x(r))}{r^{d+\alpha}} \frac{A_{\beta_1,\beta_2,\beta_3,\beta_4}(t,x(r),y)}{|x-y|^{d+\alpha}}
 \int_{K(r)}  d\wt{z} \, dr  \\
	&\ge c_4\left(1 \wedge \frac{x_d}{t^{1/\alpha}}\right)^{q}   \left(1 \wedge \frac{y_d}{t^{1/\alpha}}\right)^{q} \frac{t}{|x-y|^{d+\alpha}}\\
	&\quad\times t\int_{(y_d \vee t^{1/\alpha}) \wedge (|x-y|/4)}^{|x-y|/2}A_{\beta_1, \beta_2,  \beta_3, \beta_4}(t,x,x(r))\,A_{\beta_1, \beta_2,  \beta_3, \beta_4}(t,x(r),y)\frac{dr}{r^{\alpha+1}}. 
\end{align*} 
\end{proof}

Note that, for any $x,y \in \R^d_+$ and $|x-y|/4\le r\le|x-y|/2$, 
\begin{align}\label{e:two-jump-compare-1}
	A_{\beta_1,\beta_2,\beta_3,\beta_4}(t,x,x+r\e_d) & \asymp  \Big(1  \wedge  \frac{x_d \vee t^{1/\alpha}}{|x-y|}\Big)^{\beta_1}   \L^{\beta_3} \Big( \frac{|x-y|}{(x_d \vee t^{1/\alpha}) \wedge |x-y|}\Big) \\
	&\asymp 	A_{\beta_1,\beta_2,\beta_3,\beta_4}(t,x,x+ \frac{|x-y|}{2} \e_d) \nn
\end{align}
and,  since $y_d \le x_d +4r$, 
\begin{align}\label{e:two-jump-compare-2}
&	A_{\beta_1,\beta_2,\beta_3,\beta_4}(t,x+r\e_d,y)	\asymp  \Big(1 \wedge  \frac{y_d \vee t^{1/\alpha}}{|x-y|}\Big)^{\beta_1} \L^{\beta_3} \Big( \frac{|x-y|}{(y_d \vee t^{1/\alpha}) \wedge |x-y|} \Big)\\
	&\asymp  A_{\beta_1,\beta_2,\beta_3,\beta_4}(t,x+ \frac{|x-y|}{2} \e_d,y)
	\asymp 	A_{\beta_1,\beta_2,\beta_3,\beta_4}(t,y,y+ \frac{|x-y|}{2} \e_d). \nn
\end{align}
In particular, we have
\begin{align}\label{e:two-jump-compare-3}
	&A_{\beta_1,\beta_2,\beta_3,\beta_4}(t,x,x+\frac{|x-y|}{2}\e_d) \,	A_{\beta_1,\beta_2,\beta_3,\beta_4}(t,x+\frac{|x-y|}{2}\e_d,y)\\
	&\asymp A_{\beta_1,\beta_2,\beta_3,\beta_4}(t,x,x+\frac{|x-y|}{2}\e_d)  \,	A_{\beta_1,\beta_2,\beta_3,\beta_4}(t,y,y+\frac{|x-y|}{2}\e_d) \nn \\
	&  \asymp A_{\beta_1, \beta_1, 0, \beta_3}(t,x,y) \L^{\beta_3} \Big( \frac{|x-y|}{((x_d \wedge y_d) \vee  t^{1/\alpha}) \wedge |x-y|} \Big). \nn
\end{align}

By \eqref{e:two-jump-compare-1}--\eqref{e:two-jump-compare-3},  we have  that for all $t>0$ and $x,y \in \R^d_+$,
\begin{align}\label{e:two-jump-region}
	&\int_{|x-y|/4}^{|x-y|/2}A_{\beta_1, \beta_2,  \beta_3, \beta_4}(t,x,x+r\e_d)\,A_{\beta_1, \beta_2,  \beta_3, \beta_4}(t,x+r\e_d,y)\frac{dr}{r^{\alpha+1}} \\
	& \asymp
	A_{\beta_1, \beta_2,  \beta_3, \beta_4}(t,x,x+\frac{|x-y|}{2}\e_d)\,A_{\beta_1, \beta_2,  \beta_3, \beta_4}(t,x+\frac{|x-y|}{2}\e_d,y)
	\int_{|x-y|/4}^{|x-y|/2} \frac{dr}{r^{\alpha+1}} \nn\\
	&\asymp |x-y|^{-\alpha}  A_{\beta_1,\beta_1,0,\beta_3}(t,x,y) \L^{\beta_3} \Big( \frac{|x-y|}{((x_d \wedge y_d) \vee t^{1/\alpha}) \wedge |x-y|}\Big). \nn
\end{align}

\begin{remark}\label{r:proofCase21}
	{\rm  
Here we give a proof of the comparability of \eqref{e:Case21a} and \eqref{e:Case21}.
Let $t>0$ and $x,y \in \oR^d_+$ be such that  
 $|x-y| > 6t^{1/\alpha}$. 
For any $z \in B(x+2^{-1}|x-y|\e_d, 4^{-1}|x-y|) $,   by using the triangle inequality several times, 
we have
\begin{align}\label{e:Case2-zd}
 z_d-t^{1/\alpha} \asymp   z_d \asymp x_d \wedge y_d +|x-y|  \asymp x_d \vee y_d +|x-y|
\end{align} 
and 
\begin{align}\label{e:Case2-dist}
	|x+t^{1/\alpha}\e_d - z| \asymp |y+t^{1/\alpha}\e_d - z| \asymp   |x-z| \asymp |y-z| \asymp   |x-y|.
\end{align}

For any $t>0$ and $x,y \in \oR^d_+$ with $|x-y|>6t^{1/\alpha}$, 
 if $x_d \wedge y_d \ge |x-y|/4$, 
 then by {\bf (A3)},  \eqref{e:comp-AB}, \eqref{e:wtB-interior}, \eqref{e:Case2-zd} 
 and  \eqref{e:Case2-dist},  
\begin{align*}
&	 t |x-y|^{d+\alpha}  \int_{B(x+2^{-1}|x-y|\e_d, 
	\,4^{-1}|x-y|) }  
J(x+t^{1/\alpha}\e_d,z)\, J(z,y+t^{1/\alpha}\e_d)dz\\
& \asymp \frac{t}{|x-y|^{d+\alpha}}  \int_{B(x+2^{-1}|x-y|\e_d, 
	\,4^{-1}|x-y|) }  dz \asymp \frac{t}{|x-y|^\alpha}\\
&\asymp 	\Big(1 \wedge  \frac{t}{|x-y|^\alpha}\Big)B_{\beta_1, \beta_1, 0, \beta_3}(x+t^{1/\alpha}\e_d,y+t^{1/\alpha}\e_d)\, \L^{\beta_3} \Big( \frac{|x-y|}{((x_d \wedge y_d)+ t^{1/\alpha}) \wedge |x-y|} \Big).
\end{align*} 
 Otherwise, 
 if $x_d \wedge y_d <|x-y|/4$, then $x_d \vee y_d \le x_d \wedge y_d + |x-y|<5|x-y|/4$ so that
 $$
 z_d-t^{1/\alpha} \asymp z_d \asymp |x-y| \quad \text{for $z \in B(x+2^{-1}|x-y|\e_d, 
 	\,4^{-1}|x-y|) $}
 $$
 by \eqref{e:Case2-zd}.   Using this,   \eqref{e:Case2-dist} and  \eqref{e:comp-AB}, we get that for $z \in B(x+2^{-1}|x-y|\e_d, 
 \,4^{-1}|x-y|) $,
 \begin{align}\label{e:Case2-x}
&B_{\beta_1,\beta_2,\beta_3,\beta_4}(x+t^{1/\alpha}\e_d,z) \asymp B_{\beta_1,\beta_2,\beta_3,\beta_4}(x+t^{1/\alpha}\e_d,z+t^{1/\alpha}\e_d) \\
&\asymp A_{\beta_1,\beta_2,\beta_3,\beta_4}(t,x,z)\asymp 	A_{\beta_1,\beta_2,\beta_3,\beta_4}(t,x,x + 2^{-1}|x-y| \e_d) \nn
 \end{align}
and
 \begin{align}\label{e:Case2-y}
	&B_{\beta_1,\beta_2,\beta_3,\beta_4}(z,y+t^{1/\alpha}\e_d) \asymp B_{\beta_1,\beta_2,\beta_3,\beta_4}(z+t^{1/\alpha}\e_d,y+t^{1/\alpha}\e_d) \\
	&\asymp A_{\beta_1,\beta_2,\beta_3,\beta_4}(t,z,z)\asymp 	A_{\beta_1,\beta_2,\beta_3,\beta_4}(t,x + 2^{-1}|x-y| \e_d,y). \nn
\end{align}
By  {\bf (A3)},   \eqref{e:Case2-dist}, \eqref{e:Case2-x}, \eqref{e:Case2-y}  and \eqref{e:two-jump-compare-3}, we arrive at
\begin{align*}
	&	 t |x-y|^{d+\alpha}  \int_{B(x+2^{-1}|x-y|\e_d, 
		\,4^{-1}|x-y|) }  
	J(x+t^{1/\alpha}\e_d,z)\, J(z,y+t^{1/\alpha}\e_d)dz\\
	& \asymp \frac{tA_{\beta_1,\beta_2,\beta_3,\beta_4}(t,x,x + 2^{-1}|x-y| \e_d) \,A_{\beta_1,\beta_2,\beta_3,\beta_4}(t,x + 2^{-1}|x-y| \e_d,y)}{|x-y|^{\alpha}} \\
	&\asymp  \frac{t}{|x-y|^\alpha} A_{\beta_1, \beta_1, 0, \beta_3}(t,x,y) \L^{\beta_3} \Big( \frac{|x-y|}{((x_d \wedge y_d) \vee  t^{1/\alpha}) \wedge |x-y|} \Big)\\
	&\asymp 	\Big( 1 \wedge  \frac{t}{|x-y|^\alpha} \Big)B_{\beta_1, \beta_1, 0, \beta_3}(x+t^{1/\alpha}\e_d,y+t^{1/\alpha}\e_d) \L^{\beta_3}  \Big( \frac{|x-y|}{((x_d \wedge y_d)+ t^{1/\alpha}) \wedge |x-y|} \Big).
\end{align*} 
 Hence,  \eqref{e:Case21a} and \eqref{e:Case21} are comparable.
 }
\end{remark}

Combining \eqref{e:two-jump-region}, Proposition \ref{p:HKE-lower-1} and Lemma \ref{l:lower}, we get the following result.

\begin{prop}\label{p:lower-1}
Suppose $d \ge 2$. There exists a constant 
$C>0$ such that for all $t>0$ and $x,y \in \R^d_+$,
	\begin{align*}
 		p(t,x,y)& \ge C \Big(1 \wedge \frac{x_d}{t^{1/\alpha}}\Big)^{q} \Big(1 \wedge \frac{ y_d}{t^{1/\alpha}}\Big)^{q}
		 p_\alpha(t, x, y)
		\Big[	A_{\beta_1, \beta_2, \beta_3, \beta_4}(t,x,y)\\
		& + 	  \Big(1 \wedge  \frac{t}{|x-y|^\alpha}\Big)A_{\beta_1, \beta_2,  \beta_3, \beta_4}(t,x,x+\frac{|x-y|}{2}\e_d)\,A_{\beta_1, \beta_2,  \beta_3, \beta_4}(t,x+\frac{|x-y|}{2}\e_d,y) \Big]\\
		&\asymp \Big(1 \wedge \frac{x_d}{t^{1/\alpha}}\Big)^{q} \Big(1 \wedge \frac{ y_d}{t^{1/\alpha}}\Big)^{q}	
		p_\alpha(t, x, y)
		\Big[	A_{\beta_1, \beta_2, \beta_3, \beta_4}(t,x,y)	\nn\\
		&+  \Big(1 \wedge  \frac{t}{|x-y|^\alpha}\Big)A_{\beta_1, \beta_1, 0, \beta_3}(t,x,y) \L^{\beta_3} \Big( \frac{|x-y|}{((x_d \wedge y_d) \vee  t^{1/\alpha}) \wedge |x-y|} \Big) \Big].
	\end{align*}
\end{prop}

\begin{remark}\label{r:lowerbound}
	{\rm If $\beta_2<\alpha+\beta_1$, then there is 
	$C>0$ such that for all $t>0$ and $x,y\in \R^d$,
		\begin{align*}
		 &\Big(1 \wedge  \frac{t}{|x-y|^\alpha}\Big)A_{\beta_1, \beta_1, 0, \beta_3}(t,x,y) \L^{\beta_3} \Big( \frac{|x-y|}{((x_d \wedge y_d) \vee  t^{1/\alpha}) \wedge |x-y|} \Big)\\
		 & \le C	A_{\beta_1, \beta_2, \beta_3, \beta_4}(t,x,y).
		\end{align*}
	Indeed, for $\eps:=(\alpha+\beta_1-\beta_2)/2>0$, using
  \eqref{e:slowly-varying}, we get that for all $t>0$ and $x,y \in \R^d$,  with $x_d\le y_d$, 
	\begin{align*}
		&\Big(1 \wedge  \frac{t}{|x-y|^\alpha}\Big)A_{\beta_1, \beta_1, 0, \beta_3}(t,x,y) \L^{\beta_3} \Big( \frac{|x-y|}{((x_d \wedge y_d) \vee  t^{1/\alpha}) \wedge |x-y|} \Big)\\
		&\le c_1 \Big(1\wedge \frac{t}{|x-y|^\alpha}\Big) \Big(1\wedge \frac{x_d \vee t^{1/\alpha}}{|x-y|}\Big)^{\beta_1-\eps } \Big(1 \wedge \frac{y_d \vee t^{1/\alpha}}{|x-y|}\Big)^{\beta_1 -\eps }\\
		&= c_1 \Big(1\wedge \frac{t^{1/\alpha}}{|x-y|}\Big)^{\eps }   \Big(1\wedge \frac{x_d \vee t^{1/\alpha}}{|x-y|}\Big)^{\beta_1-\eps }  \Big(1\wedge \frac{t^{1/\alpha}}{|x-y|}\Big)^{-\beta_1+\beta_2+\eps } \Big(1 \wedge \frac{y_d \vee t^{1/\alpha}}{|x-y|}\Big)^{\beta_1 -\eps }\\
			&\le  c_1 \Big(1\wedge \frac{x_d \vee t^{1/\alpha}}{|x-y|}\Big)^{\beta_1} \Big(1 \wedge \frac{y_d \vee t^{1/\alpha}}{|x-y|}\Big)^{\beta_2} \le c_1 A_{\beta_1,\beta_2,\beta_3,\beta_4}(t,x,y).
	\end{align*}
	Therefore, in view of Proposition \ref{p:HKE-lower-1}, Proposition \ref{p:lower-1} is 	relevant only if
	$d \ge 2$ and $\beta_2 \ge \alpha+\beta_1$.
	}
\end{remark}

\begin{lemma}\label{l:lower-2}
Suppose that  $d \ge 2$ and $\beta_2=\alpha+\beta_1$.	Then  for all $t>0$ and $x,y \in \R^d_+$,
	\begin{align*}
 		&\big( |x-y|^{\alpha} \wedge t \big) \int_{(x_d \vee y_d \vee t^{1/\alpha})\wedge (|x-y|/4)}^{|x-y|/2}     A_{\beta_1, \beta_2,  \beta_3, \beta_4}(t,x,x+r\e_d)\,A_{\beta_1, \beta_2,  \beta_3, \beta_4}(t,x+r\e_d,y)
\frac{dr}{r^{\alpha+1}}
\nn\\		
		&\asymp  \Big(1 \wedge  \frac{t}{|x-y|^\alpha}\Big) A_{\beta_1, \beta_1, 0, \beta_3+\beta_4+1}(t,x,y) \L^{\beta_3} \Big( \frac{|x-y|}{((x_d \wedge y_d) \vee t^{1/\alpha}) \wedge |x-y|} \Big).
	\end{align*}
\end{lemma}
\begin{proof}  Without loss of generality, we assume $x_d \le y_d$. 

Assume first that  $y_d \vee t^{1/\alpha} < |x-y|/4$. 
Using $\beta_2=\alpha+\beta_1$ and \eqref{e:crit-1}--\eqref{e:crit-2},
since $x_d \le y_d \vee t^{1/\alpha}$,  we get
\begin{align}\label{e:crit-3} 
		&\int_{y_d \vee t^{1/\alpha}}^{|x-y|/2}A_{\beta_1, 
		\beta_2, \beta_3, \beta_4}(t,x,x+r\e_d)\,A_{\beta_1, 
		\beta_2,\beta_3, \beta_4}
		(t,x+r\e_d,y)\frac{dr}{r^{\alpha+1}}\\
		&\asymp   \int_{y_d \vee t^{1/\alpha}}^{|x-y|/2}   \Big( \frac{x_d \vee t^{1/\alpha}}{r}\Big)^{\beta_1}    \Big(\frac{y_d \vee t^{1/\alpha}}{|x-y|}\Big)^{\beta_1} \Big( \frac{ r}{|x-y|}\Big)^{
		\alpha+\beta_1}\nn\\
		& \quad \times  \L^{\beta_3} \Big( \frac{r}{x_d \vee t^{1/\alpha}} \Big)  \L^{\beta_3} \Big( \frac{r}{y_d \vee t^{1/\alpha}} \Big) \L^{\beta_4} \Big( \frac{|x-y|}{r} \Big) \frac{dr}{r^{\alpha+1}}\nn\\
		&\asymp   
	\Big(\frac{x_d \vee t^{1/\alpha}}{|x-y|}\Big)^{\beta_1}\Big(\frac{y_d \vee t^{1/\alpha}}{|x-y|}\Big)^{\beta_1}	|x-y|^{-\alpha}     \nn\\
		& \quad \times \int_{y_d \vee t^{1/\alpha}}^{|x-y|/2} \L^{\beta_3} \Big( \frac{r}{x_d \vee t^{1/\alpha}} \Big)  \L^{\beta_3} \Big(\frac{r}{y_d \vee t^{1/\alpha}} \Big) \L^{\beta_4}\Big(\frac{|x-y|}{r} \Big) \frac{dr}{r}.\nn
\end{align}
By a change of the variables and Lemma \ref{cal:3}, we see that 
\begin{align}\label{e:crit-4}
	 &\int_{y_d \vee t^{1/\alpha}}^{|x-y|/2} \L^{\beta_3} \Big( \frac{r}{x_d \vee t^{1/\alpha}} \Big)  \L^{\beta_3} \Big( \frac{r}{y_d \vee t^{1/\alpha}} \Big) \L^{\beta_4} \Big( \frac{|x-y|}{r} \Big) \frac{dr}{r} \\
	 &= \int_{4(y_d \vee t^{1/\alpha})/|x-y|}^{2}\L^{\beta_3} \Big(\frac{|x-y|s}{4(x_d \vee t^{1/\alpha})} \Big)  \L^{\beta_3} \Big( \frac{|x-y|s}{4(y_d \vee t^{1/\alpha})} \Big) \L^{\beta_4} \Big( \frac{4}{s} \Big) \frac{ds}{s}\nn\\
	 & \asymp \L^{\beta_3} \Big( \frac{|x-y|}{4(x_d \vee t^{1/\alpha})} \Big)  \L^{\beta_3+\beta_4+1} \Big( \frac{|x-y|}{4(y_d \vee t^{1/\alpha})} \Big). \nn
\end{align}
By \eqref{e:crit-3}--\eqref{e:crit-4}, 
since $y_d \vee t^{1/\alpha}<|x-y|/4$, we obtain
\begin{align}\label{e:crit-5}
& \big( |x-y|^{\alpha} \wedge t \big) 
\int_{y_d \vee t^{1/\alpha}}^{|x-y|/2}A_{\beta_1, 
		\beta_2, \beta_3, \beta_4}(t,x,x+r\e_d)\,A_{\beta_1, 
		\beta_2, \beta_3, \beta_4}
		(t,x+r\e_d,y)\frac{dr}{r^{\alpha+1}}\\
&\asymp   \Big(1 \wedge  \frac{t}{|x-y|^\alpha}\Big)\Big(\frac{x_d \vee t^{1/\alpha}}{|x-y|}\Big)^{\beta_1} \Big(\frac{y_d \vee t^{1/\alpha}}{|x-y|}\Big)^{\beta_1}\L^{\beta_3} \Big( \frac{|x-y|}{4(x_d \vee t^{1/\alpha})} \Big)  \L^{\beta_3+\beta_4+1} \Big( \frac{|x-y|}{4(y_d \vee t^{1/\alpha})} \Big). \nn
\end{align}
If $y_d \vee t^{1/\alpha} \ge |x-y|/4$, then
$\L\big(\frac{ |x-y|}{( y_d \vee t^{1/\alpha}) \wedge |x-y|} \big) \asymp 1$.
Thus, by \eqref{e:two-jump-region},   
\begin{align}\label{e:crit-6}
	&  \big( |x-y|^{\alpha} \wedge t \big) 
	\int_{|x-y|/4}^{|x-y|/2}A_{\beta_1, 
		\beta_2, \beta_3, \beta_4}(t,x,x+r\e_d)\,A_{\beta_1, 
		\beta_2,  \beta_3, \beta_4}
		(t,x+r\e_d,y)\frac{dr}{r^{\alpha+1}}\\
	&\asymp    \Big(1 \wedge  \frac{t}{|x-y|^\alpha}\Big) A_{\beta_1,\beta_1,0,\beta_3+\beta_4+1}(t,x,y) \L^{\beta_3} \Big( \frac{|x-y|}{(x_d \vee  t^{1/\alpha}) \wedge |x-y|} \Big). \nn
\end{align}
The proof is now complete. 
\end{proof}

Combining Proposition \ref{p:HKE-lower-1},  Lemma \ref{l:lower} and Lemma \ref{l:lower-2}. 
we get the following

\begin{prop}\label{p:lower-2}
Suppose that $d \ge 2$ and $\beta_2=\alpha+\beta_1$.	 
There exists a constant $C>0$ such that for all $t>0$ and $x,y \in \R^d_+$,
	\begin{align*}
 		p(t,x,y)& \ge C \Big(1 \wedge \frac{x_d}{t^{1/\alpha}}\Big)^{q} \Big(1 \wedge \frac{ y_d}{t^{1/\alpha}}\Big)^{q}
		p_\alpha(t, x, y)
		\Big[ A_{\beta_1, \beta_2, \beta_3, \beta_4}(t,x,y) \nn\\		
		&\;\; +  \Big(1 \wedge  \frac{t}{|x-y|^\alpha}\Big) A_{\beta_1, \beta_1, 0, \beta_3+\beta_4+1}(t,x,y) \L^{\beta_3} \Big( \frac{|x-y|}{((x_d \wedge y_d) \vee t^{1/\alpha}) \wedge |x-y|} \Big)\Big].
	\end{align*}
\end{prop}

\begin{remark}
	{\rm In this section, for $\rY$, {\bf (A3)}(II) is only used to get {\bf (IUBS)}.  Therefore, the results of Propositions
		\ref{p:HKE-lower-1},  \ref{p:lower-1} and \ref{p:lower-2} hold under weaker assumptions {\bf (A1)}, {\bf (A3)}(I), {\bf (A4)} and {\bf (IUBS)}.	
 }
\end{remark}

\section{Sharp heat kernel upper bounds}\label{s:sharpupper}
 
In this section we prove the sharp heat kernel upper bounds. 
The key results are Theorem \ref{t:UHK} and its Corollary \ref{c:UHK} which deals with the 
cases $d=1$ or   $\beta_2<\alpha+\beta_1$, 
and Theorem \ref{t:UHK2} which deals with the case 
$d\ge 2$ and  $\beta_2\ge \alpha+\beta_1$.

Recall that $U(r)=\{x=(\wt{x}, x_d)\in \R^d:\, |\wt{x}|<r/2,\, 0\le x_d<r/2\}$ 
for  $r>0$, $d \ge 2$ and $U(r)=[0, r/2)$ for $d=1$, and that the function $ A_{\beta_1,\beta_2, \beta_3, \beta_4}(t,x,y)$ is defined by  \eqref{def:wtB2}. 
We also recall that  we use 
$p(t,x,y)$ for both $\bar p(t,x,y)$ and  $p^\kappa(t,x,y)$.
We remind the readers that,
for an open set $D \subset \oR^d_+$ relative to the  topology on $\oR^d_+$,
 $\tau_D=\bar \tau_D= \inf\{t>0: \rY_t \notin D\}$ when we consider $\rY$, and 
 $\tau_D=\tau^\kappa_D=\inf\{t>0: Y^\kappa_t \notin D \cap \R^d_+\}$  when we consider  $Y^\kappa$.

\begin{lemma}\label{l:first-term}
Let  $b_1,b_3\ge0$ be constants with $b_1>0$ if $b_3>0$. 
Suppose that there exists a constant $
 C_0 >0$ such that 	for all $t>0$ and $z,y \in \R^d_+$,
	\begin{equation}\label{e:first-term-ass}
	p(t,z,y) \le C_0  \left(1 \wedge \frac{z_d}{t^{1/\alpha}}\right)^{q}  
	 A_{b_1,0, b_3, 0}(t,z,y)
	p_\alpha(t, x, y).
	\end{equation}
	Then there exists 
 $C=C(C_0) >0$  such that for all $t>0$ and $x=(\wt 0, x_d) \in \R^d_+$ with $x_d\le 2^{-5}$,
\begin{equation}\label{e:first-term}
\P_x(\tau_{U(1)}<t< \zeta) \le  C
t\left(1 \wedge \frac{x_d}{t^{1/\alpha}}\right)^{q}  (x_d \vee  t^{1/\alpha})^{b_1} \L^{b_3}  \Big( \frac{1}{x_d \vee t^{1/\alpha}} \Big).
	\end{equation}
\end{lemma}
\begin{proof} 
We first note that \eqref{e:first-term-ass} implies that for all $t>0$ and $y \in \R^d_+$,
\begin{align}\label{e:nestf1}
	&\P_y(|Y_t-y|>2^{-3}, \, t<\zeta)= \int_{z \in \R^d_+,\,|z-y| > 2^{-3}} p(t,y,z)dz \\
	& \le c_1 t\left(1 \wedge \frac{y_d}{t^{1/\alpha}}\right)^{q} 	\int_{z \in \R^d_+, \, |z-y| > 2^{-3}} \Big( \frac{y_d \vee t^{1/\alpha}}{|z-y|} \wedge 1 \Big)^{b_1} \L^{b_3} \Big( \frac{|z-y|}{ y_d \vee t^{1/\alpha}}\Big)	\frac{dz}{|z-y|^{d+\alpha}}\nn \\
	&\le c_2t\left(1 \wedge \frac{y_d}{t^{1/\alpha}}\right)^{q}   (y_d \vee t^{1/\alpha})^{b_1}	\L^{b_3} \Big( \frac{2^{-3}}{ y_d \vee t^{1/\alpha}} \Big)	\int_{z \in \R^d_+, \, |z-y| > 2^{-3}} \frac{dz}{|z-y|^{d+\alpha/2+b_1}}\nn\\
	&\le  c_3t\left(1 \wedge \frac{y_d}{t^{1/\alpha}}\right)^{q}   (y_d \vee t^{1/\alpha})^{b_1}	\L^{b_3} \Big( \frac{1}{ y_d \vee t^{1/\alpha}}\Big), \nn
\end{align}
where in the first inequality above we used  \eqref{e:first-term-ass} and Lemma \ref{l:kill-log-2}, 
and in the second we used \eqref{e:slowly-varying-2}.

By Proposition \ref{p:upper-heatkernel}  (see also Remark \ref{r:conti}), we have
\begin{equation}\label{e:regular-path}
	\sup_{s \le t,\, y \in \R^d_+}	\P_y\big( |Y_s - y|\ge 2^{-2}, \, s<\zeta\big) \le c_1 \sup_{s \le t, \, y \in \R^d_+} \int_{z \in \R^d_+, \, |z-y| \ge 2^{-2}} \frac{s}{|z-y|^{d+\alpha}}dz \le c_2t.
\end{equation}
If $t\ge 1/(2c_2)$, 
then  \eqref{e:first-term} follows from Corollary \ref{c:life}. 

Let $t<1/(2c_2)$.
By the strong Markov property and \eqref{e:regular-path}, we have 
\begin{align*}
       &\P_x\big(\tau_{U(1)}<t< \zeta, \, |Y_t-Y_{\tau_{U(1)}}|\ge 2^{-2}\big)
	=  \E_x\left[ \P_{Y_{\tau_{U(1)}}}  
		\left( |Y_{t- \tau_{U(1)} } - Y_0 |  \ge 2^{-2}\right) :  
	\tau_{U(1)}<t< \zeta \right] \nn\\
	& \le \P_x(\tau_{U(1)}<t< \zeta)  \sup_{s \le t,\, y \in \R^d_+}	\P_y( |Y_s - y|\ge 2^{-2}, \, s<\zeta) 
	 \le 2^{-1}\P_x(\tau_{U(1)}<t< \zeta).
\end{align*}
Thus, 
\begin{align}\label{e:first-term-1}
	\P_x(\tau_{U(1)}<t< \zeta)&=2\big(\P_x(\tau_{U(1)}<t< \zeta)-2^{-1}\P_x(\tau_{U(1)}<t< \zeta)\big) \\
	&\le 2 \big(\P_x(\tau_{U(1)}<t< \zeta)-\P_x(\tau_{U(1)}<t< \zeta, \, |Y_t-Y_{\tau_{U(1)}}|\ge 2^{-2}) \big)\nn\\
	&=  2\P_x(\tau_{U(1)}<t< \zeta, \, |Y_t-Y_{\tau_{U(1)}}|<2^{-2}). \nn
\end{align}
Note that by the triangle inequality, for any $y \in \R^d_+\setminus U(1)$ and 
$z \in B(y,2^{-2})$, we have $|z-x| \ge |y-x|-|y-z| > 15/32-1/4>2^{-3}$. 
Therefore using \eqref{e:nestf1} and \eqref{e:first-term-1}, we have
\begin{align*}
	&\P_x(\tau_{U(1)}<t< \zeta)\le  2\P_x(\tau_{U(1)}<t< \zeta, \, |Y_t-Y_{\tau_{U(1)}}|<2^{-2})\\
	& \le 2\P_x(|Y_t-x|>2^{-3}, \, t<\zeta)\le  c_4t\left(1 \wedge \frac{x_d}{t^{1/\alpha}}\right)^{q}   (x_d \vee t^{1/\alpha})^{b_1}	\L^{b_3} \Big(\frac{1}{ x_d \vee t^{1/\alpha}}\Big).
\end{align*}
The proof is complete.\end{proof}

 Note that for any $t,k,r>0$ and $a \ge 0$,
\begin{equation}\label{e:equiv-form}
	\Big(1 \wedge \frac{r}{t^{1/\alpha}} \Big)^{a} (r \vee t^{1/\alpha})^{k} = 	r^{k} \Big(1 \wedge \frac{r}{t^{1/\alpha}} \Big)^{a-k}. 
\end{equation}

\begin{lemma}\label{l:UHK-case1-induction}
There exists a constant $C>0$ such that for all $t>0$ and $x,y \in \R^d_+$,
	\begin{equation*}
		p(t,x,y) \le C \left(1 \wedge \frac{x_d \wedge y_d}{t^{1/\alpha}}\right)^{q}  
	A_{\beta_1,0, \beta_3, 0}(t,x,y)
		p_\alpha(t, x, y).
	\end{equation*}
\end{lemma}
\begin{proof}    
By Proposition \ref{p:UHK-rough}, the lemma holds for $\beta_1=0$.

We assume $\beta_1>0$ and set $a_n= \beta_1 \wedge \frac{n\alpha}{2}$ for $n\ge0$. 
Below, we prove by induction that for any $n \ge 0$, there exists a constant $C>0$ such that for all $t>0$ and $x,y \in \R^d_+$,
	\begin{align}\label{e:UHK-induction-0}
	p(t,x,y) \le C \Big(1 \wedge \frac{x_d \wedge y_d}{t^{1/\alpha}}\Big)^{q}  
	A_{a_n,0, \beta_3, 0}(t,x,y)
	p_\alpha(t, x, y).
\end{align} The lemma is a direct consequence of \eqref{e:UHK-induction-0}.  

By Proposition \ref{p:UHK-rough} and the fact that the logarithmic term in $A_{\beta_1,0, \beta_3, 0}(t,x,y)$ is always larger than 1, \eqref{e:UHK-induction-0} holds for $n=0$.
Suppose \eqref{e:UHK-induction-0} holds for $n-1$. By symmetry and  \eqref{e:kernel-scaling}, we can assume without loss of generality that $x_d \le y_d$, $\wt x=0$ and  $|x-y|=5$. If $t > 1$ or $x_d > 2^{-5}$, then  \eqref{e:UHK-induction-0} follows from Proposition \ref{p:UHK-rough} and \eqref{e:wtB-interior}.

Let $t \le 1$ and $x_d\le 2^{-5}$. Then $y_d\le x_d+|x-y|\le 4+2^{-5}$ by the triangle inequality. Our goal is to show that there exists a constant $c_1>0$ independent of $t,x,y$ such that 
\begin{equation}\label{e:UHK-case1}
	p(t,x,y) \le c_1t\Big(1 \wedge \frac{x_d }{t^{1/\alpha}}\Big)^{q}   (x_d \vee  t^{1/\alpha})^{a_n} \L^{\beta_3}\Big(\frac{ y_d \vee t^{1/\alpha} }{ x_d \vee t^{1/\alpha}} \Big).
\end{equation}

Set $V_1=U(1)$, $V_3=B(y,2) \cap \oR^d_+$ and $V_2=\oR^d_+ \setminus (V_1 \cup V_3)$. Similarly to \eqref{e:UHK-rough-0} and \eqref {e:UHK-rough-1}, we get from Proposition \ref{p:upper-heatkernel} and the triangle inequality that 
\begin{equation}\label{e:UHK-case1-2}
	\sup_{s \le t, \, z \in V_2} p(s,z,y) \le c_2 \sup_{s \le t, \, z \in \R^d_+, |z-y|\ge 2} \frac{s}{|z-y|^{d+\alpha}} \le 2^{-d-\alpha}c_2t
\end{equation}
and
\begin{equation}\label{e:UHK-case1-3}
	{\rm dist}(V_1, V_3) \ge \sup_{u\in V_1, \, w \in V_3} (4-|x-u|-|y-w|) \ge 1.
\end{equation}  
By  the induction hypothesis,  
condition \eqref{e:first-term-ass} in Lemma  \ref{l:first-term} holds  with  $b_1= a_{n-1}$ and $b_3=\beta_3$.
Thus, since $ a_{n}-a_{n-1} \le \alpha/2$, we get from Lemma  \ref{l:first-term} and \eqref{e:UHK-case1-2} that
\begin{align}\label{e:UHK-case1-term1}
	&\P_x(\tau_{V_1}<t< \zeta) \sup_{s \le t, \, z \in V_2} p(s,z,y)\\
	&\le 
	c_3 t  \Big(1 \wedge \frac{x_d}{t^{1/\alpha}}\Big)^{q }   (x_d \vee t^{1/\alpha})^{a_{n-1}}    (t^{1/\alpha})^{\alpha/2} t^{1/2} \L^{\beta_3} \Big( \frac{1}{ t^{1/\alpha}} \Big)\,\nn\\
	&\le 
	c_3 t  \Big(1 \wedge \frac{x_d}{t^{1/\alpha}}\Big)^{q }   (x_d \vee t^{1/\alpha})^{a_{n-1}}    (t^{1/\alpha})^{a_{n}-a_{n-1}} \,
	\Big(\sup_{s\le 1}s^{1/2} \L^{\beta_3} \Big( \frac{1}{ s^{1/\alpha}} \Big)\Big)  \nn\\
	&  
\le c_4t\Big(1 \wedge \frac{x_d }{t^{1/\alpha}}\Big)^{q} (x_d \vee t^{1/\alpha})^{a_n }. \nn
\end{align}

In order to apply Lemma \ref{l:general-upper-2} and get the desired result, it remains  to bound $\int_0^t\int_{V_3} \int_{V_1} p^{V_1}(t-s, x, u) \sB(u,w) p(s, y,w) du dw ds$. We consider the following two cases separately.

\smallskip

\textbf{(Case 1)} $q \ge \alpha + a_n$ and $10x_d< t^{1/\alpha}$. 

 Pick $\eps\in (0,\beta_1)$ such that  $q<\alpha+\beta_1-\eps$. Using {\bf (A3)}(II),    
 Lemmas \ref{l:kill-log}(i)--(ii), \ref{l:kill-log-2} (see Remark \ref{r:kill-log}), and
\eqref{e:UHK-case1-3}, we see that for all $u \in V_1$ and $w \in V_3$,
\begin{align}\label{e:case1-B}
	\sB(u,w) \le c_5 B_{\beta_1-\eps, 0, 0, 0}
	(u,w) \le c_6 \left(  \frac{u_d}{|u-w|} \right)^{\beta_1-\eps} \le c_6 u_d^{\beta_1-\eps}.
\end{align}
By  \eqref{e:case1-B} and Lemma \ref{l:6.10},  we have
\begin{align}\label{e:UHK-case1-term2}
	&\int_0^t\int_{V_3} \int_{V_1} p^{V_1}(t-s, x, u) \sB(u,w) p(s, y,w) du dw ds\\
	&\le c_6\int_0^t \int_{V_1} p^{V_1}(t-s, x, u) u_d^{\beta_1-\eps} du   \int_{V_3} p(s, y,w)  dw ds\le c_6\int_0^\infty \int_{V_1} p^{V_1}(
		s, x, u) u_d^{\beta_1-\eps} du ds\nn\\
		&\le c_7 x_d^q = c_7 t\left(\frac{x_d }{t^{1/\alpha}}\right)^{q} (t^{1/\alpha})^{q-\alpha}\le c_8 t\left(1 \wedge \frac{x_d }{t^{1/\alpha}}\right)^{q} (x_d \vee t^{1/\alpha})^{a_n}. \nn
\end{align}
In the last inequality above, we used the facts that $q-\alpha \ge a_n $, $10x_d< t^{1/\alpha}$ and $x_d \vee t^{1/\alpha} \le 1$. 
Now, using Lemma  \ref{l:general-upper-2},  \eqref{e:UHK-case1-term1} and \eqref{e:UHK-case1-term2}, we get  \eqref{e:UHK-case1} in this case.

\medskip

\textbf{(Case 2)} $q < \alpha + a_n $ or $10x_d\ge t^{1/\alpha}$. 

By {\bf (A3)}(II),   
\eqref{e:UHK-case1-3}, Lemma \ref{l:kill-log}(ii) and Lemma \ref{l:kill-log-2}  
(see Remark \ref{r:kill-log}), 
 it holds that  for all $u \in V_1$ and $w\in V_3$,
\begin{align}\label{e:case2-B}
	\sB(u,w) \le c_9B_{\beta_1, 0, \beta_3, 0}(u,w) \le c_{10}\Big(  \frac{u_d}{|u-w|} \Big)^{\beta_1}\L^{\beta_3} \Big( \frac{w_d \wedge |u-w|}{u_d \wedge |u-w|}\Big) \le c_{10} u_d^{\beta_1}\L^{\beta_3} \Big( \frac{w_d}{u_d} \Big). 
\end{align}
By \eqref{e:slowly-varying-2} (or using that $u\mapsto u^{\beta_1}\L^{\beta_3}(t/u)$ is almost increasing) 
and Corollary \ref{c:life}, since  $ \beta_1>0$ and $a_n \le \beta_1$, we get that for any $0<s < t$ and $w \in V_3$, 
\begin{align}\label{e:case2-1}
	& \int_{u \in V_1:u_d < x_d} p(s, x, u) u_d^{\beta_1} \L^{\beta_3} \Big( \frac{ w_d}{u_d} \Big) du \\
	& \le  \int_{u \in V_1:u_d < x_d \vee s^{1/\alpha}} p(s, x, u) u_d^{\beta_1} \L^{\beta_3} \Big( \frac{ w_d}{u_d} \Big) du \nn \\
	&\le c_{11} (x_d \vee s^{1/\alpha})^{\beta_1} \L^{\beta_3} \Big( \frac{ w_d}{x_d \vee s^{1/\alpha}} \Big) \int_{u \in V_1:u_d < x_d \vee s^{1/\alpha}} p(s, x, u) du \nn\\
	 	&\le c_{11}\P_x(\zeta>s) (x_d \vee s^{1/\alpha})^{a_n} \L^{\beta_3} \Big( \frac{ w_d}{x_d \vee s^{1/\alpha}} \Big)  \nn\\
	 	& \le c_{12} \Big(1 \wedge \frac{x_d}{s^{1/\alpha}} \Big)^{q} (x_d \vee s^{1/\alpha})^{a_n} \L^{\beta_3} \Big( \frac{ w_d}{x_d \vee s^{1/\alpha}} \Big). \nn
\end{align}
Next, using the induction hypothesis  and  Lemma \ref{cal:2}, 
since $a_n  \le \beta_1$ and $a_n<\alpha+ a_{n-1}$, we get  that for any $0<s< t$ and $w \in V_3$,
\begin{align}\label{e:case2-2}
	& \int_{u \in V_1:u_d \ge x_d} p(s, x, u) u_d^{\beta_1} \L^{\beta_3} \Big( \frac{ w_d}{u_d} \Big) du \\
	&\le c_{13} \Big(1 \wedge \frac{x_d}{s^{1/\alpha}} \Big)^{q}  \int_{u \in V_1:u_d \ge x_d}  \Big( \frac{x_d \vee s^{1/\alpha}}{|x-u|} \wedge 1 \Big)^{a_{n-1}}  
	 \left(s^{-d/\alpha} \wedge \frac{s}{|x-u|^{d+\alpha}}\right) \nn\\
	&\hspace{3cm}\times  \L^{\beta_3} \Big(\frac{  |x-u|}{(x_d \vee s^{1/\alpha}) \wedge |x-u|} \Big)u_d^{a_n} \L^{\beta_3} \Big( \frac{ w_d}{u_d} \Big)du \nn\\
	&\le c_{14} \Big(1 \wedge \frac{x_d}{s^{1/\alpha}} \Big)^{q} (x_d \vee s^{1/\alpha})^{ a_n} \L^{\beta_3} \Big( \frac{ w_d}{x_d \vee s^{1/\alpha}} \Big). \nn
\end{align}
Similarly, again splitting 
 the integration into two parts $w_d< y_d$ and $w_d\ge y_d$, and using the induction hypothesis  and 
Lemma \ref{cal:2} again,  we also get that for any $0<s<t$,
\begin{align}\label{e:case2-3}
	&\int_{V_3}  p(s, y,w) \L^{\beta_3} \Big( \frac{ w_d}{ x_d \vee (t-s)^{1/\alpha} } \Big) dw  \\
	&\le c_{15} \Big(1 \wedge \frac{y_d}{s^{1/\alpha}} \Big)^{q}  \L^{\beta_3} \Big(\frac{ y_d \vee s^{1/\alpha}}{x_d \vee (t-s)^{1/\alpha} } \Big)
	\le c_{15}  \L^{\beta_3} \Big( \frac{ y_d \vee t^{1/\alpha}}{x_d \vee (t-s)^{1/\alpha} } \Big). \nn
\end{align}

By \eqref{e:case2-B}--\eqref{e:case2-3} and \eqref{e:equiv-form}, we have
\begin{align}\label{e:UHK-case1-7}
	&\int_0^t\int_{V_3} \int_{V_1} p^{V_1}(t-s, x, u) \sB(u,w) p(s, y,w) du dw ds\\
	&\le c_{10}\int_0^t \int_{V_3} p(s, y,w)   \int_{V_1} p(t-s, x, u) u_d^{\beta_1} \L^{\beta_3} \Big( \frac{w_d}{u_d} \Big) du   dw ds 	\nn\\
	&\le c_{16}   \int_0^t \Big(1 \wedge \frac{x_d}{(t-s)^{1/\alpha}} \Big)^{q} (x_d \vee (t-s)^{1/\alpha})^{a_n} \int_{V_3} p(s, y,w) \L^{\beta_3} \Big( \frac{ w_d}{x_d \vee (t-s)^{1/\alpha}} \Big)   dw ds \nn\\
	&\le c_{17}x_d^{a_n} \int_0^t  \Big(1 \wedge \frac{x_d}{(t-s)^{1/\alpha}} \Big)^{q- a_n }  \L^{\beta_3} \Big( \frac{ y_d \vee t^{1/\alpha}}{x_d \vee (t-s)^{1/\alpha} } \Big)     ds=:I. \nn
\end{align}
When $q<\alpha+a_n$, we get from Lemma \ref{cal:00} that
\begin{align}
I \le  c_{18}tx_d^{a_n} \Big(1 \wedge \frac{x_d}{t^{1/\alpha}} \Big)^{q-a_n}  \L^{\beta_3} \Big( \frac{ y_d \vee t^{1/\alpha}}{x_d \vee t^{1/\alpha} } \Big) .
\end{align}
When $10x_d \ge t^{1/\alpha}$, we also get from Lemma \ref{cal:00} that
\begin{align}\label{e:case2-4}
	I &\le c_{17}x_d^{a_n} \int_0^t  \L^{\beta_3} \Big( \frac{ y_d \vee t^{1/\alpha}}{x_d \vee (t-s)^{1/\alpha} } \Big) ds \le c_{19}tx_d^{a_n}  \L^{\beta_3} \Big( \frac{ y_d \vee t^{1/\alpha}}{x_d \vee t^{1/\alpha} } \Big) \\
	&\le 10^{q-a_n}c_{19}tx_d^{a_n} \Big(1 \wedge \frac{x_d}{t^{1/\alpha}} \Big)^{q- a_n }  \L^{\beta_3} \Big( \frac{ y_d \vee t^{1/\alpha}}{x_d \vee t^{1/\alpha} } \Big). \nn
\end{align}
Now, \eqref{e:UHK-case1}  follows from \eqref{e:equiv-form},  \eqref{e:UHK-case1-term1}, \eqref{e:UHK-case1-7}, \eqref{e:case2-4} and Lemma \ref{l:general-upper-2}. The proof is complete.
\end{proof}

Now we use Lemma \ref{l:UHK-case1-induction} to
improve the bound in \eqref{e:UHK-case1-term1}.

\begin{lemma}\label{l:firstpart}
	Let $0<t \le 1$ and $x,y \in \R^d_+$ be such that $\wt x =\wt 0$,  $x_d \le 2^{-5}$ 
	and $|x-y|=5$. Set $V_1=U(1)$, $V_3=B(y,2) \cap \oR^d_+$ and $V_2=\oR^d_+ \setminus (V_1 \cup V_3)$. There exists a constant $C>0$ independent of $t,x$ and $y$ such that
\begin{align*}
	&\P_x(\tau_{V_1}<t< \zeta) \sup_{s \le t, \, z \in V_2} p(s,z,y)\\
	&\le Ct^2 \Big(1 \wedge \frac{x_d}{t^{1/\alpha}}\Big)^{q} \Big( 1 \wedge \frac{ y_d}{t^{1/\alpha}}\Big)^{q} (x_d \vee t^{1/\alpha})^{\beta_1}   (y_d \vee t^{1/\alpha})^{\beta_1} \L^{\beta_3} \Big(\frac{1}{ x_d \vee t^{1/\alpha}} \Big)\L^{\beta_3} \Big(\frac{1}{ y_d \vee t^{1/\alpha}} \Big).
\end{align*}
\end{lemma}
\begin{proof} 
Note that, since  $y_d\le x_d+|x-y|\le 5+2^{-5}<6$, for any $0<s  \le 1$ and  $z \in \R^d_+$ with $|z-y|\ge 2$,
\begin{align}\label{e:first-kill-log-00}
(y_d\vee s^{1/\alpha})  \wedge |z-y| \ge \frac13(y_d\vee s^{1/\alpha}).
\end{align}
By  Lemmas \ref{l:UHK-case1-induction} and \ref{l:kill-log-2}, and applying \eqref{e:first-kill-log-00},
\begin{align}\label{e:first-kill-log-0}
	&\sup_{s \le t, \, z \in V_2} p(s,z,y) \le c_1  \sup_{s \le t, \, z \in \R^d_+, |z-y|\ge 2}  \Big( 1 \wedge \frac{z_d \wedge y_d}{s^{1/\alpha}}\Big)^{q}  
	A_{\beta_1,0,\beta_3,0}(s,z,y) 
	\frac{s}{|z-y|^{d+\alpha}} \\
&\le c_1  \sup_{s \le t, \, z \in \R^d_+, |z-y|\ge 2}  \Big( 1 \wedge \frac{ y_d}{s^{1/\alpha}}\Big)^{q}  
\Big(\frac{y_d\vee s^{1/\alpha}}{|z-y|} \Big)^{\beta_1} \L^{\beta_3} \Big(\frac{ 3 |z-y|}{y_d\vee s^{1/\alpha}} \Big)
	\frac{s}{|z-y|^{d+\alpha}}. \nn
\end{align}
Using that $u\mapsto u^{\beta_1}\L^{\beta_3}(1/u)$ is almost increasing, we have that  for any $0<s \le t \le 1$ and  $z \in \R^d_+$ with $|z-y|\ge 2$,
\begin{align}\label{e:first-kill-log-1}
	&\Big( 1 \wedge \frac{ y_d}{s^{1/\alpha}}\Big)^{q}  
\Big(\frac{y_d\vee s^{1/\alpha}}{|z-y|} \Big)^{\beta_1} \L^{\beta_3} \Big(\frac{ 3 |z-y|}{y_d\vee s^{1/\alpha}} \Big)
	\frac{s}{|z-y|^{d+\alpha}} \\ 
	& \le c_2
	 \Big( 1 \wedge \frac{ y_d}{s^{1/\alpha}}\Big)^{q}     
	(y_d \vee  s^{1/\alpha})^{\beta_1}  \L^{\beta_3} \Big(\frac{1}{ y_d \vee s^{1/\alpha}} \Big)	\frac{s}{|z-y|^{d+\alpha}}
	\nn\\
	&\le 2^{-d-\alpha}c_2s \Big( 1 \wedge \frac{ y_d}{s^{1/\alpha}}\Big)^{q}     
	(y_d \vee  s^{1/\alpha})^{\beta_1}  \L^{\beta_3} \Big(\frac{1}{ y_d \vee s^{1/\alpha}} \Big). \nn
\end{align}

Let $\eps>0$ be such that $q<\alpha+\beta_1-\eps$.
Using \eqref{e:equiv-form}, we see that 
\begin{align*}
	&s \Big( 1 \wedge \frac{ y_d}{s^{1/\alpha}}\Big)^{q}     
	(y_d \vee  s^{1/\alpha})^{\beta_1}  \L^{\beta_3} \Big(\frac{
	  1 }{ y_d \vee s^{1/\alpha}} \Big)\nn\\
	& = sy_d^{\beta_1-\eps} \Big( 1 \wedge \frac{ y_d}{s^{1/\alpha}}\Big)^{q-\beta_1+\eps}  	(y_d \vee  s^{1/\alpha})^{\eps}     \L^{\beta_3} \Big(\frac{1 }{ y_d \vee s^{1/\alpha}} \Big)
	\nn\\
		& = s^{(\alpha+\beta_1-q-\eps)/\alpha}y_d^{\beta_1-\eps} (y_d \wedge s^{1/\alpha})^{q-\beta_1+\eps}  	(y_d \vee  s^{1/\alpha})^{\eps}     \L^{\beta_3} \Big(\frac{1 }{ y_d \vee s^{1/\alpha}} \Big).
\end{align*}
Thus, the map $s \mapsto s \big( 1 \wedge \frac{ y_d}{s^{1/\alpha}}\big)^{q}     
(y_d \vee  s^{1/\alpha})^{\beta_1}  \L^{\beta_3}\big(\frac{
	1}{ y_d \vee s^{1/\alpha}} \big)$ is almost increasing on $(0,t]$
 by \eqref{e:slowly-varying-2}. Using this and  \eqref{e:first-kill-log-0}--\eqref{e:first-kill-log-1},   we get  that
\begin{align}\label{e:first-kill-log-2}
	\sup_{s \le t, \, z \in V_2} p(s,z,y) &\le c_3 \sup_{s \le t} s \Big( 1 \wedge \frac{ y_d}{s^{1/\alpha}}\Big)^{q}     
	(y_d \vee  s^{1/\alpha})^{\beta_1}  \L^{\beta_3} \Big(\frac{
	1}{ y_d \vee s^{1/\alpha}} \Big) \\
	&\le c_4t \Big( 1 \wedge \frac{ y_d}{t^{1/\alpha}}\Big)^{q}     
	(y_d \vee  t^{1/\alpha})^{\beta_1}  \L^{\beta_3} \Big(\frac{
	1}{ y_d \vee t^{1/\alpha}} \Big). \nn
\end{align}
Note that  \eqref{e:first-term-ass} is satisfied with  $a_1=\beta_1$ and 
$a_3=\beta_3$  by Lemma \ref{l:UHK-case1-induction}.  
Now combining Lemma \ref{l:first-term} and \eqref{e:first-kill-log-2},  we obtain the conclusion of the lemma.
\end{proof}

\begin{lemma}\label{l:UHK-case1-main1}
Let $\eta_1,\eta_2, \gamma \ge 0$.
There exists a constant $C>0$ such that for any $x \in \R^d_+$ and any $s,k,l>0$,
\begin{align*}
	& \int_{B_+(x, 2)} p(s, x, z) 
	z_d^{\gamma} \L^{\eta_1} \Big( \frac{k}{z_d} \Big) \L^{\eta_2} \Big( \frac{z_d}{l} \Big)  dz\\
		& \le Cx_d^{\gamma} \Big(1 \wedge \frac{x_d}{s^{1/\alpha}} \Big)^{q-\gamma}  \L^{\eta_1} \Big( \frac{k}{x_d \vee s^{1/\alpha}} \Big) \L^{\eta_2} \Big( \frac{x_d \vee s^{1/\alpha}}{l} \Big)\\
	&\quad + C\1_{\{\gamma>\alpha+\beta_1\}}  sx_d^{\beta_1} \Big(1 \wedge \frac{x_d}{s^{1/\alpha}} \Big)^{q-\beta_1}  \L^{\beta_3} \Big(\frac{  2}{x_d \vee s^{1/\alpha} } \Big) \L^{\eta_1} ( k ) \L^{\eta_2} \Big( \frac{1}{l} \Big) \\
	&\quad  + C\1_{\{\gamma=\alpha+\beta_1,\, x_d \vee s^{1/\alpha} <2 \}}  sx_d^{\beta_1} \Big(1 \wedge \frac{x_d}{s^{1/\alpha}} \Big)^{q-\beta_1} 
	\int_{x_d \vee s^{1/\alpha}}^2  \L^{\beta_3} \Big(\frac{  r}{x_d \vee s^{1/\alpha} } \Big) \L^{\eta_1}\Big(\frac{ k}{r} \Big) \L^{\eta_2}\Big(\frac{r}{l} \Big)  \frac{dr	}{r}.
\end{align*}
\end{lemma}
\begin{proof} For any $x,z \in \R^d_+$ and $s>0$, by Lemmas \ref{l:UHK-case1-induction} and \ref{l:kill-log-2}, 
\begin{align*}
&p(s,x,z) \le c_1 \Big(1 \wedge \frac{x_d}{s^{1/\alpha}}\Big)^q A_{\beta_1,0,\beta_3,0}(s,x,z) 
p_\alpha(s, x, z)\\
& \le  c_2   \Big(  1 \wedge \frac{x_d}{s^{1/\alpha}}\Big)^q   \Big(   1 \wedge \frac{x_d \vee s^{1/\alpha}}{|x-z|} \Big)^{\beta_1}            \L^{\beta_3}
 \Big(  \frac{|x-z|}{(x_d \vee s^{1/\alpha}) \wedge |x-z|} \Big)   
p_\alpha(s, x, z).
\end{align*}
Now combining \eqref{e:equiv-form} and  Lemma \ref{cal:2}, we get the desired result.
\end{proof}

We now state the first main result of this section. 

\begin{thm}\label{t:UHK}
For any $\eps \in (0,\alpha/2]$, there exists a constant $C>0$ such that
	\begin{align*}
		p(t,x,y) &\le C \Big(1 \wedge \frac{x_d}{t^{1/\alpha}}\Big)^{q} \Big(1 \wedge \frac{ y_d}{t^{1/\alpha}}\Big)^{q}
		A_{\beta_1, \beta_2 \wedge (\alpha+\beta_1-\eps), \beta_3, \beta_4}(t,x,y)
		p_\alpha(t, x, y),
	\end{align*}
	for all $t>0$ and $x,y \in \R^d_+$.
\end{thm}

This theorem will proved by using several lemmas. 
We first introduce some additional notation:
For $b_1,b_2,b_3,b_4 \ge 0$, $t>0$ and $x,y \in \R^d_+$, define
\begin{align}
	\sI_{t,1}(x,y;b_1,b_2,b_3,b_4)&:=\int_0^{t/2} 
\int_{B_+(y, 2)} \int_{B_+(x, 2)} \1_{\{u_d\le w_d\}} 
 p(t-s, x,u) 	p(s,y,w) 
  \label{e:sI-1} \\
 &\hspace{1cm} \times u_d^{b_1}w_d^{b_2} \L^{b_3} \Big( \frac{w_d}{u_d} \Big) \L^{b_4} \Big( \frac{1}{w_d} \Big)\, du\, dw\, ds,  \nn\\
 	\sI_{t,2}(x,y;b_1,b_2,b_3,b_4)&:=\int_0^{t/2} 
 \int_{B_+(y, 2)} \int_{B_+(x, 2)}
 \1_{\{w_d\le u_d\}} 
 p(t-s, x,u) 	p(s,y,w) 
 \label{e:sI-2}\\
 &\hspace{1cm} \times u_d^{b_2}w_d^{b_1} \L^{b_3} \Big( \frac{u_d}{w_d} \Big) \L^{b_4} \Big( \frac{1}{u_d} \Big)\, du\, dw\, ds. \nn
\end{align}

\begin{lemma}\label{l:newlemma0}
	Let $b_1,b_2,b_3,b_4 \ge 0$.
Let $b_1' \in [0,b_1]$ and $b_2':=b_1+b_2-b_1'$. Then  $$ \sI_{t,i}(x,y; b_1,b_2, b_3, b_4)\le \sI_{t,i}(x,y; b_1', b_2', b_3, b_4), \quad i=1,2.$$
\end{lemma}
\begin{proof}
 If $u_d\le w_d$, then $u_d^{b_1}w_d^{b_2}\le u_d^{b_1'}w_d^{b_2'}$, implying that 
 $\sI_{t,1}(x,y; b_1,b_2, b_3, b_4)\le \sI_{t,1}(x,y; b_1', b_2', b_3, b_4).$  The other inequality is analogous. \end{proof}

\begin{lemma}\label{l:newlemma1} 
	Let $b_1,b_2,b_3,b_4 \ge 0$ with $b_1 \vee b_2<\alpha+\beta_1$,  $x, y\in \R^d_+$ with  $|x-y|=5$, and $t\in (0,1]$. 
	There exists a constant $C>0$ independent of $t,x,y$ such that the following hold. 
	
	\smallskip

\noindent	(i)  If $b_2>q-\alpha$ or $y_d \ge t^{1/\alpha}$, then
	\begin{align*}
		&\sI_{t,1}(x,y; b_1,b_2,b_3,b_4)\\ 
		&\le C tx_d^{b_1}y_d^{b_2}\Big(1 \wedge \frac{x_d}{t^{1/\alpha}} \Big)^{q-b_1} \Big(1 \wedge \frac{y_d}{t^{1/\alpha}} \Big)^{q-b_2} 
		\L^{b_3} \Big( \frac{y_d \vee t^{1/\alpha}}{x_d \vee t^{1/\alpha}} \Big)\L^{b_4} \Big( \frac{1}{y_d \vee t^{1/\alpha}} \Big).
	\end{align*}
\noindent	(ii)  If $b_1>q-\alpha$ or $y_d \ge t^{1/\alpha}$, then
	\begin{align*}
                 &\sI_{t,2}(x,y; b_1,b_2,b_3,b_4)\\	
                 &\le C tx_d^{b_2}y_d^{b_1}\Big(1 \wedge \frac{x_d}{t^{1/\alpha}} \Big)^{q-b_2} \Big(1 \wedge \frac{y_d}{t^{1/\alpha}} \Big)^{q-b_1}  
		\L^{b_3} \Big( \frac{x_d \vee t^{1/\alpha}}{y_d \vee t^{1/\alpha}} \Big)\L^{b_4} \Big( \frac{1}{x_d \vee t^{1/\alpha}} \Big).
	\end{align*}
\end{lemma}
\begin{proof} We only give the proof of (i). (ii) can be proved similarly.

By using Lemma \ref{l:UHK-case1-main1} together with the fact that  
$t-s\asymp t$  if  $0\le s\le t/2$,  we see that  for $0\le s\le t/2$,  
\begin{align*}
 \int_{B_+(x, 2)} 
 p(t-s, x,u)   u_d^{b_1}\L^{b_3} \Big( \frac{w_d}{u_d} \Big) du
 	\le c_1 x_d^{b_1}\Big(1 \wedge \frac{x_d}{t^{1/\alpha}} \Big)^{q- b_1}   \L^{b_3} \Big( \frac{w_d}{x_d \vee t^{1/\alpha}} \Big).   
 \end{align*}
 Thus, using Lemma  \ref{l:UHK-case1-main1} again, we get that for $0\le s\le t/2$,  
 \begin{align*}
&\int_{B_+(y, 2)} \int_{B_+(x, 2)}
 p(t-s, x,u) 	p(s,y,w)  u_d^{b_1}w_d^{b_2} \L^{b_3} \Big( \frac{w_d}{u_d} \Big) \L^{b_4} \Big( \frac{1}{w_d} \Big)\, du\, dw\nn\\
 &\le c_1 x_d^{b_1}\Big(1 \wedge \frac{x_d}{t^{1/\alpha}} \Big)^{q- b_1}
	\int_{B_+(y, 2)} 
	p(s, y,w) w_d^{b_2}    \L^{b_3} \Big( \frac{w_d}{x_d \vee t^{1/\alpha}} \Big)     \L^{b_4} \Big( \frac{1}{ w_d} \Big)   dw \nn \\ 
	&\le c_2 x_d^{b_1}y_d^{b_2} \Big(1 \wedge \frac{x_d}{t^{1/\alpha}} \Big)^{q-b_1}   \Big(1 \wedge 
	\frac{y_d}{s^{1/\alpha}} \Big)^{q-b_2}   \L^{b_3} \Big( \frac{y_d \vee s^{1/\alpha}}{x_d \vee t^{1/\alpha}} \Big)   \L^{b_4} \Big( \frac{1}{y_d \vee s^{1/\alpha}} \Big)
 \end{align*}
From this and Lemma \ref{cal:0}, we get that
\begin{align*}
			&\sI_{t,1}(x,y; b_1,b_2,b_3,b_4)x_d^{-b_1}y_d^{-b_2} \Big(1 \wedge \frac{x_d}{t^{1/\alpha}} \Big)^{b_1-q} \nn\\
	&\le c_2    \int_0^{t/2} \Big(1 \wedge 
	\frac{y_d}{s^{1/\alpha}} \Big)^{q-b_2}   \L^{b_3} \Big( \frac{y_d \vee s^{1/\alpha}}{x_d \vee t^{1/\alpha}} \Big)   \L^{b_4} \Big( \frac{1}{y_d \vee s^{1/\alpha}} \Big) ds \nn \\
	&\le c_2 \int_0^{t} \Big(1 \wedge 
	\frac{y_d}{s^{1/\alpha}} \Big)^{q-b_2}   \L^{b_3} \Big( \frac{y_d \vee s^{1/\alpha}}{x_d \vee (t-s)^{1/\alpha}} \Big)   \L^{b_4} \Big( \frac{1}{y_d \vee s^{1/\alpha}} \Big) ds \nn \\
	&\le c_3 t \Big(1 \wedge \frac{y_d}{t^{1/\alpha}} 
	\Big)^{q-b_2}    \L^{b_3} \Big( \frac{y_d \vee t^{1/\alpha}}{x_d \vee t^{1/\alpha}} \Big)\L^{b_4} \Big( \frac{1}{y_d 
		\vee t^{1/\alpha}} \Big).
\end{align*}

\end{proof}

\begin{lemma}\label{l:newlemma1+} 
	Let $b_1,b_2,b_3,b_4 \ge 0$ be such that $b_1>q-\alpha$ and  $b_1 \vee b_2<\alpha+\beta_1$.  Assume that $b_2>0$ if $b_4>0$. Then, there exists a constant $C>0$ such that for all $x, y\in \R^d_+$ with $x_d \le y_d$ and   $|x-y|=5$, and all $t\in (0,1]$,
	\begin{align}\label{e:new-1}
		&\sI_{t,1}(x,y; b_1,b_2,b_3,b_4) + 	\sI_{t,2}(x,y; b_1,b_2,b_3,b_4)\\
		\le& C tx_d^{b_1}y_d^{b_2}\Big(1 \wedge \frac{x_d}{t^{1/\alpha}} \Big)^{q-b_1} \Big(1 \wedge \frac{y_d}{t^{1/\alpha}} \Big)^{q-b_2}   
		\L^{b_3} \Big( \frac{y_d \vee t^{1/\alpha}}{x_d \vee t^{1/\alpha}} \Big)\L^{b_4} \Big( \frac{1}{y_d \vee t^{1/\alpha}} \Big) \nn
	\end{align}
and
	\begin{align}\label{e:new-2}
              &\sI_{t,1}(y,x; b_1,b_2,b_3,b_4) +  \sI_{t,2}(y,x; b_1,b_2,b_3,b_4)	\\
              \le C& tx_d^{b_1}y_d^{b_2}\Big(1 \wedge \frac{x_d}{t^{1/\alpha}} \Big)^{q-b_1} \Big(1 \wedge \frac{y_d}{t^{1/\alpha}} \Big)^{q-b_2}  
	      \L^{b_3} \Big( \frac{y_d \vee t^{1/\alpha}}{x_d \vee t^{1/\alpha}} \Big)\L^{b_4} \Big( \frac{1}{y_d \vee t^{1/\alpha}} \Big). \nn
	\end{align}
\end{lemma}
\begin{proof}  Let $\delta \in (0,1)$ be such that \begin{align}\label{e:new-delta}
 	 b_1 \vee b_2 + (1-\delta)(b_1 \wedge b_2) <\alpha+\beta_1,
 \end{align}
and let $b_1':=\delta(b_1 \wedge b_2)$ and $b_2':=b_1+ b_2-b_1'$. 
We see that $b_1' \in [0,b_1 \wedge b_2]$, $b_2' \ge b_1>q-\alpha$ and
\begin{align}\label{e:new-delta'}
	b_1' \le b_2' = b_1 \vee b_2  + (1-\delta)(b_1 \wedge b_2)<\alpha+ \beta_1 
\end{align}
by \eqref{e:new-delta}. Moreover, since $x_d \le y_d$ and $b_2'-b_1= b_2-b_1'$,  we see that 
\begin{align}\label{e:newlemma1-0}
	(x_d \vee t^{1/\alpha})^{b_2'-b_1} (y_d  \vee t^{1/\alpha})^{b_1'-b_2} \L^{b_4} \Big(\frac{1}{x_d \vee t^{1/\alpha}} \Big) \Big/\L^{b_4} \Big(\frac{1}{y_d \vee t^{1/\alpha}} \Big) \le c_1.
\end{align}
Indeed, when $b_4=0$, \eqref{e:newlemma1-0} clearly holds with $c_1=1$. If $b_4>0$, then $b_2>0$ so that $b_2>\delta b_2 \ge b_1'$. Hence, we get \eqref{e:newlemma1-0} from \eqref{e:slowly-varying-2}.

\medskip

(i) We first prove \eqref{e:new-1}. For this, we distinguish between two cases: $y_d \ge t^{1/\alpha}$ and $y_d<t^{1/\alpha}$.

Assume first that $y_d \ge t^{1/\alpha}$.  The desired bound for $\sI_{t,1}(x,y;b_1,b_2,b_3,b_4)$ follows from Lemma \ref{l:newlemma1}(i). On the other hand, by using Lemma 
 \ref{l:newlemma0}  
in the first inequality, Lemma \ref{l:newlemma1}(ii) in the second inequality (which uses  \eqref{e:new-delta'} and  $y_d \ge t^{1/\alpha}$), \eqref{e:equiv-form} in the equality, and \eqref{e:newlemma1-0} in the last inequality, we get that
\begin{align*}
	&\sI_{t,2}(x,y; b_1,b_2,b_3,b_4)\le \sI_{t,2}(x,y; b_1',b_2',b_3,b_4)\nn \\
	&\le c_2 t x_d^{b_2'}y_d^{b_1'} \Big(1 \wedge \frac{x_d}{t^{1/\alpha}} \Big)^{q-b_2'}   \Big(1 \wedge \frac{y_d}{t^{1/\alpha}} \Big)^{q-b_1'}   \L^{b_3} \Big( \frac {y_d \vee t^{1/\alpha}}{x_d \vee t^{1/\alpha}}\Big)      \L^{b_4} \Big( \frac{1}{x_d \vee t^{1/\alpha}} \Big)\nn\\
	&=c_2 tx_d^{b_1}y_d^{b_2} \Big(1 \wedge \frac{x_d}{t^{1/\alpha}} \Big)^{q-b_1}   \Big(1 \wedge \frac{y_d}{t^{1/\alpha}} \Big)^{q-b_2}   \L^{b_3} \Big( \frac {y_d \vee t^{1/\alpha}}{x_d \vee t^{1/\alpha}}\Big)      \L^{b_4} \Big( \frac{1}{y_d \vee t^{1/\alpha}} \Big) \nn\\
	&\quad \;\times  (x_d \vee t^{1/\alpha})^{b_2'-b_1}  (y_d \vee t^{1/\alpha})^{b_1'-b_2}\L^{b_4} \Big( \frac{1}{x_d \vee t^{1/\alpha}} \Big)\Big/\L^{b_4} \Big( \frac{1}{y_d \vee t^{1/\alpha}} \Big)\nn\\
		&\le c_1c_2 tx_d^{b_1}y_d^{b_2} \Big(1 \wedge \frac{x_d}{t^{1/\alpha}} \Big)^{q-b_1}   \Big(1 \wedge \frac{y_d}{t^{1/\alpha}} \Big)^{q-b_2}   \L^{b_3} \Big( \frac {y_d \vee t^{1/\alpha}}{x_d \vee t^{1/\alpha}}\Big)      \L^{b_4} \Big( \frac{1}{y_d \vee t^{1/\alpha}} \Big).
\end{align*}

Assume now that $y_d <t^{1/\alpha}$. Using Lemma
\ref{l:newlemma0}  and Lemma \ref{l:newlemma1}(i) (which uses $b_2'>q-\alpha$), we get
\begin{align}\label{e:new-case11}
	& \sI_{t,1}(x,y; b_1,b_2,b_3,b_4)\le  \sI_{t,1}(x,y; b_1',b_2',b_3,b_4) \\
	\le c_3& tx_d^{b_1'}y_d^{b_2'} \Big(1 \wedge \frac{x_d}{t^{1/\alpha}} \Big)^{q-b_1'}  \Big(1 \wedge \frac{y_d}{t^{1/\alpha}} \Big)^{q-b_2'}   \L^{b_3} \Big( \frac{y_d \vee t^{1/\alpha}}{x_d \vee t^{1/\alpha}}    \Big) \L^{b_4} \Big(\frac{1}{y_d \vee t^{1/\alpha}}  \Big). \nn
\end{align}
Also, since $b_1>q-\alpha$, we get from Lemma \ref{l:newlemma1}(ii) that
\begin{align}\label{e:new-case12}
	 	& \sI_{t,2}(x,y; b_1,b_2,b_3,b_4)\\
	 \le c_4& tx_d^{b_2}y_d^{b_1} \Big(1 \wedge \frac{x_d}{t^{1/\alpha}} \Big)^{q-b_2}  \Big(1 \wedge \frac{y_d}{t^{1/\alpha}} \Big)^{q-b_1}   \L^{b_3} \Big( \frac{x_d \vee t^{1/\alpha}}{y_d \vee t^{1/\alpha}}    \Big) \L^{b_4} \Big( \frac{1}{x_d \vee t^{1/\alpha}}  \Big). \nn
\end{align}
Since $x_d \le y_d<t^{1/\alpha}$ and $b_1'+b_2'=b_1+b_2$, it holds that $x_d \vee t^{1/\alpha}=y_d \vee t^{1/\alpha}=t^{1/\alpha}$ and 
\begin{align*}
	&x_d^{b_1'}y_d^{b_2'} \Big(1 \wedge \frac{x_d}{t^{1/\alpha}} \Big)^{q-b_1'}  \Big(1 \wedge \frac{y_d}{t^{1/\alpha}} \Big)^{q-b_2'}= x_d^{b_2}y_d^{b_1} \Big(1 \wedge \frac{x_d}{t^{1/\alpha}} \Big)^{q-b_2}  \Big(1 \wedge \frac{y_d}{t^{1/\alpha}} \Big)^{q-b_1}\\
	&=  \frac{x_d^qy_d^q}{t^{(2q-b_1-b_2)/\alpha}}=	x_d^{b_1}y_d^{b_2} \Big(1 \wedge \frac{x_d}{t^{1/\alpha}} \Big)^{q-b_1}  \Big(1 \wedge \frac{y_d}{t^{1/\alpha}} \Big)^{q-b_2}.
\end{align*}
Thus,  \eqref{e:new-1} follows from \eqref{e:new-case11} and \eqref{e:new-case12}.

 (ii) Now, we prove \eqref{e:new-2}.
By using Lemma
\ref{l:newlemma0}   in the first inequality, Lemma \ref{l:newlemma1}(i) in the second inequality (which uses \eqref{e:new-delta'} and $b_2'>q-\alpha$), \eqref{e:equiv-form} together with the fact that $x_d \le y_d$ in the third inequality, \eqref{e:newlemma1-0} in the last inequality, we obtain
\begin{align*}
	& \sI_{t,1}(y,x; b_1,b_2,b_3,b_4) \le  \sI_{t,1}(y,x; b_1',b_2',b_3,b_4) \nn\\
	&\le c_5 tx_d^{b_2'}y_d^{b_1'} \Big(1 \wedge \frac{x_d}{t^{1/\alpha}} \Big)^{q-b_2'}  \Big(1 \wedge \frac{y_d}{t^{1/\alpha}} \Big)^{q-b_1'}   \L^{b_3} \Big( \frac{x_d \vee t^{1/\alpha}}{y_d \vee t^{1/\alpha}}    \Big) \L^{b_4} \Big( \frac{1}{x_d \vee t^{1/\alpha}}  \Big)\\
	&\le c_5 tx_d^{b_1}y_d^{b_2}\Big(1 \wedge \frac{x_d}{t^{1/\alpha}} \Big)^{q-b_1} \Big(1 \wedge \frac{y_d}{t^{1/\alpha}} 
	\Big)^{q-b_2}    \L^{b_3} \Big( \frac{y_d \vee t^{1/\alpha}}{x_d \vee t^{1/\alpha}} \Big)\L^{b_4} \Big( \frac{1}{y_d 
		\vee t^{1/\alpha}} \Big)\nn\\
	&\quad \times (x_d \vee t^{1/\alpha})^{b_2'-b_1} (y_d \vee t^{1/\alpha})^{b_1'-b_2} \L^{b_4} \Big( \frac{1}{x_d \vee t^{1/\alpha}}  \Big)\Big/\L^{b_4} \Big( \frac{1}{y_d \vee t^{1/\alpha}}  \Big)\nn\\
	&\le c_1c_5 tx_d^{b_1}y_d^{b_2}\Big(1 \wedge \frac{x_d}{t^{1/\alpha}} \Big)^{q-b_1} \Big(1 \wedge \frac{y_d}{t^{1/\alpha}} 
	\Big)^{q-b_2}    \L^{b_3} \Big( \frac{y_d \vee t^{1/\alpha}}{x_d \vee t^{1/\alpha}} \Big)\L^{b_4} \Big( \frac{1}{y_d 
		\vee t^{1/\alpha}} \Big).
\end{align*}
On the other hand, by Lemma \ref{l:newlemma1}(ii), it holds that 
\begin{align*}
	& \sI_{t,2}(y,x; b_1,b_2,b_3,b_4)  \nn\\
	&\le c_3 tx_d^{b_1}y_d^{b_2}\Big(1 \wedge \frac{x_d}{t^{1/\alpha}} \Big)^{q-b_1} \Big(1 \wedge \frac{y_d}{t^{1/\alpha}} 
	\Big)^{q-b_2}    \L^{b_3} \Big( \frac{y_d \vee t^{1/\alpha}}{x_d \vee t^{1/\alpha}} \Big)\L^{b_4} \Big( \frac{1}{y_d 
		\vee t^{1/\alpha}} \Big).
\end{align*}
The proof is complete. \end{proof}

\noindent
\textsc{Proof of Theorem \ref{t:UHK}.} 
Set $\wh\beta_2:=\beta_2 \wedge (\alpha+\beta_1-\eps)$. As in the proof of Lemma \ref{l:UHK-case1-induction}, by symmetry, Proposition \ref{p:UHK-rough}, \eqref{e:kernel-scaling}, 
and \eqref{e:wtB-interior},
we can assume without loss of generality that $x_d \le y_d \wedge 2^{-5}$, $\widetilde{x}=0$ and $|x-y|=5$, and  then  it is enough to show that there exists a constant $c_1>0$ independent of $x$ and $y$  such that for any $t \le 1$,
\begin{align}\label{e:UHK1-claim}
	p(t,x,y) &\le c_1 t\left(1 \wedge \frac{x_d }{t^{1/\alpha}}\right)^{q} \left(1 \wedge \frac{y_d }{t^{1/\alpha}}\right)^{q}  (x_d \vee t^{1/\alpha})^{\beta_1}  (y_d \vee t^{1/\alpha})^{\wh\beta_2} \\
	& \quad \times \L^{\beta_3} \Big(\frac{ y_d \vee t^{1/\alpha} }{ x_d \vee t^{1/\alpha}} \Big) \,\L^{\beta_4}\Big(\frac{1}{ y_d \vee t^{1/\alpha} } \Big). \nn
\end{align}

Let $t \le 1$.  Set $V_1=U(1)$, $V_3=B(y,2) \cap \oR^d_+$ and $V_2=\oR^d_+ \setminus (V_1 \cup V_3)$.
By Lemma \ref{l:firstpart},
\begin{align}\label{e:firstpart}
&	\P_x(\tau_{V_1}<t< \zeta) \sup_{s \le t, \, z \in V_2} p(s,z,y) \\
&\le c_2t^2 \left(1 \wedge \frac{x_d}{t^{1/\alpha}}\right)^{q} \left( 1 \wedge \frac{ y_d}{t^{1/\alpha}}\right)^{q} (x_d \vee t^{1/\alpha})^{\beta_1}   (y_d \vee t^{1/\alpha})^{\beta_1} \L^{2\beta_3}\Big(\frac{1}{ t^{1/\alpha}} \Big)\nn\\
&\le c_2t \left(1 \wedge \frac{x_d}{t^{1/\alpha}}\right)^{q} \left( 1 \wedge \frac{ y_d}{t^{1/\alpha}}\right)^{q} (x_d \vee t^{1/\alpha})^{\beta_1}   (y_d \vee t^{1/\alpha})^{\wh\beta_2}   t^{\eps/\alpha}
\L^{2\beta_3}\Big(\frac{1}{ t^{1/\alpha}} \Big)\nn\\
&\le c_3t \left(1 \wedge \frac{x_d}{t^{1/\alpha}}\right)^{q} \left( 1 \wedge \frac{ y_d}{t^{1/\alpha}}\right)^{q} (x_d \vee t^{1/\alpha})^{\beta_1}   (y_d \vee t^{1/\alpha})^{\wh\beta_2}, \nn
\end{align}
where in the last inequality above we used \eqref{e:slowly-varying}.

Next, we show that there exists a constant $C'>0$ such that 
\begin{align}\label{e:secondpart}
	&I:= \int_0^t \int_{V_3}\int_{V_1} p^{V_1}(t-s,x,u) \sB(u,w) p(s,y,w)du\,dw\,ds \\
	\le C'	 &t  x_d^{\beta_1}  y_d^{\wh\beta_2} \left(1 \wedge \frac{x_d}{t^{1/\alpha}} \right)^{q-\beta_1}    \left(1 \wedge \frac{y_d}{t^{1/\alpha}} \right)^{q- \wh\beta_2}  \L^{\beta_3} \Big( \frac{ y_d \vee t^{1/\alpha}}{ x_d \vee t^{1/\alpha} } \Big)   \L^{\beta_4} \Big( \frac{1}{ y_d \vee t^{1/\alpha} } \Big). \nn
\end{align} 
Once we get \eqref{e:secondpart}, by \eqref{e:equiv-form} and  \eqref{e:firstpart},  we can apply Lemma \ref{l:general-upper-2}  to get \eqref{e:UHK1-claim} and finish the proof.  

By {\bf (A3)}(II), since  $|u-w|\asymp 1$
and $u_d\vee w_d \le y_d+2 \le  7 +2^{-5}$ for $u \in V_1$ and $w \in V_3$,  we have
\begin{align*}
	I	&\le c_{4} \int_0^{t} \int_{V_3} \int_{V_1}\1_{\{u_d \le w_d\}} p(t-s,x,u) p(s, y,w)  u_d^{\beta_1} w_d^{\wh\beta_2}  \L^{\beta_3} \Big( \frac{w_d}{u_d} \Big) \L^{\beta_4} \Big( \frac{1}{ w_d} \Big)     du    dw ds\\
	&\; +  c_{4} \int_0^{t} \int_{V_3} \int_{V_1} \1_{\{w_d \le u_d\}} p(t-s,x,u) p(s, y,w)   
u_d^{\wh \beta_2} w_d^{\beta_1}  	\L^{\beta_3} \Big( \frac{u_d}{w_d} \Big) \L^{\beta_4} \Big( \frac{1}{ u_d} \Big)     du   dw ds\\
 	& = c_4 \Big(  \big(\int_0^{t/2}+ 
	\int_{t/2}^{t} 
	\big) \int_{V_3} \int_{V_1} \1_{\{u_d \le w_d\}} \dots +
	\big(\int_0^{t/2}+ 
	\int_{t/2}^{t} \big) \int_{V_3} \int_{V_1} \1_{\{w_d \le u_d\}} \dots\Big).
\end{align*}
Thus, by the change of variable $\wt s = t-s$ in integrals $\int_{t/2}^t$,
\begin{align*}
	I	 \le c_4\sum_{i=1}^2\sI_{t,i}(x,y; \beta_1, \wh{\beta}_2, \beta_3, \beta_4)+c_4\sum_{i=1}^2\sI_{t,i}(y,x; \beta_1, \wh{\beta}_2, \beta_3, \beta_4),\end{align*}
where the functions $\sI_{t,i}(x,y; \beta_1, \wh{\beta}_2, \beta_3, \beta_4)$,  $1\le i \le 2$, are defined in \eqref{e:sI-1}--\eqref{e:sI-2}. 
Then by using 
Lemma \ref{l:newlemma1+}, we conclude that \eqref{e:secondpart} holds.  The proof is complete. 
\qed 

Combining Remark \ref{r:d=1} and Theorem \ref{t:UHK}, we immediately get the
sharp heat kernel upper bound for $d=1$ or
  $\beta_2<\alpha+\beta_1$. 
\begin{cor}\label{c:UHK}
 	Suppose that 
 $d=1$ or 
 	 $\beta_2<\alpha+\beta_1$.  There exists a constant $C>0$ such that for all $t>0$ and $x,y \in \R^d_+$,
 	\begin{align*}
 		p(t,x,y) &\le C \left(1 \wedge \frac{x_d}{t^{1/\alpha}}\right)^{q} \left(1 \wedge \frac{ y_d}{t^{1/\alpha}}\right)^{q}
 		A_{\beta_1, \beta_2, \beta_3, \beta_4}(t,x,y)
		p_\alpha(t, x, y).
 	\end{align*}
 	
\end{cor}

\smallskip

Here is the second main result of the section.
\begin{thm}\label{t:UHK2}
Suppose that $d \ge 2$ and
$\beta_2 \ge \alpha+\beta_1$.  There exists a constant $C>0$ such that for all $t>0$ and $x,y \in \R^d_+$,
		\begin{align*}
 		 p(t,x,y) &\le C \left(1 \wedge \frac{x_d}{t^{1/\alpha}}\right)^{q} \left(1 \wedge \frac{ y_d}{t^{1/\alpha}}\right)^{q}
		p_\alpha(t, x, y)
 		\Big[	A_{\beta_1, \beta_2, \beta_3, \beta_4}(t,x,y)	\nn\\
 		&\quad + \left(1 \wedge  \frac{t}{|x-y|^\alpha}\right)\L^{\beta_3} \Big( \frac{|x-y|}{((x_d \wedge y_d)+ t^{1/\alpha}) \wedge |x-y|} \Big)\\
 		&\quad\quad \times \Big( \1_{\{\beta_2>\alpha+\beta_1\}}A_{\beta_1, \beta_1, 0, \beta_3}(t,x,y)+ \1_{\{\beta_2=\alpha+\beta_1\}}A_{\beta_1, \beta_1, 0, \beta_3+\beta_4+1}(t,x,y)\Big) \Big].
 	\end{align*}

\end{thm}

Again, we first introduce additional notation and prove a lemma.

 For $b_1,b_2,b_3,b_4 \ge 0$, $t>0$ and $x,y \in \R^d_+$, define
\begin{align}\label{e:sI}
 	&\sI_t(x,y; b_1,b_2,b_3,b_4)\\:=&
	\int_0^t 
	\int_{B_+(y, 2)} \int_{B_+(x, 2)}
 p(t-s, x,u) 	p(s,y,w)u_d^{b_1}w_d^{b_2} \L^{b_3} \Big( \frac{w_d}{u_d} \Big) \L^{b_4} \Big( \frac{1}{w_d} \Big)\, du\, dw\, ds. \nn
\end{align}

\begin{lemma}\label{l:newlemma2} 
Let $b_1,b_2,b_3,b_4\ge 0$ be such that $b_1 \wedge b_2>q-\alpha$ and $b_1<\alpha+\beta_1$. There exists a constant $C>0$ such that for all $x, y\in \R^d_+$ with $|x-y|=5$, and all $t\in (0,1]$,
\begin{align*}
& \sI_t(x,y; b_1,b_2,b_3,b_4) 
\nn\\
&\le C t x_d^{b_1}y_d^{b_2} \Big(1 \wedge \frac{x_d}{t^{1/\alpha}} \Big)^{q-b_1}   \Big(1 \wedge \frac{y_d}{t^{1/\alpha}} \Big)^{q-b_2}   \L^{b_3} \Big( \frac{y_d \vee t^{1/\alpha}}{x_d \vee t^{1/\alpha}} \Big)   \L^{b_4} \Big(  + \frac{1}{y_d \vee t^{1/\alpha}} \Big)\\
 	& + \1_{\{b_2>\alpha+\beta_1\}} C t^2 x_d^{b_1}y_d^{\beta_1}\Big(1 \wedge \frac{x_d}{t^{1/\alpha}} \Big)^{q-b_1}  \Big(1 \wedge \frac{y_d}{t^{1/\alpha}} \Big)^{q-\beta_1}   \L^{b_3} \Big(  \frac{1}{x_d \vee t^{1/\alpha}} \Big)     \L^{b_3} \Big(  \frac{1}{y_d \vee t^{1/\alpha}} \Big)\\
 	& +\1_{\{b_2=\alpha+\beta_1, y_d<2 \}} C t x_d^{b_1}y_d^{ \beta_1 }  \Big(1 \wedge \frac{x_d}{t^{1/\alpha}} \Big)^{q-b_1}  \Big(1 \wedge \frac{y_d}{t^{1/\alpha}} \Big)^{q- \beta_1} \\
	&\quad\times  \int_{y_d}^{2}(r^\alpha \wedge t)\L^{b_3} \Big( \frac{r}{x_d \vee t^{1/\alpha}} \Big)   \L^{b_3} \Big( \frac{r}{y_d \vee t^{1/\alpha}} \Big)\L^{b_4} \Big( \frac{1}{r}\Big)    \frac{dr}{r}.
\end{align*}
\end{lemma}
\begin{proof} 
 By using Lemma \ref{l:UHK-case1-main1} in the first inequality (integration with respect to $u$; note that 
$b_1<\alpha+\beta_1$) and the second inequality (integration with respect to $w$),
we get
\begin{align*}
 	&\sI_t(x,y; b_1,b_2,b_3,b_4)\le c_1 x_d^{b_1}\int_0^{t} \Big(1 \wedge \frac{x_d}{(t-s)^{1/\alpha}} \Big)^{q-b_1}\\
 	&\hspace{2.5cm}\times 
         \int_{B_+(y, 2)} 
         p(s, y,w) 
	w_d^{b_2}    \L^{b_3} \Big( \frac{w_d}{x_d \vee (t-s)^{1/\alpha}} \Big)     \L^{b_4} \Big( \frac{1}{ w_d} \Big)   dw ds \\ 
 	&\le c_2\left (x_d^{b_1}y_d^{b_2}I_1+
 	\1_{\{b_2>\alpha+\beta_1\}}  x_d^{b_1}y_d^{\beta_1}I_2+ 
	\1_{\{b_2=\alpha+\beta_1, y_d<2 \}} x_d^{b_1}y_d^{\beta_1} I_3\right), 
\end{align*} 
	where
\begin{align*}
I_1:=  \int_0^{t} \Big(1 \wedge \frac{x_d}{(t-s)^{1/\alpha}} \Big)^{q-b_1}  \Big(1 \wedge \frac{y_d}{s^{1/\alpha}} \Big)^{q-b_2}  \L^{b_3} \Big(\frac{y_d \vee s^{1/\alpha}}{x_d \vee (t-s)^{1/\alpha}} \Big)  \L^{b_4} \Big(\frac{1}{y_d \vee s^{1/\alpha}} \Big) ds,
\end{align*}
\begin{align*}
I_2:= \int_0^{t} s\Big(1 \wedge \frac{x_d}{(t-s)^{1/\alpha}} \Big)^{q-b_1}  \Big(1 \wedge 
 	\frac{y_d}{s^{1/\alpha}} \Big)^{q-\beta_1 }   
 	\L^{b_3} \Big( \frac{2}{y_d \vee s^{1/\alpha}} \Big)\ \L^{b_3} \Big( \frac{1}{x_d \vee (t-s)^{1/\alpha}} \Big)   ds,
\end{align*}
	and
\begin{align*}	
I_3:=&\int_0^{t}  s\Big(1 \wedge \frac{x_d}{(t-s)^{1/\alpha}} \Big)^{q-  b_1}  \Big(1 \wedge \frac{y_d}{s^{1/\alpha}} \Big)^{q-\beta_1}   \\
	&   \times   \int_{y_d \vee s^{1/\alpha}}^2  \L^{b_3} \Big( \frac{r}{y_d \vee s^{1/\alpha}} \Big)\L^{b_3} \Big( \frac{r}{x_d \vee (t-s)^{1/\alpha}} \Big) \L^{b_4} \Big( \frac{1}{r}\Big) \frac{dr}{r} \, ds. 
 \end{align*} 
Applying
 Lemma \ref{cal:0} to $I_1$ and  $I_2$, 
 Lemma \ref{l:cal1} to $I_3$, 
 and \eqref{e:slowly-varying-2}, wee see that
$$
 I_1 \le 
  c_3 t  \Big(1 \wedge \frac{x_d}{t^{1/\alpha}} \Big)^{q-b_1}  \Big(1 \wedge \frac{y_d}{t^{1/\alpha}} \Big)^{q-b_2}   \L^{b_3} \Big( \frac{y_d \vee t^{1/\alpha}}{x_d \vee t^{1/\alpha}} \Big)   \L^{b_4} \Big( \frac{1}{y_d \vee t^{1/\alpha}} \Big),
 $$
 $$	
	I_2  \le 
	c_3 t^2 \Big(1 \wedge \frac{x_d}{t^{1/\alpha}} \Big)^{q-b_1} \Big(1 \wedge \frac{y_d}{t^{1/\alpha}} \Big)^{q-\beta_1}    \L^{b_3} \Big( \frac{1}{x_d \vee t^{1/\alpha}} \Big)     \L^{b_3} \Big( \frac{1}{y_d \vee t^{1/\alpha}} \Big),
	$$
 	and
	 \begin{align*} 
I_3  \le &
	 c_3 t   \Big(1 \wedge \frac{x_d}{t^{1/\alpha}} \Big)^{q-b_1} \Big(1 \wedge \frac{y_d}{t^{1/\alpha}} \Big)^{q- \beta_1}   \int_{y_d}^{2} (r^\alpha \wedge t)\L^{b_3} \Big( \frac{r}{x_d \vee t^{1/\alpha}} \Big)   \L^{b_3} \Big( \frac{r}{y_d \vee t^{1/\alpha}} \Big)\L^{b_4} \Big( \frac{1}{r}\Big)    \frac{dr}{r}.
 \end{align*}  
This proves the lemma. \end{proof}

\noindent
\textsc{Proof of Theorem \ref{t:UHK2}.}   As in the proof of Lemma \ref{l:UHK-case1-induction}, by symmetry, \eqref{e:two-jump-compare-3},
  Proposition \ref{p:UHK-rough}, \eqref{e:kernel-scaling} and \eqref{e:wtB-interior}, 
 we can assume without loss of generality that $x_d \le y_d \wedge 2^{-5}$, $\widetilde{x}=0$ and  $|x-y|=5$, 
 and then  it is enough to show that there exists a constant $c_1>0$ independent of $x$ and $y$  such that for any $t \le 1$,
  \begin{align} 
  	p(t,x,y) &\le  c_1 t\left(1 \wedge \frac{x_d }{t^{1/\alpha}}\right)^{q} \left(1 \wedge \frac{y_d }{t^{1/\alpha}}\right)^{q} (x_d \vee t^{1/\alpha})^{\beta_1}  
  \label{e:UHK2-claim0}\\
		&\quad\times \Big[ (y_d \vee t^{1/\alpha})^{\beta_2}
  \L^{\beta_3} \Big(\frac{ y_d \vee t^{1/\alpha} }{ x_d \vee t^{1/\alpha}} \Big) \,\L^{\beta_4}\Big(\frac{1}{ y_d \vee t^{1/\alpha} } \Big) \nn\\
 	 &
	\qquad 
	 +  \1_{\{\beta_2>\alpha+\beta_1\}} t  (y_d \vee t^{1/\alpha})^{\beta_1} \L^{\beta_3} \Big(\frac{1}{ x_d \vee t^{1/\alpha}} \Big)\L^{\beta_3} \Big(\frac{1}{ y_d \vee t^{1/\alpha}} \Big)\nn \\
 	 &
	\qquad 
	 +  \1_{\{\beta_2=\alpha+\beta_1\}}t   (y_d \vee t^{1/\alpha})^{\beta_1} \L^{\beta_3}\Big(\frac{1}{ x_d \vee t^{1/\alpha}}\Big)\L^{\beta_3+\beta_4+1}\Big(\frac{1}{ y_d \vee t^{1/\alpha}} \Big) \Big]. 	  \nn
 \end{align}

 Let $t \le 1$.  Set $V_1=U(1)$, $V_3=B(y,2) \cap \oR^d_+$ and $V_2=\oR^d_+ \setminus (V_1 \cup V_3)$.
 By  Lemmas \ref{l:general-upper-2} and  \ref{l:firstpart}  
 it remains to prove that $I:= \int_0^t\int_{V_3} \int_{V_1} p^{V_1}(t-s, x, u) \sB(u,w) p(s, y,w) du dw ds$ is bounded above by 
  the last term on the right-hand side of \eqref{e:UHK2-claim0}.
  
By {\bf (A3)}(II), since  $|u-w|\asymp 1$ for $u \in V_1$ and $w \in V_3$, using the change of variables $\wt s = t-s$ we have 
\begin{align}\label{e:two-jump-decomp0}
	 I\le c_2 \left(\sI_t(x,y;\beta_1,\beta_2,\beta_3,\beta_4) + \sI_t(y,x;\beta_1,\beta_2,\beta_3,\beta_4)\right),
\end{align}
where the functions $\sI_t(x,y;\beta_1,\beta_2,\beta_3,\beta_4)$ is defined in \eqref{e:sI}. 
 By Lemma \ref{l:newlemma2} and \eqref{e:equiv-form}, the right hand side of \eqref{e:two-jump-decomp0} is less than or equal to 
 $ c_3 t (1 \wedge \frac{x_d }{t^{1/\alpha}})^{q} (1 \wedge \frac{y_d }{t^{1/\alpha}})^{q}$ times   
\begin{align}
& (x_d \vee t^{1/\alpha})^{\beta_1}  (y_d \vee t^{1/\alpha})^{\beta_2}
 \L^{\beta_3} \Big(\frac{ y_d \vee t^{1/\alpha} }{ x_d \vee t^{1/\alpha}} \Big) \,\L^{\beta_4}\Big(\frac{1}{ y_d \vee t^{1/\alpha} } \Big)  \nn\\
& \quad +   (x_d \vee t^{1/\alpha})^{\beta_2}  (y_d \vee t^{1/\alpha})^{\beta_1}  \L^{\beta_3} \Big( \frac{x_d \vee t^{1/\alpha}}{y_d \vee t^{1/\alpha}}\Big)  
\L^{\beta_4} \Big( \frac{1}{x_d \vee t^{1/\alpha}} \Big)  \nn\\
&\quad +  \1_{\{\beta_2>\alpha+\beta_1\}}  t
(x_d \vee t^{1/\alpha})^{\beta_1}   (y_d \vee t^{1/\alpha})^{\beta_1}  \L^{\beta_3} \Big(\frac{1}{ x_d \vee t^{1/\alpha}} \Big)\L^{\beta_3} \Big(\frac{1}{ y_d \vee t^{1/\alpha}} \Big) \nn\\
&\quad + \1_{\{\beta_2=\alpha+\beta_1\}} 
(x_d \vee t^{1/\alpha})^{\beta_1}   (y_d \vee t^{1/\alpha})^{\beta_1}\nn\\
&\qquad \;  \times\int_{x_d}^{2}(r^\alpha \wedge t)\L^{\beta_3} \Big( \frac{r}{x_d \vee t^{1/\alpha}} \Big)   \L^{\beta_3} \Big( \frac{r}{y_d \vee t^{1/\alpha}} \Big)\L^{\beta_4} \Big( \frac{1}{r}\Big)    \frac{dr}{r}. \nn\\ &=:I+II+III+IV.\nn
\end{align} 
For $IV$ we used $x_d\le y_d \wedge 2$ to get that $ \1_{\{y_d<2\}} \int_{y_d}^2\le \int_{x_d}^2$.

Since $\beta_2>\beta_1$ and $x_d \le y_d$,   by choosing $\epsilon \in (0,\beta_2-\beta_1)$ we see that
\begin{align}\label{e:two-jump-decomp'}
&(x_d \vee t^{1/\alpha})^{\beta_2}  (y_d \vee t^{1/\alpha})^{\beta_1} \L^{\beta_3} 
		\bigg(  \frac{x_d \vee t^{1/\alpha}}{y_d \vee t^{1/\alpha}}\bigg)
         \L^{\beta_4} \bigg( \frac{1}{x_d \vee t^{1/\alpha}} \bigg)\\
         &=(x_d \vee t^{1/\alpha})^{\beta_2-\beta_1-\epsilon} (x_d \vee t^{1/\alpha})^{\beta_1} (y_d \vee t^{1/\alpha})^{\beta_1}\nn \\
         &\quad \times \L^{\beta_3} \bigg(  \frac{x_d \vee t^{1/\alpha}}{y_d \vee t^{1/\alpha}}\bigg)(x_d \vee t^{1/\alpha})^{\epsilon}\L^{\beta_4} \bigg(  \frac{1}{x_d \vee t^{1/\alpha}} \bigg) \nn \\
         &\le (y_d \vee t^{1/\alpha})^{\beta_2-\beta_1-\epsilon} (x_d \vee t^{1/\alpha})^{\beta_1} (y_d \vee t^{1/\alpha})^{\beta_1}\nn \\
         &\quad \times \L^{\beta_3} \bigg( \frac{y_d \vee t^{1/\alpha}}{x_d \vee t^{1/\alpha}}\bigg)(x_d \vee t^{1/\alpha})^{\epsilon}\L^{\beta_4} \bigg( \frac{1}{x_d \vee t^{1/\alpha}} \bigg) \nn \\
         & \le c_4 (x_d \vee t^{1/\alpha})^{\beta_1}  (y_d \vee t^{1/\alpha})^{\beta_2} 
		\L^{\beta_3} \bigg(  \frac{y_d \vee t^{1/\alpha}}{x_d \vee t^{1/\alpha}}\bigg)         
         \L^{\beta_4} \bigg(  \frac{1}{y_d \vee t^{1/\alpha}} \bigg).\nn
\end{align}
The last inequality follows from almost increase of the function $u\mapsto u^{\epsilon}\L^{\beta_4}(1/u)$. Applying  \eqref{e:two-jump-decomp'}  to $II$
and  combining it with  $I$ and $III$,
we arrive at the result in case $\beta_2>\alpha+\beta_1$.

Assume now that $\beta_2=\alpha+\beta_1$.  From the above argument, 
to prove the result,  
in view of \eqref{e:UHK2-claim0} 
and $IV$, 
it suffices to show that 
 \begin{align}\label{e:estimate-of-I'}
& \int_{x_d}^{2}(r^\alpha \wedge t)\L^{\beta_3} \Big( \frac{r}{x_d \vee t^{1/\alpha}} \Big)   \L^{\beta_3} \Big( \frac{r}{y_d \vee t^{1/\alpha}} \Big)\L^{\beta_4} \Big( \frac{1}{r}\Big)    \frac{dr}{r} \\
 &\le c_5 t   \L^{\beta_3} \Big(\frac{1}{ x_d \vee t^{1/\alpha}}\Big)
 	  \L^{\beta_3+\beta_4+1}\Big(
	  \frac{1}{ y_d \vee t^{1/\alpha}} \Big)\nn\\
	  &\quad+c_5 (y_d \vee t^{1/\alpha})^{ \alpha}\L^{\beta_3} \Big( \frac{y_d \vee t^{1/\alpha}}{x_d \vee t^{1/\alpha}} \Big)\L^{\beta_4} \Big( \frac{1}{y_d \vee t^{1/\alpha}}\Big). \nn
 \end{align}

By Lemma \ref{cal:3} (with $b_1=\beta_3$, $b_2=\beta_4$, $k=x_d \vee t^{1/\alpha}$ and $l=y_d \vee t^{1/\alpha}$), it holds that 
\begin{align}\label{e:crit1}
	&\int_{y_d \vee t^{1/\alpha}}^{2}(r^\alpha \wedge t)\L^{\beta_3} \Big( \frac{r}{x_d \vee t^{1/\alpha}} \Big)   \L^{\beta_3} \Big( \frac{r}{y_d \vee t^{1/\alpha}} \Big)\L^{\beta_4} \Big( \frac{1}{r}\Big)    \frac{dr}{r} \\
	&\le t\int_{y_d \vee t^{1/\alpha}}^{2}\L^{\beta_3} \Big( \frac{r}{x_d \vee t^{1/\alpha}} \Big)   \L^{\beta_3} \Big( \frac{r}{y_d \vee t^{1/\alpha}} \Big)\L^{\beta_4} \Big( \frac{1}{r}\Big)    \frac{dr}{r}\nn\\
	&\le c_6t\L^{\beta_3} \Big( \frac{1}{x_d \vee t^{1/\alpha}}\Big) \L^{\beta_3+\beta_4+1} \Big(\frac{1}{y_d \vee t^{1/\alpha}}\Big). \nn
\end{align}
On the other hand, using \eqref{e:slowly-varying-2},  we see that
\begin{align}\label{e:crit2}
	&\int_{x_d}^{y_d \vee t^{1/\alpha}} (r^\alpha \wedge t)\L^{\beta_3} \Big( \frac{r}{x_d \vee t^{1/\alpha}} \Big)   \L^{\beta_3} \Big( \frac{r}{y_d \vee t^{1/\alpha}} \Big)\L^{\beta_4} \Big( \frac{1}{r}\Big)    \frac{dr}{r} \\
		&\le \L^{\beta_3}(1)\L^{\beta_3} \Big( \frac{y_d \vee t^{1/\alpha}}{x_d \vee t^{1/\alpha}} \Big)  \int_{x_d}^{y_d \vee t^{1/\alpha}} \L^{\beta_4} \Big( \frac{1}{r}\Big)    \frac{dr}{r^{1-\alpha}}\nn\\
	&\le c_7\L^{\beta_3} \Big( \frac{y_d \vee t^{1/\alpha}}{x_d \vee t^{1/\alpha}} \Big)\L^{\beta_4} \Big( \frac{1}{y_d \vee t^{1/\alpha}}\Big)\int_{x_d}^{y_d \vee t^{1/\alpha}} \Big( \frac{y_d \vee t^{1/\alpha}}{r} \Big)^{\alpha/2} \frac{dr}{r^{1-\alpha}}\nn\\
	&\le c_8(y_d \vee t^{1/\alpha})^\alpha\L^{\beta_3} \Big( \frac{y_d \vee t^{1/\alpha}}{x_d \vee t^{1/\alpha}} \Big)\L^{\beta_4} \Big( \frac{1}{y_d \vee t^{1/\alpha}}\Big).\nn
\end{align}
Combining \eqref{e:crit1}--\eqref{e:crit2}, we show that \eqref{e:estimate-of-I'} holds true. The proof is complete. 
\qed

\section{Proofs of Theorems \ref{t:HKE} and  \ref{t:HKE-kappa}}  \label{s:pmain}

\noindent \textsc{Proof of Theorem \ref{t:HKE}.} Since $\bar p(t,x,y)$ is jointly continuous 
(see  Remark \ref{r:conti}), 
it suffices to prove that 
\eqref{e:allcase1}--\eqref{e:Case3} hold for $(t,x,y) \in (0,\infty) \times \R^d_+ \times \R^d_+$.

  Using \eqref{e:comp-AB} and {\bf (A3)}, we get the lower heat kernel estimate in \eqref{e:allcase1} from  Proposition \ref{p:HKE-lower-1} and the upper heat kernel estimate from   Corollary \ref{c:UHK}.

 For remainder of the proof, we assume that $d \ge 2$. 
We first note 
that by  {\bf (A3)} and \eqref{e:wtB-interior}, 
\begin{align}\label{e:Case111}
& t^{-d/\alpha} \wedge
		\big(tJ(x+t^{1/\alpha}\e_d,y+t^{1/\alpha}\e_d )\big)
		 \asymp t^{-d/\alpha} \wedge \frac{t 	B_{\beta_1, \beta_2, \beta_3, \beta_4}(x+t^{1/\alpha}\e_d,y+t^{1/\alpha}\e_d)	}{|x-y|^{d+\alpha}} \\
		& \asymp  
		  \Big( t^{-d/\alpha} \wedge \frac{t}{|x-y|^{d+\alpha}}\Big)	B_{\beta_1, \beta_2, \beta_3, \beta_4}(x+t^{1/\alpha}\e_d,y+t^{1/\alpha}\e_d), \nn
	\end{align}
	which implies the second comparison in \eqref{e:Case1}.

(i) 
 Using \eqref{e:comp-AB}, we get the lower heat kernel estimate in the first comparison  in \eqref{e:Case1} from  Proposition \ref{p:HKE-lower-1} and the upper heat kernel estimate from   Corollary \ref{c:UHK}.

(ii) For  \eqref{e:Case2},  using \eqref{e:comp-AB}, we get the lower heat kernel estimate from  Proposition \ref{p:lower-1} (see Remark \ref{r:lowerbound}) and the upper heat kernel estimate from  Theorem \ref{t:UHK2}.  

(iii)  Using \eqref{e:comp-AB}, we get the lower heat kernel estimate in \eqref{e:Case3} from  Proposition \ref{p:lower-2}.
The upper heat kernel estimate in \eqref{e:Case3} follows  from  Theorem \ref{t:UHK2}. 

From the comparisons in (i)-(iii) and 
\eqref{e:wtB-interior}, we have that $\bar p(t,x,y) \asymp 	 	 t^{-d/\alpha}$ when $t^{1/\alpha} \ge  |x-y|/8$ and 
this implies that \eqref{e:allcase} holds for $t^{1/\alpha} \ge  |x-y|/8$.
 Moreover,  \eqref{e:allcase} for  $\beta_2<\alpha+\beta_1$ follows from \eqref{e:Case1}, \eqref{e:comp-AB} and \eqref{e:Case111}. 

By  {\bf (A3)} and \eqref{e:comp-AB}, we have that when  $t^{1/\alpha} < |x-y|/8$,  
\begin{align}\label{e:intJ1}
&t\int_{(x_d \vee y_d \vee t^{1/\alpha})\wedge (|x-y|/4)}^{|x-y|/2}  J(x+t^{1/\alpha}\e_d, 
		 x+r \e_d  ) J(x+ r \e_d  ,y+t^{1/\alpha}\e_d)r^{d-1} dr \\
&\asymp \frac{t}{|x-y|^{d+\alpha}}	\int_{(x_d \vee y_d \vee t^{1/\alpha})\wedge (|x-y|/4)}^{|x-y|/2}  A_{\beta_1, \beta_2,  \beta_3, \beta_4}(t,x,x+r\e_d)\,A_{\beta_1, \beta_2,  \beta_3, \beta_4}(t,x+r\e_d,y)
\frac{dr}{r^{\alpha+1}}. \nn
\end{align}
Thus, we see 
that, for $\beta_2=\alpha+\beta_1$ and $t^{1/\alpha} <  |x-y|/8$, \eqref{e:allcase}
follows from   \eqref{e:Case3}, \eqref{e:intJ1} and  
Lemma \ref{l:lower-2}, and 
that, for $\beta_2>\alpha+\beta_1$ and $t^{1/\alpha} <  |x-y|/8$, the lower bound in \eqref{e:allcase}
follows from     \eqref{e:intJ1} and 
Proposition \ref{p:HKE-lower-1} and Lemma \ref{l:lower}.

We have from  \eqref{e:intJ1}, \eqref{e:two-jump-region} and \eqref{e:comp-AB} that,  when  $t^{1/\alpha} < |x-y|/8$,  
\begin{align}\label{e:intJ2}
&t\int_{(x_d \vee y_d \vee t^{1/\alpha})\wedge (|x-y|/4)}^{|x-y|/2}  J(x+t^{1/\alpha}\e_d, 
		 x+r \e_d  ) J(x+ r \e_d  ,y+t^{1/\alpha}\e_d)r^{d-1} dr \\
&\ge \frac {c_1t}{|x-y|^{d+2\alpha}}	  B_{\beta_1, \beta_1, 0, \beta_3}(x+t^{1/\alpha}\e_d,y+t^{1/\alpha}\e_d) \L^{\beta_3}\Big( \frac{|x-y|}{((x_d \wedge y_d) \vee t^{1/\alpha}) \wedge |x-y|}\Big). \nn
\end{align}
Now,   
for $\beta_2>\alpha+\beta_1$ and $t^{1/\alpha} <  |x-y|/8$, the upper bound in \eqref{e:allcase} follows from  \eqref{e:intJ1},  
\eqref{e:intJ2}
and the upper bound in  \eqref{e:Case2}.

Finally, from the joint continuity of $\bar p(t,x,y)$ and upper heat kernel estimates, we deduce that $\rY$ is a Feller process and finish the proof by Remark \ref{r:conti}. 
\qed

 \noindent \textsc{Proof of Theorem  \ref{t:HKE-kappa}.}  The second comparison in \eqref{e:factor} follows from Corollary \ref{c:life-2}. By \eqref{e:comp-AB}, Theorem \ref{t:HKE} and  Propositions \ref{p:HKE-lower-1}, \ref{p:lower-1} and \ref{p:lower-2},  $p^\kappa(t,x,y) \ge c_1  \big(1 \wedge \frac{x_d}{t^{1/\alpha}}\big)^{q_\kappa} \big(1 \wedge \frac{ y_d}{t^{1/\alpha}}\big)^{q_\kappa} \bar p(t,x,y)$ 
 for all $(t,x,y) \in (0,\infty) \times \R^d_+ \times \R^d_+$. On the other hand,  by \eqref{e:comp-AB}, Theorem \ref{t:HKE}, Corollary \ref{c:UHK}, and  Theorem \ref{t:UHK2}, 
  $p^\kappa(t,x,y) \le c_2 \big(1 \wedge \frac{x_d}{t^{1/\alpha}}\big)^{q_\kappa} \big(1 \wedge \frac{ y_d}{t^{1/\alpha}}\big)^{q_\kappa} \bar p(t,x,y)$ 
 for all $(t,x,y) \in (0,\infty) \times \R^d_+ \times \R^d_+$ and hence  \eqref{e:factor} holds true.
 
Note that for each $(t,y) \in (0,\infty) \times \R^d_+$, the map $x\mapsto p^\kappa(t,x,y)$ vanishes continuously on $\partial \R^d_+$. Hence,  using the joint continuity of $p^\kappa(t,x,y)$ and upper heat kernel estimates, we deduce that $Y^\kappa$ is a Feller process. By Remark \ref{r:conti}, the proof is complete.  
\qed

\section{Green function estimates}\label{s:Green}

In this section, we give proofs of Theorems \ref{t:green-1} and \ref{t:green-2}.

\medskip

\noindent \textsc{Proof of Theorem \ref{t:green-1}.} When $d>\alpha$, we get the upper bound of \eqref{e:green-reflect} from  Corollary \ref{c:green-upper}. On the other hand, by Lemma \ref{l:NDL-bar-1} and Remark \ref{r:conti}, we have
\begin{align*}
	\bar G(x,y) \ge  \int_{|x-y|^\alpha}^\infty \bar p(t,x,y) dt \ge c_1 \int_{|x-y|^\alpha}^\infty  t^{-d/\alpha}dt = \begin{cases}
		\frac{c_1\alpha}{d-\alpha}|x-y|^{-d+\alpha} \;\; &\mbox{if } d>\alpha;\\[2pt]
		\infty  \;\; &\mbox{if } d\le \alpha.
	\end{cases}
\end{align*}
The proof is complete. 
\qed

In the remainder of this section, we assume the setting of Theorem \ref{t:green-2} holds and  denote by $q$ the constant $q_\kappa$ in \eqref{e:killing-potential}, which is strictly positive. 

Let $x,y \in \R^d_+$ be such that $x_d \le y_d$ and $|x-y|=1$. From Theorem \ref{t:UHK},  Proposition \ref{p:HKE-lower-1} and \eqref{e:wtB-interior}, we have for $t \le 1$,
\begin{align}\label{e:green-small-1}
	p^\kappa(t,x,y)& \le C t \left( 1 \wedge \frac{x_d}{t^{1/\alpha}}\right)^{q}   \left( 1 \wedge \frac{y_d}{t^{1/\alpha}}\right)^{q} \big((x_d \vee t^{1/\alpha}) \wedge 1\big)^{\beta_1}
		\\
	&\;\; \times  
	\big((y_d \vee t^{1/\alpha}\big)\wedge 1\big)^{\beta_2 \wedge (\alpha/2+\beta_1)} \L^{\beta_3} \Big( \frac{(y_d \vee t^{1/\alpha}) \wedge 1}{(x_d \vee t^{1/\alpha}) \wedge 1}\Big) \L^{\beta_4}\Big( \frac{1}{(y_d \vee t^{1/\alpha})}	\Big),\nn
\end{align}
\begin{align}\label{e:green-small-2}
	p^\kappa(t,x,y) &\ge c t \left( 1 \wedge \frac{x_d}{t^{1/\alpha}}\right)^{q}   \left( 1 \wedge \frac{y_d}{t^{1/\alpha}}\right)^{q} \big((x_d \vee t^{1/\alpha}) \wedge 1\big)^{\beta_1}\big((y_d \vee t^{1/\alpha}\big)\wedge 1\big)^{\beta_2}	\\
	&\quad \times  \L^{\beta_3} \Big( \frac{(y_d \vee t^{1/\alpha}) \wedge 1}{(x_d \vee t^{1/\alpha}) \wedge 1}\Big) \L^{\beta_4} \Big( \frac{1}{(y_d \vee t^{1/\alpha})}	\Big),\nn
\end{align}
and for $t>1$,
\begin{align}\label{e:green-large}
	p^\kappa(t,x,y) \asymp t^{-d/\alpha} \left( 1 \wedge \frac{x_d}{t^{1/\alpha}}\right)^q   \left( 1 \wedge \frac{y_d}{t^{1/\alpha}}\right)^q.
\end{align}

\begin{lemma}\label{l:int-smalltime}
Let $x,y \in \R^d_+$ be such that $x_d \le y_d$ and $|x-y|=1$. Set $\wh q:=2\alpha+\beta_1+\beta_2-q$. Then  we have
	\begin{equation*}
		\int_0^1 p^\kappa(t,x,y)dt \asymp \begin{cases}
			(x_d \wedge 1)^q (y_d \wedge 1)^{q} &\mbox{ if } q <\wh q,\\[2pt]
			\displaystyle	(x_d \wedge 1)^q (y_d \wedge 1)^{q} \L^{\beta_4 +1} (1/y_d) &\mbox{ if } q =\wh q,\\[2pt]
				\displaystyle(x_d \wedge 1)^q (y_d \wedge 1)^{\wh q} \L^{\beta_4}(1/y_d)&\mbox{ if } q > \wh q,
		\end{cases}
	\end{equation*}
where the comparison constant is independent of $x,y$. 
In particular, when $d=1$, we have $\int_0^1 p^\kappa(t,x,y) dt\asymp (x \wedge 1)^q$.
\end{lemma}
\begin{proof} 
Set $G_1:=	\int_0^1 p^\kappa(t,x,y)dt$. 
Note that, when $d=1$, we have $y = 1+x\ge 1$.
  
\textbf{(Case 1)} $x_d \ge 1$: We have from \eqref{e:green-small-1} and \eqref{e:green-small-2} 
that $G_1\asymp \int_0^1  tdt \asymp 1.$

\textbf{(Case 2)} $y_d \ge 1>x_d$:  By  \eqref{e:green-small-2} and \eqref{e:equiv-form},
\begin{align*}
	G_1&\ge c_1x_d^{\beta_1}\int_{1/2}^1 t \left( 1 \wedge \frac{x_d}{t^{1/\alpha}}\right)^{q-\beta_1} dt \ge  2^{-1}c_1x_d^{\beta_1}\int_{1/2}^1 (1 \wedge x_d)^{q-\beta_1} dt = 2^{-2}c_1x_d^q.
\end{align*}
Besides, since $\alpha+\beta_1-q>0$, we get from \eqref{e:green-small-1} and \eqref{e:equiv-form} that
\begin{align*}
	G_1\le c_2x_d^{\beta_1} \int_0^1  t \Big( 1\wedge \frac{x_d}{t^{1/\alpha}}\Big)^{q-\beta_1}  
	\L^{\beta_3} (t^{-1/\alpha}) dt \le c_3x_d^{q}\int_{0}^1 t^{(\alpha+\beta_1-q)/\alpha}  \L^{\beta_3} (t^{-1/\alpha}) dt \le  c_4x_d^q.
\end{align*}

\textbf{(Case 3)}  
$1>y_d \ge x_d$  and $d\ge 2$: 
 Set $\wh \beta_2:=\beta_2 \wedge (\alpha/2+\beta_1)$. 
Note that 
\begin{align}
	\int_{y_d^\alpha}^1 t
	^{-1+(\wh q -q)/\alpha} \L^{\beta_4} (t^{-1/\alpha}) dt\asymp
 g(y_d):=
	 \begin{cases}
		1 &\mbox{ if } \,q<\wh q;\\[2pt]
		 \L^{\beta_4+1} (1/y_d) &\mbox{ if } \,q=\wh q;\\[2pt]
	y_d^{\wh q - q} \L^{\beta_4}(1/y_d) &\mbox{ if } \,q>\wh q.
	\end{cases}\label{e:Green-case3}
\end{align}

For the lower bound, we get from  \eqref{e:green-small-2} and \eqref{e:Green-case3} that
\begin{align*}
	G_1&\ge c_5 x_d^{q}y_d^{q}
	\int_{y_d^\alpha}^1 
	t^{(\alpha+\beta_1+\beta_2-2q)/\alpha} \L^{\beta_4} (t^{-1/\alpha}) dt= c_5 x_d^{q}y_d^{q}
	\int_{y_d^\alpha}^1 t
	^{-1+(\wh q -q)/\alpha} \L^{\beta_4} (t^{-1/\alpha}) dt\nn\\
	&\asymp x_d^{q}y_d^{q} 
	 g(y_d).
\end{align*}

For the upper bound, we see from \eqref{e:green-small-1}  that
\begin{align*}
	G_1&\le c_6x_d^{\beta_1} y_d^{\wh \beta_2} \L^{\beta_3}( {y_d}/{x_d}) \L^{\beta_4}({1}/{y_d}) \int_0^{x_d^\alpha} t  dt  \\
	&\quad + c_6x_d^{q} y_d^{\wh\beta_2}  \L^{\beta_4} (1/y_d)\int_{x_d^\alpha}^{y_d^\alpha} t^{(\alpha+\beta_1-q)/\alpha} \L^{\beta_3}({y_d}/{t^{1/\alpha}}) dt \\
	&\quad + c_6x_d^q y_d^q  \int_{y_d^\alpha}^1  t^{(\alpha+ \beta_1+\wh\beta_2-2q)/\alpha} \L^{\beta_4} (t^{-1/\alpha}) dt =:c_6(I+II+III).
\end{align*}
Using \eqref{e:slowly-varying}, we obtain
\begin{align}\label{e:Green-I}
	I \le c_7x_d^{2\alpha+\beta_1} y_d^{\wh \beta_2}\Big(\frac{y_d}{x_d}\Big)^{2\alpha+\beta_1-q}\L^{\beta_4} (1/y_d) =  c_7 x_d^{q} y_d^{2\alpha+\beta_1+\wh \beta_2-q} \L^{\beta_4} (1/y_d).
\end{align}
Since $\alpha+\beta_1>q$, the map  $t\mapsto t^{(\alpha+\beta_1-q)/\alpha} 
\L^{\beta_3}(y_d/t^{1/\alpha})$ is almost increasing. Therefore, 
\begin{align}\label{e:Green-II}
	II \le c_8 x_d^{q} y_d^{\alpha+\beta_1+\wh \beta_2-q} \L^{\beta_4} (1/y_d)\int_{x_d^\alpha}^{y_d^\alpha} dt \le c_8 x_d^{q} y_d^{2\alpha+\beta_1+\wh \beta_2-q} \L^{\beta_4} (1/y_d).
\end{align}

We consider the cases $q<\wh q$ and $q \ge \wh q$ separately. First, suppose that $q<\wh q$, which is equivalent to $q<\alpha + (\beta_1+\beta_2)/2$. 
 Since $q<\alpha+\beta_1$, it follows that
$q<\alpha+(\beta_1+\wh\beta_2)/2$. 
Using \eqref{e:slowly-varying}, since $y_d< 1$, we see from \eqref{e:Green-I}--\eqref{e:Green-II} that
\begin{align*}
	I+II \le c_9 x_d^q y_d^q  y_d^{2\alpha+ \beta_1+\wh \beta_2-2q}\L^{\beta_4} (1/y_d) \le c_{10} x_d^q y_d^q.
\end{align*}
Moreover, since $(\alpha+ \beta_1+\wh \beta_2-2q)/\alpha >-1$, we get
\begin{align*}
	III \le  c_7x_d^q y_d^q  \int_{0}^1  t^{(\alpha+ \beta_1+\wh\beta_2-2q)/\alpha} \L^{\beta_4} (t^{-1/\alpha}) dt \le c_{11} x_d^qy_d^q.
\end{align*}
Therefore, we arrive at the result in this case.

Suppose that $q \ge \wh q$. In this case, we have $2\alpha + \beta_1 + \beta_2 = q+\wh q \le 2q <2\alpha+2\beta_1$. Thus, $\beta_2<\beta_1$ and $\wh \beta_2=\beta_2$. Then we deduce the desired  upper bound from \eqref{e:Green-case3}--\eqref{e:Green-II}.
The proof is complete. \end{proof}

\begin{lemma}\label{l:int-largetime}
	Let $x,y \in \R^d_+$ be such that $x_d \le y_d$ and $|x-y|=1$. 
	Then  we have
	\begin{equation*}
	\int_1^\infty p^\kappa(t,x,y) dt\asymp \begin{cases}
			(x_d \wedge 1)^q (y_d \wedge 1)^{q} &\mbox{ if } d >\alpha,\\[2pt]
			(x_d \wedge 1)^q (y_d \wedge 1)^{q}
			\L(x_d) 	
			&\mbox{ if } d=1=\alpha,\\[2pt]
			\displaystyle(x_d \wedge 1)^q (y_d \wedge 1)^{q} (x_d \vee 1)^{\alpha-1}
			&\mbox{ if } d=1<\alpha,
		\end{cases}
	\end{equation*}
	where the comparison constant is independent of $x,y$. 
\end{lemma}
\begin{proof} By \eqref{e:green-large}, we have
$$	\int_1^\infty p^\kappa(t,x,y) dt \asymp \int_1^\infty t^{-d/\alpha} \left( 1 \wedge \frac{x_d}{t^{1/\alpha}}\right)^q   \left( 1 \wedge \frac{y_d}{t^{1/\alpha}}\right)^q dt =:G_2.$$

If $d>\alpha$, then by Lemma \ref{cal:green}, $G_2 \asymp (x_d \wedge 1)^q (y_d \wedge 1)^q$.

If $d=1=\alpha$, then using Lemma \ref{cal:green}, we get
\begin{align*}
	G_2 &\asymp \int_1^{ x_d \vee 1} t^{-1} dt +  x_d^q \int_{ x_d \vee 1}^\infty t^{-1-q}   \left( 1 \wedge \frac{y_d}{t}\right)^q dt \\
	&\asymp \1_{\{x_d>1\}} \log(x_d) + x_d^q ( x_d \vee 1)^{-q}   \left( 1 \wedge \frac{y_d}{x_d \vee 1}\right)^q\\
	& \asymp \begin{cases}
		\L(x_d) &\text{ if } x_d>1\\
		x_d^q (y_d \wedge 1)^q  &\text{ if } x_d\le 1
	\end{cases} \asymp (x_d \wedge 1)^q(y_d \wedge 1)^q 		\L(x_d). 	
\end{align*}

If $d=1<\alpha$, then since  $y_d>2$ implies  $x_d \ge y_d-|x-y| \ge y_d/2 >1$, and $\alpha-q-1 \le 0$,  we get
\begin{align*}
	G_2 &\asymp \1_{\{y_d>2\}} \Big( \int_1^{ x_d^\alpha}    t^{-1/\alpha}dt  + x_d^q \int_{ x_d^\alpha}^{(2y_d)^\alpha} 	t^{-(q+1)/\alpha}dt  + x_d^q y_d^q  	\int_{(2y_d)^\alpha}^\infty 	t^{-(2q+1)/\alpha} dt \Big)\nn\\
&\quad + \1_{\{y_d \le 2\}} x_d^qy_d^q \int_1^\infty t^{-(2q+1)/\alpha}dt   \\
	&\asymp \1_{\{y_d>2\}} \big(  x_d^{\alpha-1}  + x_d^{q+\alpha - (q+1)} + x_d^q y_d^{\alpha-q-1}\big) +  \1_{\{y_d\le 2\}} x_d^q y_d^q \\
	&  \asymp \1_{\{y_d>2\}} x_d^{\alpha-1} +  \1_{\{y_d\le 2\}} x_d^q y_d^q\asymp (x_d \wedge 1)^q(y_d \wedge 1)^q (x_d \vee 1)^{\alpha-1}.
\end{align*}
The proof is complete. 
\end{proof}

\noindent \textsc{Proof of Theorem \ref{t:green-2}.} From
scaling property \eqref{e:kernel-scaling}, we obtain
\begin{align}\label{e:green-scaling}
G^\kappa(x,y) = |x-y|^{-d+\alpha} \, G^\kappa(x/|x-y|, y/|x-y|), \quad x,y \in \R^d_+.
\end{align}
Using \eqref{e:green-scaling}, symmetry and Lemmas \ref{l:int-smalltime} and \ref{l:int-largetime}, we deduce the desired result. 
\qed

\section{Appendix: Estimates of some integrals}\label{appendix B}
Note that for any $\eps>0$, 
\begin{equation}\label{e:slowly-varying}
	\L(r) <
	(2+\eps^{-1}) r^\eps \quad \text{for all} \; r \ge 1,
\end{equation}
\begin{equation}\label{e:slowly-varying-2}	
\L(ar)/\L(r) <	(1+\eps^{-1}) a^{\eps} \quad  \quad \text{for all }  a\ge1 \text{ and } r>0.
\end{equation}

Recall the definition of $ A_{b_1,b_2, b_3, b_4}(t,x,y)$  from \eqref{def:wtB2}. 
\begin{lemma}\label{l:kill-log}
	Let $b_1, b_2, b_3, b_4 \ge 0$.
	
\noindent (i) If $b_1>0$, then for any $\eps \in (0, b_1]$, 
there exists $c_1>0$ such that
\begin{equation*}
	A_{b_1, b_2, b_3, b_4} (t, x,y) \le c_1	A_{b_1-\eps, b_2, 0, b_4}(t, x,y) \quad \text{ for all } t\ge 0,\, x,y \in \R^d_+ .
\end{equation*}

\noindent (ii) If $b_2>0$, then for any $\eps \in (0, b_2]$, 
there exists $c_2>0$ such that
\begin{equation*}
	A_{b_1, b_2, b_3, b_4}(t, x,y)  \le c_2	A_{b_1, b_2-\eps,  b_3,0}(t, x,y) \quad \text{ for all } t\ge 0, \, x,y \in \R^d_+.
\end{equation*}
\end{lemma}
\begin{proof} The results follow from \eqref{e:slowly-varying}. \end{proof} 

\begin{lemma}\label{l:kill-log-2}
Let $b_1, b_3 \ge 0$. Suppose that $b_1>0$ if $b_3>0$. Then  there exists a constant $C>0$ such that 	for all $t\ge 0$ and  $x,y \in \R^d_+$,
	\begin{align*}
	A_{b_1, 0, b_3, 0} (t, x,y)	&\le  C\, \Big(\frac{x_d\vee t^{1/\alpha}}{|x-y|} \wedge 1 \Big)^{b_1} \L^{b_3} \Big(\frac{ (y_d\vee t^{1/\alpha}) \wedge |x-y|}{(x_d\vee t^{1/\alpha})  \wedge |x-y|} \Big).
	\end{align*}
\end{lemma} 
\begin{proof} When $x_d\le y_d$, the result is clear. Assume that $x_d>y_d$. 
Set $r=\frac{x_d\vee t^{1/\alpha}}{|x-y|} \wedge 1$ and $s=\frac{y_d\vee t^{1/\alpha}}{|x-y|} \wedge 1$.
Then $0<s \le r \le 1$ and the desired inequality is equivalent to 
\begin{equation}\label{e:kill-log-21}	
 \L^{b_3}(r/s) /\L^{b_3}(s/r)  \le C(r/s)^{b_1}.  
\end{equation}
If $b_3=0$, then \eqref{e:kill-log-21} clearly holds with $C=1$. If $b_3>0$ and $b_1>0$, then  we get from \eqref{e:slowly-varying} that
$
\L^{b_3}(r/s) /\L^{b_3}(s/r)  \le  \L^{b_3}(r/s) \le c_1 (r/s)^{b_1}.
$
This proves the lemma. \end{proof}

\begin{remark}\label{r:kill-log}
	{\rm Recall that $B_{b_1,b_2,b_3,b_4}(x,y)=A_{b_1,b_2,b_3,b_4}(0,x,y)$ for $x,y \in \R^d_+$. 
		Thus, 
		the results of Lemmas \ref{l:kill-log} and \ref{l:kill-log-2} hold true with
	$B_{b_1,0,b_3,0}(x,y)$ instead of $A_{b_1,0,b_3,0}(t,x,y)$.
	} 
\end{remark}

\begin{lemma}\label{cal:00}	
Let $\gamma \in \R$, $b \ge 0$ and $t,k,l>0$. Suppose that either $\gamma<\alpha$ or $k\ge t^{1/\alpha}$. Then we have	$$	\int_0^t \Big(1 \wedge \frac{k}{s^{1/\alpha}} \Big)^{\gamma}  \L^{b} \Big( \frac{l}{k \vee s^{1/\alpha} } \Big)   ds \le 	Ct \Big(1 \wedge \frac{k}{t^{1/\alpha}} \Big)^{\gamma}  \L^{b} \Big( \frac{l}{k \vee t^{1/\alpha} } \Big),	$$
where $C>0$ is a constant which depends only on $\gamma$ and $b$. \end{lemma}
\begin{proof}  If $k \ge t^{1/\alpha}$, then 	the result holds since the left hand side is  $t \L^b(l/k)$.

Suppose that $k<t^{1/\alpha}$ and  $\gamma<\alpha$. 
Let $\eps>0$ be such that  $\gamma+b\eps<\alpha$. Then the left hand side is
\begin{align*}
	&\L^{b} \Big(\frac{l}{k}\Big)\int_0^{k^\alpha}  ds + k^\gamma \int_{k^\alpha}^t \frac{1}{s^{\gamma/\alpha}} \L^{b} \Big(\frac{l}{s^{1/\alpha}}\Big)  ds \le k^\alpha \L^{b} \Big(\frac{l}{k}\Big) +c_1 k^\gamma \L^{b} \Big(\frac{l}{t^{1/\alpha}}\Big)\int_{k^\alpha}^t \frac{	t^{b\eps/\alpha}}{ s^{\gamma/\alpha + b\eps/\alpha}} ds\\
	&\le c_2 k^\alpha \Big(\frac{t^{1/\alpha}}{k}\Big)^{\alpha-\gamma} \L^{b} \Big(\frac{l}{t^{1/\alpha}}\Big) +  c_2 k^\gamma t^{1-\gamma/\alpha}  \L^{b} \Big(\frac{l}{t^{1/\alpha}}\Big) \\
	&= 2c_2 t \Big( \frac{k}{t^{1/\alpha}}\Big)^\gamma\L^{b} \Big(\frac{l}{t^{1/\alpha}}\Big) = 2c_2 t \Big( 1 \wedge  \frac{k}{t^{1/\alpha}}\Big)^\gamma\L^{b} \Big(\frac{l}{t^{1/\alpha}}\Big).
	\end{align*}
We used \eqref{e:slowly-varying-2} in both inequalities above.\end{proof}

\begin{lemma}\label{cal:0}
	Let $b_1,b_2 \in \R$, $b_3,b_4\ge 0$ and $t,x_d,y_d>0$.
	Suppose that (1) either $b_1>q-\alpha$ or $x_d\ge t^{1/\alpha}$, and (2) either $b_2>q-\alpha$ or $y_d\ge t^{1/\alpha}$. Then we have
	\begin{align*}
&	\int_0^t \Big(1 \wedge \frac{x_d}{(t-s)^{1/\alpha}} \Big)^{q-b_1} \Big(1 \wedge \frac{y_d}{s^{1/\alpha}} \Big)^{q- b_2}  \L^{b_3} \Big( \frac{ y_d \vee s^{1/\alpha}}{ x_d \vee (t-s)^{1/\alpha} } \Big) \L^{b_4} \Big( \frac{1}{ y_d \vee s^{1/\alpha} } \Big)  ds\\
&\le C t\Big(1 \wedge \frac{x_d}{t^{1/\alpha}} \Big)^{q-b_1} \Big(1 \wedge \frac{y_d}{t^{1/\alpha}} \Big)^{q- b_2}  \L^{b_3} \Big( \frac{ y_d \vee t^{1/\alpha}}{ x_d \vee t^{1/\alpha} } \Big) \L^{b_4} \Big( \frac{1}{ y_d \vee t^{1/\alpha} } \Big)
	\end{align*}
where $C>0$ is a constant which depends only on $b_1,b_2,b_3$ and $b_4$.
\end{lemma}
\begin{proof} Observe that 
\begin{align*}
	&\int_0^{\frac{t}{2}}  \Big(1 \wedge \frac{x_d}{(t-s)^{1/\alpha}} \Big)^{q-b_1} \Big(1 \wedge \frac{y_d}{s^{1/\alpha}} \Big)^{q-b_2}\L^{b_3} \Big( \frac{ y_d \vee s^{1/\alpha}}{x_d \vee (t-s)^{1/\alpha} } \Big)    \L^{b_4} \Big( \frac{1}{y_d \vee  s^{1/\alpha} } \Big)  ds\nn\\
	&\le c_1 \Big(1 \wedge \frac{x_d}{t^{1/\alpha}} \Big)^{q-b_1} \L^{b_3} \Big( \frac{ y_d \vee t^{1/\alpha}}{x_d \vee t^{1/\alpha} } \Big)   \int_0^{\frac{t}{2}}   \Big(1 \wedge \frac{y_d}{s^{1/\alpha}} \Big)^{q-b_2} \L^{b_4} \Big( \frac{1}{y_d \vee  s^{1/\alpha} } \Big)  ds
\end{align*}
and
\begin{align*}
	&\int_{\frac{t}{2}}^t  \Big(1 \wedge \frac{x_d}{(t-s)^{1/\alpha}} \Big)^{q-b_1} \Big(1 \wedge \frac{y_d}{s^{1/\alpha}} \Big)^{q-b_2}\L^{b_3} \Big( \frac{ y_d \vee s^{1/\alpha}}{x_d \vee (t-s)^{1/\alpha} } \Big)    \L^{b_4} \Big( \frac{1}{y_d \vee  s^{1/\alpha} } \Big)  ds\nn\\
	&\le c_2\Big(1 \wedge \frac{y_d}{t^{1/\alpha}} \Big)^{q-b_2} \L^{b_4} \Big( \frac{1}{y_d \vee  t^{1/\alpha} } \Big)   \int_0^{\frac{t}{2}}   \Big(1 \wedge \frac{x_d}{s^{1/\alpha}} \Big)^{q-b_1} \L^{b_3} \Big( \frac{ y_d \vee t^{1/\alpha}}{x_d \vee s^{1/\alpha} } \Big)    ds.
\end{align*}
Using Lemma \ref{cal:00} twice, we arrive at the result. \end{proof}

\begin{lemma}\label{l:cal1}
Let $b_1,b_2 \in \R$, $b_3,b_4\ge 0$ and $t,x_d,y_d>0$. Suppose that (1) either $b_1>q-\alpha$ or $x_d\ge t^{1/\alpha}$, and (2) either $b_2>q-\alpha$ or $y_d\ge t^{1/\alpha}$. Then we have 
that for $y_d \vee t^{1/\alpha} <2$, 
\begin{align*}
	&\int_0^{t}  s\Big(1 \wedge \frac{x_d}{(t-s)^{1/\alpha}} \Big)^{q-  b_1}  \Big(1 \wedge \frac{y_d}{s^{1/\alpha}} \Big)^{q- b_2}   \\
	&\quad\times   \int_{y_d \vee s^{1/\alpha}}^2  \L^{b_3} \Big( \frac{r}{y_d \vee s^{1/\alpha}} \Big)\L^{b_3} \Big( \frac{r}{x_d \vee (t-s)^{1/\alpha}} \Big) \L^{b_4} \Big( \frac{1}{r}\Big) \frac{dr}{r} \, ds \\
	& \le C t \Big(1 \wedge \frac{x_d}{t^{1/\alpha}} \Big)^{q-  b_1}  \Big(1 \wedge \frac{y_d}{t^{1/\alpha}} \Big)^{q-b_2 } \int_{y_d}^{2}(r^\alpha \wedge t)\L^{b_3} \Big( \frac{r}{x_d \vee t^{1/\alpha}} \Big)   \L^{b_3} \Big( \frac{r}{y_d \vee t^{1/\alpha}} \Big)\L^{b_4} \Big( \frac{1}{r}\Big)    \frac{dr}{r},
\end{align*}
where $C>0$ is a constant which depends only on $b_1,b_2,b_3$ and $b_4$. 
\end{lemma}
\begin{proof} Using Fubini's theorem and Lemma \ref{cal:00} twice as in the proof of Lemma \ref{cal:0}, we get
\begin{align*}
	&\int_0^{t}  s\Big(1 \wedge \frac{x_d}{(t-s)^{1/\alpha}} \Big)^{q-b_1}  \Big(1 \wedge \frac{y_d}{s^{1/\alpha}} \Big)^{q-b_2}    \\
	&\quad\times   \int_{y_d \vee s^{1/\alpha}}^2  \L^{b_3} \Big( \frac{r}{y_d \vee s^{1/\alpha}} \Big)\L^{b_3} \Big( \frac{r}{x_d \vee (t-s)^{1/\alpha}} \Big) \L^{b_4} \Big( \frac{1}{r}\Big) \frac{dr}{r} \, ds \\
	= &\int_{y_d}^{2}\int_{0}^{r^\alpha \wedge t}    s\Big(1 \wedge \frac{x_d}{(t-s)^{1/\alpha}} \Big)^{q-b_1}  \Big(1 \wedge \frac{y_d}{s^{1/\alpha}} \Big)^{q- b_2 } \\
	&\quad\times    \L^{b_3} \Big( \frac{r}{y_d \vee s^{1/\alpha}} \Big)\L^{b_3} \Big( \frac{r}{x_d \vee (t-s)^{1/\alpha}} \Big) \L^{b_4} \Big( \frac{1}{r}\Big)   \, ds \frac{dr}{r}\\
	\le &\int_{y_d}^{2} (r^\alpha \wedge t)\int_{0}^{ t}  \Big(1 \wedge \frac{x_d}{(t-s)^{1/\alpha}} \Big)^{q-b_1}  \Big(1 \wedge \frac{y_d}{s^{1/\alpha}} \Big)^{q-b_2} \\
	&\quad\times    \L^{b_3} \Big( \frac{r}{y_d \vee s^{1/\alpha}} \Big)\L^{b_3} \Big( \frac{r}{x_d \vee (t-s)^{1/\alpha}} \Big) \L^{b_4} \Big( \frac{1}{r}\Big)   \, ds \frac{dr}{r}\\
		\le &\, C t \Big(1 \wedge \frac{x_d}{t^{1/\alpha}} \Big)^{q-b_1}  \Big(1 \wedge \frac{y_d}{t^{1/\alpha}} \Big)^{q- b_2 } \int_{y_d}^{2}(r^\alpha \wedge t)\L^{b_3} \Big( \frac{r}{x_d \vee t^{1/\alpha}} \Big)   \L^{b_3} \Big( \frac{r}{y_d \vee t^{1/\alpha}} \Big)\L^{b_4} \Big( \frac{1}{r}\Big)    \frac{dr}{r}.
\end{align*}
\end{proof}

\begin{lemma}\label{cal:new1} There is a constant $C>0$ such that for all $x \in \R^d_+$ and $A>0$,
	\begin{align*}
		\int_{B_+(x, A)}\frac{dz}{z_d^{1/2}} \le CA^{d} (x_d \vee A)^{-1/2}.
	\end{align*}
\end{lemma}
\begin{proof} We have
\begin{align}\label{e:new1-1}		
\int_{B_+(x, A)}
\frac{dz}{z_d^{1/2}}  \le 	\int_{\wt z \in \R^{d-1}, \, |\wt x-\wt z| < A} d \wt z \int_{(x_d-A) \vee 0}^{x_d+A}  \frac{ds}{s^{1/2}}  \le c_1 A^{d-1} \int_{(x_d-A) \vee 0}^{x_d+A}  \frac{ds}{s^{1/2}} .
\end{align}
If $x_d \ge 2A$, then
\begin{align}\label{e:new1-2}
	 \int_{(x_d-A) \vee 0}^{x_d+A}  \frac{ds}{s^{1/2}}  \le  \frac{1}{(x_d/2)^{1/2} }  \int_{x_d-A}^{x_d+A}ds =2^{3/2}A x_d^{-1/2}.
\end{align} 
If $x_d<2A$, then 
\begin{align}\label{e:new1-3}
	\int_{(x_d-A) \vee 0}^{x_d+A} \frac{ds}{s^{1/2}}  \le  \int_{0}^{4A}\frac{ds}{s^{1/2}}  = 4A^{1/2}.
\end{align} 
Combining \eqref{e:new1-1} with \eqref{e:new1-2}--\eqref{e:new1-3}, we arrive at the result. \end{proof}

\begin{lemma}\label{cal:new2}
\noindent (i) 	There is a constant $C>0$ such that for all $x \in \R^d_+$ and $0<A < x_d$,
\begin{align*}
	\int_{z \in \R^d_+, \, x_d \ge |x-z|>A} \frac{dz}{z_d^{1/2} |x-z|^{d+\alpha}} \le Cx_d^{-1/2}A^{-\alpha}.
\end{align*}
	
	\noindent (ii) 	Let $\eps \in (0,1)$ and $\delta >0$. There is a constant $C'>0$ such that for all $x \in \R^d_+$ and $A \ge x_d$,
	\begin{align*}
		\int_{z \in \R^d_+, \, |x-z|>A} \frac{dz}{z_d^{\eps}\, |x-z|^{d+\delta}} \le C'A^{-\eps-\delta }.
	\end{align*}
\end{lemma}
\begin{proof} Without loss of generality, we assume $x=(\wt 0, x_d)$.

(i) Note that
\begin{align*}
	&\int_{z \in \R^d_+, \, x_d \ge|x-z|>A} \frac{dz}{z_d^{1/2} |x-z|^{d+\alpha}} \nn\\
	&= \int_{z \in \R^d_+, \, x_d \ge|x-z|>A, \, |\wt z| \le |x_d-z_d|} \frac{dz}{z_d^{1/2} |x-z|^{d+\alpha}} + \int_{z \in \R^d_+, \, x_d \ge|x-z|>A, \, |\wt z|>|x_d-z_d|} \frac{dz}{z_d^{1/2} |x-z|^{d+\alpha}}\nn\\
	&=:I + II. 
\end{align*}
First, we have
\begin{align*}
	I &\le  \int_{x_d \ge |x_d-z_d|>\frac{A}{2}} \frac{1}{z_d^{1/2} |x_d-z_d|^{d+\alpha}} \int_{\wt z\in \R^{d-1}, \, |\wt z| \le |x_d-z_d|} d \wt zdz_d\\
	&\le  c_1\Big(\int_{0}^{\frac{x_d}{2}} \frac{dz_d}{z_d^{1/2} |x_d-z_d|^{1+\alpha}}  +\int_{\frac{x_d}{2}}^{x_d- \frac{A}{2}} \frac{dz_d}{z_d^{1/2} |x_d-z_d|^{1+\alpha}}  +  \int_{x_d+\frac{A}{2}}^{2x_d} \frac{dz_d}{z_d^{1/2} |x_d-z_d|^{1+\alpha}}  \Big) \\
		&\le  c_1\Big(\frac{2^{1+\alpha}}{x_d^{1+\alpha}}\int_{0}^{\frac{x_d}{2}} {z_d^{-1/2}} dz_d +\frac{1}{(x_d/2)^{1/2}}\int_{\frac{x_d}{2}}^{x_d- \frac{A}{2}} \frac{dz_d}{ |x_d-z_d|^{1+\alpha}}   + {x_d^{-1/2}} \int_{x_d+\frac{A}{2}}^{2x_d} \frac{dz_d}{|x_d-z_d|^{1+\alpha}}  \Big) \nn\\
		&\le c_2 \left( x_d^{-1/2-\alpha} + x_d^{-1/2} A^{-\alpha} + x_d^{-1/2}A^{-\alpha} \right)\le c_3x_d^{-1/2} A^{-\alpha}.
\end{align*}
On the other hand, we see that for any $z \in \R^d_+$ with $x_d \ge|x-z|>A$ and $|\wt z|>|x_d-z_d|$,
\begin{align*}
	z_d \ge x_d - |x_d-z_d| \ge x_d - \frac{1}{2} (|\wt z| + |x_d - z_d|)  \ge x_d - \frac{1}{\sqrt 2} |x - z| \ge \left(1- \frac{1}{\sqrt 2}\right)x_d.
\end{align*}
Hence, we also have that
\begin{align*}
	II&\le  c_4 x_d^{-1/2}\int_{z \in \R^d_+, \,|x-z|>A} \frac{dz}{|x-z|^{d+\alpha}} \le c_5 x_d^{-1/2}A^{-\alpha}.
\end{align*}

(ii) Observe that
\begin{align*}
	 	&\int_{z \in \R^d_+, \, |x-z|>A} \frac{dz}{z_d^{\eps}\, |x-z|^{d+\delta}} \nn\\
	 	&\le 	\int_{z \in \R^d_+, \,  |x_d-z_d|\ge |\wt z| \vee \frac{A}{2}} \frac{d \wt z\,dz_d}{z_d^{\eps}\, |x_d-z_d|^{d+\delta}} + 	\int_{z \in \R^d_+, \, |\wt z|\ge |x_d-z_d| \vee \frac{A}{2}} \frac{d \wt z\,dz_d}{z_d^{\eps}\, |\wt z|^{d+\delta}}=:I + II. 
\end{align*}
Using Fubini's  theorem, since $A \ge x_d$, we see that
\begin{align*}
	II&\le \int_{\wt z \in \R^{d-1}, \, |\wt z|\ge \frac{A}{2}} \frac{1}{|\wt z|^{d+\delta}}  \int_0^{|\wt z|+x_d} z_d^{-\eps }dz_d\, d \wt z \le c_1  \int_{\wt z \in \R^{d-1}, \, |\wt z|\ge \frac{A}{2}} \frac{(|\wt z|+A)^{1-\eps}}{|\wt z|^{d+\delta}} d \wt z\\
	&\le c_2  \int_{\wt z \in \R^{d-1}, \, |\wt z|\ge \frac{A}{2}} \frac{d \wt z }{|\wt z|^{d-1+ \eps + \delta }} \le c_3 A^{-\eps-\delta}.
\end{align*}
On the other hand, using 
$\int_{\wt z \in \R^{d-1}, \, |\wt z|\le |x_d-z_d|} d\wt z \le c_4 |x_d-z_d|^{d-1}$, we also see that  
\begin{align*}
	I &\le c_5\Big( \int_{0}^{(x_d-\frac{A}{2}) \vee 0} \frac{dz_d}{z_d^{\eps}\,|x_d-z_d|^{1+\delta}} +  \int_{x_d+\frac{A}{2}}^{\infty} \frac{dz_d}{z_d^{\eps}\,|x_d-z_d|^{1+\delta}} \Big) \nn\\
	&\le c_6\Big( \frac{1}{A^{1+\delta }} \int_{0}^{(x_d-\frac{A}{2}) \vee 0} z_d^{-\eps}dz_d + \frac{1}{A^{\eps}} \int_{x_d+\frac{A}{2}}^{\infty} \frac{dz_d}{|x_d-z_d|^{1+\delta}} \Big)\nn\\
		&\le c_7\Big( \frac{1}{A^{1+\delta }}\Big(\big(x_d - \frac{A}{2} \big) \vee 0 \Big)^{1-\eps} + \frac{1}{A^{\eps+\delta}} \Big) \le   \frac{c_6}{A^{\eps + \delta }},
\end{align*}
where we used the fact that $A\ge x_d$ in the last inequality. The proof is complete. \end{proof}

For $\gamma, \eta_1, \eta_2 \ge 0$ 
and $k,l>0$, define
\begin{equation*}
	f_{\gamma,\eta_1,\eta_2,k,l}(r):=r^\gamma \L^{\eta_1}\Big(\frac{k}{r} \Big) \L^{\eta_2}\Big( \frac{r}{l} \Big).
\end{equation*}

\begin{lemma}\label{cal:basic}
	Let $\gamma, \eta_1, \eta_2 \ge 0$.
	
	\noindent (i) For any $\eps>0$, there exist constants $C,C'>0$ such that for any $k,l,r>0$ and any $a \ge 1$,  
	\begin{equation}\label{e:cal:basic}
		C a^{\gamma-\eps} \le 	\frac{f_{\gamma,\eta_1,\eta_2,k,l}(ar)}{ f_{\gamma,\eta_1,\eta_2,k,l}(r)} \le C' a^{\gamma+\eps}. 
	\end{equation}

	\noindent (ii) Assume that $\gamma>0$. Then  there exists a constant $C>0$  such that for any $k,l,r>0$ and any $a \ge 1$,
	\begin{equation*}
	\frac{f_{\gamma,\eta_1,\eta_2,k,l}(ar)}{ f_{\gamma,\eta_1,\eta_2,k,l}(r)} \ge C.  
	\end{equation*}
\end{lemma}
\begin{proof} (i)
If $\eta_2=0$, the second inequality in \eqref {e:cal:basic} is true with any $\eps \ge 0$. In case $\eta_2>0$, 
for any given $\eps>0$, let $\eps':=\eps/\eta_2$. We get from \eqref{e:slowly-varying-2} that for all $a \ge 1$ and $r>0$,
\begin{align*}
	\frac{f_{\gamma,\eta_1,\eta_2,k,l}(ar)}{ f_{\gamma,\eta_1,\eta_2,k,l}(r)} \le a^\gamma \left( \frac{\L(ar/l)}{\L(r/l)}\right)^{\eta_2 } \le a^{\gamma} \big((1+1/\eps')a^{\eps'} \big)^{\eta_2}=c(\eps, \eta_2)a^{\gamma+\eps}.
\end{align*}
The first inequality can be proved by a similar argument.

\noindent(ii) 
The desired result follows from the first inequality in \eqref{e:cal:basic} with $\eps=\gamma$. 
\end{proof}

\begin{lemma}\label{cal:2}
	Let $b_1,b_2,\eta_1,\eta_2, \gamma \ge 0$. 
	 There exists a constant $C>0$ such that for any 
	$x \in \R^d_+$ and $s,k,l>0$,
	\begin{align*}
		&\int_{B_+(x, 2)}
		\Big(1 \wedge \frac{ x_d \vee s^{1/\alpha} }{|x-z|} \Big)^{b_1}  \L^{b_2}\Big(\frac{|x-z|}{ (x_d \vee s^{1/\alpha}) \wedge |x-z|} \Big)\\
		&\hspace{1.8cm}\times \Big( s^{-d/\alpha} \wedge \frac{s}{|x-z|^{d+\alpha}}\Big)z_d^\gamma \L^{\eta_1} \Big( \frac{k}{z_d} \Big) \L^{\eta_2} \Big( \frac{z_d}{l} \Big) dz\\
		& \le C(x_d \vee s^{1/\alpha})^\gamma \L^{\eta_1} \Big( \frac{k}{x_d \vee s^{1/\alpha}} \Big) \L^{\eta_2} \Big( \frac{x_d \vee s^{1/\alpha}}{l} \Big)\\
		&\;\;   + C 
             \1_{\{x_d \vee s^{1/\alpha} <2  \}} 
		s(x_d \vee s^{1/\alpha})^{b_1}\Big[
		\1_{\{\gamma>\alpha+b_1 \}}    \L^{b_2}\Big(\frac{  2}{x_d \vee s^{1/\alpha} } \Big) \L^{\eta_1} (k ) \L^{\eta_2} \Big( \frac{1}{l} \Big) \\
		&\hspace{1.9cm} + \1_{\{\gamma=\alpha+b_1\}} \int_{x_d \vee s^{1/\alpha}}^2   \L^{b_2}\Big(\frac{  r}{x_d \vee s^{1/\alpha} } \Big) \L^{\eta_1} \Big(\frac{ k}{r} \Big) \L^{\eta_2} \Big(\frac{r}{l} \Big)  \frac{dr	}{r}
		\Big].
	\end{align*}
\end{lemma}
\begin{proof}
 Using the triangle inequality, we see that for any $z \in \R^d_+$,
\begin{equation}\label{e:cal2-triangle}
	z_d \le x_d + |x-z| \le 2(x_d \vee |x-z|).
\end{equation}
Therefore, using  Lemma \ref{cal:basic}(i)-(ii) and Lemma \ref{cal:new1}, we get that
\begin{align*}
	&\int_{B_+(x, s^{1/\alpha})}	 \Big(1 \wedge \frac{ x_d \vee s^{1/\alpha} }{|x-z|} \Big)^{b_1}  \L^{b_2}\Big(\frac{|x-z|}{ (x_d \vee s^{1/\alpha}) \wedge |x-z|} \Big)\, s^{-d/\alpha} f_{\gamma,\eta_1,\eta_2,k,l}(z_d) dz\nn\\
	&=\L^{b_2}(1)s^{-d/\alpha}
	\int_{B_+(x, s^{1/\alpha})}
	z_d^{-1/2} f_{\gamma + \frac12,\eta_1,\eta_2,k,l}(z_d) dz\nn\\
	&\le 	c_1s^{-d/\alpha} f_{\gamma + \frac12,\eta_1,\eta_2,k,l}(2(x_d  \vee s^{1/\alpha}))  
	\int_{B_+(x, s^{1/\alpha})}
	\frac{dz}{z_d^{1/2}} \nn\\
	&\le c_2 (x_d \vee s^{1/\alpha})^{-1/2} f_{\gamma + \frac12,\eta_1,\eta_2,k,l}(x_d  \vee s^{1/\alpha}) =  c_2  f_{\gamma,\eta_1,\eta_2,k,l}(x_d  \vee s^{1/\alpha}).
\end{align*}
When $x_d>s^{1/\alpha}$,  we  get from \eqref{e:cal2-triangle}, Lemma \ref{cal:basic}(i)-(ii) and Lemma \ref{cal:new2}(i) that 
\begin{align*}
	&	s\int_{z \in \R^d_+, \,  x_d\ge |x-z| > s^{1/\alpha}}\Big(1 \wedge \frac{x_d\vee s^{1/\alpha} }{|x-z|} \Big)^{b_1}  \L^{b_2}\Big(\frac{  |x-z|}{(x_d  \vee s^{1/\alpha}) \wedge |x-z|} \Big)  \frac{f_{\gamma,\eta_1,\eta_2,k,l}(z_d)}{|x-z|^{d+\alpha}} dz \nn\\
	&= \L^{b_2}(1)	s
	\int_{z \in \R^d_+, \,  x_d\ge |x-z| > s^{1/\alpha}} \frac{f_{\gamma+\frac12,\eta_1,\eta_2,k,l}(z_d)}{z_d^{1/2}|x-z|^{d+\alpha}} dz \nn\\
	&\le  c_3	s f_{\gamma + \frac12,\eta_1,\eta_2,k,l}(2x_d)
	\int_{z \in \R^d_+, \,  x_d\ge |x-z| > s^{1/\alpha}} \frac{1}{z_d^{1/2}|x-z|^{d+\alpha}} dz \nn\\
	&\le c_4x_d^{-1/2}f_{\gamma + \frac12,\eta_1,\eta_2,k,l}(x_d) = c_4f_{\gamma,\eta_1,\eta_2,k,l}(x_d \vee s^{1/\alpha}).
\end{align*}

It remains to bound the integral over $\{z \in \R^d_+:x_d \vee s^{1/\alpha}<|x-z|<2\}$
under the assumption $x_d \vee s^{1/\alpha} <2$. 
 For this, we consider the following three cases separately.

\smallskip

(i) Case $\gamma<\alpha+b_1$: Fix $\eps\in(0,1)$ such that $\gamma+3\eps<\alpha+b_1$. Using \eqref{e:cal2-triangle}, \eqref{e:slowly-varying}, Lemma \ref{cal:basic}(i)-(ii) and Lemma \ref{cal:new2}(ii), we get
\begin{align*}
	&	s\int_{z \in \R^d_+, \, |x-z| > x_d \vee s^{1/\alpha}} \Big(1 \wedge \frac{x_d\vee s^{1/\alpha} }{|x-z|} \Big)^{b_1}  \L^{b_2}\Big(\frac{  |x-z|}{(x_d  \vee s^{1/\alpha}) \wedge |x-z|} \Big)  \frac{f_{\gamma,\eta_1,\eta_2,k,l}(z_d)}{|x-z|^{d+\alpha}} dz \\
		&=	s(x_d \vee s^{1/\alpha})^{b_1}\int_{z \in \R^d_+, \, |x-z| > x_d \vee s^{1/\alpha}} \L^{b_2}\Big(\frac{  |x-z|}{x_d  \vee s^{1/\alpha}} \Big)  \frac{f_{\gamma+\eps,\eta_1,\eta_2,k,l}(z_d)}{z_d^{\eps}|x-z|^{d+\alpha + b_1}} dz \\
	&\le c_5 s(x_d \vee s^{1/\alpha})^{b_1}
	\int_{z \in \R^d_+, \, |x-z| >x_d \vee s^{1/\alpha}}   \L^{b_2}\Big(\frac{  |x-z|}{x_d \vee s^{1/\alpha} } \Big)  \frac{f_{\gamma + \eps,\eta_1,\eta_2,k,l}(2|x-z|)	}{z_d^{\eps}|x-z|^{d+\alpha+b_1}}  dz\\
	&\le c_6 s(x_d \vee s^{1/\alpha})^{b_1}
	\int_{z \in \R^d_+, \, |x-z| >x_d \vee s^{1/\alpha}}  
	 \Big(\frac{  |x-z|}{x_d \vee s^{1/\alpha} } \Big)^{\eps + \gamma+2\eps} \frac{f_{\gamma+\eps ,\eta_1,\eta_2,k,l}(x_d \vee s^{1/\alpha})	}{z_d^{\eps}|x-z|^{d+\alpha+b_1}}  dz\\
	&= c_6 s(x_d \vee s^{1/\alpha})^{b_1-\gamma-3\eps}   f_{\gamma+\eps ,\eta_1,\eta_2,k,l}(x_d \vee s^{1/\alpha})
	\int_{z \in \R^d_+, \, |x-z| >x_d \vee s^{1/\alpha}}  \frac{dz}{z_d^\eps |x-z|^{d+\alpha+b_1-\gamma-3\eps}}\\
	&\le c_7 s(x_d \vee s^{1/\alpha})^{-\alpha} (x_d \vee s^{1/\alpha})^{-\eps}   f_{\gamma+\eps ,\eta_1,\eta_2,k,l}(x_d \vee s^{1/\alpha}) \le c_7 f_{\gamma,\eta_1,\eta_2,k,l}(x_d \vee s^{1/\alpha}). 
\end{align*}

(ii) Case $\gamma>\alpha+b_1$: Fix $\eps>0$ such that $\gamma-\eps>\alpha+b_1$. Using  \eqref{e:cal2-triangle}, Lemma \ref{cal:basic}(i)-(ii) and \eqref{e:slowly-varying-2}, we get
\begin{align*}
	&	s\int_{z \in \R^d_+, \, x_d \vee s^{1/\alpha}<|x-z|<2} \Big(1 \wedge \frac{x_d\vee s^{1/\alpha} }{|x-z|} \Big)^{b_1}  \L^{b_2}\Big(\frac{  |x-z|}{(x_d  \vee s^{1/\alpha}) \wedge |x-z|} \Big)  \frac{f_{\gamma,\eta_1,\eta_2,k,l}(z_d)}{|x-z|^{d+\alpha}} dz \\
	&\le c_8 s(x_d \vee s^{1/\alpha})^{b_1}
	\int_{z \in \R^d_+, \, x_d \vee s^{1/\alpha}<|x-z|<2}   \L^{b_2}\Big(\frac{  |x-z|}{x_d \vee s^{1/\alpha} } \Big)  \frac{f_{\gamma,\eta_1,\eta_2,k,l}(2|x-z|)	}{|x-z|^{d+\alpha+b_1}}  dz\\
	&\le c_9 s(x_d \vee s^{1/\alpha})^{b_1} f_{\gamma,\eta_1,\eta_2,k,l}(4) \L^{b_2}\Big(\frac{ 2}{x_d \vee s^{1/\alpha} } \Big) 
	\int_{z \in \R^d_+, \,|x-z|<2}  \frac{dz	}{|x-z|^{d+\alpha+b_1-\gamma+\eps}}  \\
	&\le c_{10} s(x_d \vee s^{1/\alpha})^{b_1}   f_{\gamma,\eta_1,\eta_2,k,l}(1) \L^{b_2}\Big(\frac{  2}{x_d \vee s^{1/\alpha} } \Big) . 
\end{align*}

(iii) Case $\gamma=\alpha+b_1$: In this case, we see that 
\begin{align*}
	&	s\int_{z \in \R^d_+, \, x_d \vee s^{1/\alpha}<|x-z|<2} \Big(1 \wedge \frac{x_d\vee s^{1/\alpha} }{|x-z|} \Big)^{b_1}  \L^{b_2}\Big(\frac{  |x-z|}{(x_d  \vee s^{1/\alpha}) \wedge |x-z|} \Big)  \frac{f_{\gamma,\eta_1,\eta_2,k,l}(z_d)}{|x-z|^{d+\alpha}} dz \\
	&\le c_{11} s(x_d \vee s^{1/\alpha})^{b_1}
	\int_{z \in \R^d_+, \, x_d \vee s^{1/\alpha}<|x-z|<2}   \L^{b_2}\Big(\frac{  |x-z|}{x_d \vee s^{1/\alpha} } \Big)  \frac{f_{\gamma,\eta_1,\eta_2,k,l}(|x-z|)	}{|x-z|^{d+\alpha+b_1}}  dz\\
&= c_{11} s(x_d \vee s^{1/\alpha})^{b_1}	\int_{z \in \R^d_+, \, x_d \vee s^{1/\alpha}<|x-z|<2}   \L^{b_2}\Big(\frac{  |x-z|}{x_d \vee s^{1/\alpha} } \Big)\\
	&\hspace{5.3cm}\times  \L^{\eta_1} \Big(\frac{ k}{|x-z|} \Big)\L^{\eta_2} \Big(\frac{  |x-z|}{l} \Big)  \frac{dz	}{|x-z|^{d}} \\
	&\le c_{12}s(x_d \vee s^{1/\alpha})^{b_1}\int_{x_d \vee s^{1/\alpha}}^2   \L^{b_2}\Big(\frac{  r}{x_d \vee s^{1/\alpha} } \Big) \L^{\eta_1} \Big(\frac{ k}{r} \Big) \L^{\eta_2} \Big(\frac{r}{l} \Big)  \frac{dr	}{r}.
\end{align*}
The proof is complete.
\end{proof}

\begin{lemma}\label{cal:3}
	Let $b_1,b_2 \ge 0$. For any $0<k\le l<1$, 
	$$
	\int_l^2 \L^{b_1} \Big( \frac{r}{k}\Big) \L^{b_1} \Big( \frac{r}{l}\Big) \L^{b_2} \Big( \frac{1}{r} \Big) \frac{dr}{r} \asymp \L^{b_1} \Big( \frac{1}{k}\Big) \L^{b_1+b_2+1} \Big(\frac{1}{l}\Big),
	$$
	with comparison constants independent of $k$ and $l$.
\end{lemma}
\begin{proof} Note that
\begin{equation}\label{e:log-equiv}
\L\Big( \frac{1}{r} \Big) \asymp \log\Big( \frac{2e}{r} \Big) \asymp \L\Big(\frac{1}{\sqrt r} \Big), \quad 0<r<2.
\end{equation}
Hence, we get 
\begin{align*}
		&\int_l^2 \L^{b_1} \Big( \frac{r}{k}\Big) \L^{b_1} \Big( \frac{r}{l}\Big) \L^{b_2} \Big( \frac{1}{r} \Big) \frac{dr}{r} \ge c_1 \int_{\sqrt l}^{2} \log^{b_1} \Big( \frac{2er}{k}\Big) \log^{b_1} \Big( \frac{2er}{l}\Big) \log^{b_2} \Big( \frac{2e}{r} \Big) \frac{dr}{r}\\
			&\ge c_1\log^{b_1} \Big( \frac{2e}{\sqrt k}\Big) \log^{b_1} \Big( \frac{2e}{\sqrt l}\Big) \int_{\sqrt l}^{2}  \log^{b_2} \Big( \frac{2e}{r} \Big) \frac{dr}{r}\\
			&=  \frac{c_1}{b_2+1}\log^{b_1} \Big( \frac{2e}{\sqrt k}\Big) \log^{b_1} \Big( \frac{2e}{\sqrt l}\Big) \Big( \log^{b_2+1}\Big( \frac{2e}{\sqrt l} \Big) - 1   \Big) 
			\ge c_2 \L^{b_1} \Big( \frac{1}{k}\Big) 
   \L^{b_1+b_2+1} \Big(\frac{1}{l}\Big).
\end{align*}
On the other hand, using \eqref{e:log-equiv} and \eqref{e:slowly-varying-2}, we get
\begin{align*}
		\int_l^2 \L^{b_1} \Big( \frac{r}{k}\Big) \L^{b_1} \Big( \frac{r}{l}\Big) \L^{b_2} \Big( \frac{1}{r} \Big) \frac{dr}{r}&\le c_3\L^{b_1} \Big( \frac{2}{k}\Big) \L^{b_1} \Big( \frac{2}{l}\Big) \int_l^2 \log^{b_2} \Big(\frac{2e}{r} \Big) \frac{dr}{r}\\
		&	\le c_4\L^{b_1} \Big( \frac{1}{k}\Big) 
		\L^{b_1+b_2+1} \Big(\frac{1}{l}\Big).
\end{align*}
\end{proof}

\begin{lemma}\label{cal:green}
Let $\gamma>1$ and $b_1,b_2\ge 0$. For any $a,k,l>0$,
$$
\int_a^\infty t^{-\gamma} \Big(1 \wedge \frac{k}{t} \Big)^{b_1}  \Big(1 \wedge \frac{l}{t} \Big)^{b_2}  dt \asymp  a^{1-\gamma} \Big(1 \wedge \frac{k}{a} \Big)^{b_1}\Big(1 \wedge \frac{l}{a} \Big)^{b_2},
$$
with comparison constants independent of $a,k$ and $l$.
\end{lemma}
\begin{proof} We have
\begin{align*}
	\int_a^\infty t^{-\gamma} \Big(1 \wedge \frac{k}{t} \Big)^{b_1}  \Big(1 \wedge \frac{l}{t} \Big)^{b_2}  dt 
	&\le \frac{a^{1-\gamma}}{\gamma-1} \Big(1 \wedge \frac{k}{a} \Big)^{b_1}\Big(1 \wedge \frac{l}{a} \Big)^{b_2}
\end{align*}
and
\begin{align*}
		\int_a^\infty t^{-\gamma} \Big(1 \wedge \frac{k}{t} \Big)^{b_1}  \Big(1 \wedge \frac{l}{t} \Big)^{b_2}  dt  &\ge  \Big(1 \wedge \frac{k}{2a} \Big)^{b_1}  \Big(1 \wedge \frac{l}{2a} \Big)^{b_2}\int_a^{2a} t^{-\gamma}  dt\\
		 & \ge \frac{(1-2^{1-\gamma})a^{1-\gamma}}{2^{b_1+b_2}(\gamma-1)} \Big(1 \wedge \frac{k}{a} \Big)^{b_1}\Big(1 \wedge \frac{l}{a} \Big)^{b_2}.
\end{align*}
\end{proof}

\bigskip
\noindent
{\bf Acknowledgements:}
We thank the referees for helpful comments and suggestions. A significant part of this work was done while Zoran Vondra\v{c}ek was
Guest Professor of the College of Natural Sciences at Seoul National University within the
Brain Pool Program of the National Research Foundation of Korea (NRF). The hospitality of
Department of Mathematical Sciences and Research Institute of Mathematics, Seoul National
University, is gratefully acknowledged.

\end{document}